\documentclass[opre,nonblindrev]{informs3opt}

\usepackage{endnotes}
\let\footnote=\endnote

\usepackage{natbib}
\bibpunct[, ]{(}{)}{,}{a}{}{,}
\def\bibfont{\small}

\TheoremsNumberedThrough 
\ECRepeatTheorems

\EquationsNumberedThrough

\usepackage{mathtools}

\usepackage{graphicx}
\usepackage{subfigure}
\usepackage{multirow,enumerate}

\usepackage[hidelinks,bookmarksnumbered=true]{hyperref}
\usepackage{tikz}
\usepackage{url}
\usepackage{amssymb}
\usepackage{psfrag}
\usepackage{amsmath}
\usepackage{setspace}
\newcommand{\exclude}[1]{}
\usepackage{bm}
\usepackage[flushleft]{threeparttable}
\usepackage{algorithm}
\usepackage{algpseudocode}
\algdef{SE}[DOWHILE]{Do}{doWhile}{\algorithmicdo}[1]{\algorithmicwhile\ #1}
\algnewcommand{\Or}{\textbf{or}}
\algnewcommand{\And}{\textbf{and}}

\usepackage{titlesec}
\usepackage[titletoc,toc,title]{appendix}

\usepackage{thmtools} 
\usepackage{thm-restate}

\usepackage{cleveref}

\def\E{{\mathbb E}}

\def\Pr{{\mathbb{P}}}

\def\Re{\mathbb{R}}
\def\Qe{\mathbb{Q}}
\def\I{\mathbb{I}}

\def\hat{\widehat}

\def \F{\mathcal{F}}

\def \Ze{{\mathbb{Z}}}

\def\F{{\mathcal F}}

\def\Re{{\mathbb R}}

\def\X{{\mathcal X}}

\DeclareMathOperator{\proj}{proj}

\newcommand{\bfx}{\bm{x}}

\newcommand{\bfz}{\bm{z}}
\newcommand{\bfb}{\bm{b}}
\newcommand{\bfy}{\bm{y}}

\usepackage{soul}
\usepackage{xcolor, soul}
\usepackage{comment}
\newcommand{\vect}[1]{\boldsymbol{\bm{#1}}}

\algtext*{EndWhile}
\algtext*{EndIf}
\algdef{SE}[DOWHILE]{Do}{doWhile}{\algorithmicdo}[1]{\algorithmicwhile\ #1}
\algnewcommand\INPUT{\item[\textbf{Input:}]}
\algnewcommand\OUTPUT{\item[\textbf{Output:}]}
\def\F{{\mathcal F}}

\def\KSD{{\mathrm{KSD}}}
\def\WD{{\mathrm{WD}}}
\def\DP{{\mathrm{DP}}}
\def\MC{{\mathrm{MC}}}
\newcommand*{\QEDA}{\hfill\ensuremath{\square}}

\usepackage{cleveref}

\usepackage{graphics}
\usepackage{xcolor}

\usepackage{pgfplots}
\pgfplotsset{width=8cm,compat=1.9}

\usepackage{tikz}
\usetikzlibrary{shapes,decorations,arrows,calc,arrows.meta,fit,positioning}
\tikzset{
state/.style ={ellipse, draw, minimum width = 0.7 cm},
}
\usetikzlibrary{decorations.pathreplacing}

\allowdisplaybreaks

\usepackage[flushleft]{threeparttable}

\begin{document}

\RUNAUTHOR{Qing Ye, Grani A. Hanasusanto, and Weijun Xie}
	
\RUNTITLE{Distributionally Fair Stochastic Optimization using Wasserstein Distance}
	
\TITLE{Distributionally Fair Stochastic Optimization using Wasserstein Distance}
	
\ARTICLEAUTHORS{
\AUTHOR{Qing Ye}
\AFF{H. Milton Stewart School of Industrial and Systems Engineering, Georgia Institute of Technology, Atlanta, GA, USA, \EMAIL{qye40@gatech.edu}} 
\AUTHOR{Grani A. Hanasusanto}
\AFF{Department of Industrial \& Enterprise Systems Engineering, University of Illinois Urbana-Champaign, Urbana, IL, USA, 
\EMAIL{gah@illinois.edu}} 
\AUTHOR{Weijun Xie}
\AFF{H. Milton Stewart School of Industrial and Systems Engineering, Georgia Institute of Technology,  Atlanta, GA, USA, \EMAIL{wxie@gatech.edu}}
} 

\ABSTRACT{
A traditional stochastic program under a finite population typically seeks to optimize efficiency by maximizing the expected profits or minimizing the expected costs, subject to a set of constraints. However, implementing such optimization-based decisions can have varying impacts on individuals, and when assessed using the individuals' utility functions, these impacts may differ substantially across demographic groups delineated by sensitive attributes, such as gender, race, age, and socioeconomic status. As each group comprises multiple individuals, a common remedy is to enforce group fairness, which necessitates the measurement of disparities in the distributions of utilities across different groups. This paper introduces the concept of Distributionally Fair Stochastic Optimization (DFSO) based on the Wasserstein fairness measure. The DFSO aims to minimize distributional disparities among groups, quantified by the Wasserstein distance, while adhering to an acceptable level of inefficiency. Our analysis reveals that: (i) the Wasserstein fairness measure recovers the demographic parity fairness prevalent in binary classification literature; (ii) this measure can approximate the well-known Kolmogorov–Smirnov fairness measure with considerable accuracy; and (iii) despite DFSO's biconvex nature, the epigraph of the Wasserstein fairness measure is generally Mixed-Integer Convex Programming Representable (MICP-R). Additionally, we introduce two distinct lower bounds for the Wasserstein fairness measure: the Jensen bound, applicable to the general Wasserstein fairness measure, and the Gelbrich bound, specific to the type-2 Wasserstein fairness measure. We establish the exactness of the Gelbrich bound and quantify the theoretical difference between the Wasserstein fairness measure and the Gelbrich bound. Lastly, the theoretical underpinnings of the Wasserstein fairness measure enable us to design efficient algorithms to solve DFSO problems. Our numerical studies validate the effectiveness of these algorithms, confirming their practical use in achieving distributional fairness in several societally pertinent real-world stochastic optimization problems. 
}
\maketitle

\section{Introduction}

Optimization empowers decision-making by providing an efficient solution to address complex problems in many domains. Its widespread use has motivated research studies focusing on the societal impact of optimization-based decisions. Since the traditional approach optimizes efficiency relevant to profits or costs, 
the optimization outcomes can have varying impacts across demographic groups delineated by sensitive attributes, including gender, race, age, and socioeconomic status. As each group comprises multiple individuals, enforcing group fairness necessitates the measurement of disparities of probability distributions of individual utilities between different groups. Traditional fairness measures are often based on summary statistics, such as minimum, mean, or deviation, which can be insufficient to quantify distributional disparities since each notion only characterizes a particular aspect of the probability distributions. On the other hand, statistical distance metrics, such as the Wasserstein distance, can be employed to quantify distributional fairness accurately. However, these metrics introduce significant computational challenges, and hence they remain largely unexplored in the field of fair decision-making. This motivates us to study distributional fairness.

\subsection{Setting}
The conventional decision-making problem under uncertainty is to optimize the total expected cost efficiency. Such an optimization problem can be formulated as the stochastic program
\begin{align}
 V^*=\min_{\bfx \in \mathcal X} \E_{\Pr}[Q(\bfx,\tilde{\bm\xi})],\label{eq_sp}
\end{align}
where $\mathcal X\subseteq \Re^n$ specifies a mixed-integer convex representable decision space \citep{lubin2022mixed}, $Q(\cdot,\cdot)$ is a recourse function in stochastic programming or a loss function in machine learning, and $\tilde{\bm\xi}\in \Re^{\kappa}$ are the random problem parameters governed by a probability distribution $\Pr$ with support $\Xi$. The stochastic program \eqref{eq_sp} and its variants with risk aversion and distributional robustness have been a prevailing modeling paradigm for numerous decision-making problems (see the survey paper \citealt{rahimian2019distributionally}).

Many real-life decision-making problems may often involve a sensitive attribute such as gender, race, or age in the random parameters $\tilde{\bm\xi}$,  designated by the component $\tilde{\xi}_{\kappa}\in A$, where the set $A$  denotes a finite collection of possible outcomes in the sensitive attribute (e.g., $A=\{\textrm{male, female}\}$). This sensitive attribute partitions the outcome space into groups. Thus,  by invoking the law of total expectation, we can rewrite the stochastic program \eqref{eq_sp} equivalently as
\begin{align}\label{eq_sp_v2}
 V^*=\min_{\bfx \in \mathcal X} \sum_{a\in A} \Pr(\tilde{\xi}_{\kappa}=a) \E_{\Pr_a}[Q(\bfx,\tilde{\bm\xi}_a)],
\end{align}
where $\Pr_a$ is a shorthand for the conditional distribution of $\tilde{\bm\xi}$ given $\tilde{\xi}_{\kappa}=a$. Observe that the objective function constitutes a weighted sum of conditional expectations, where the weights correspond to the marginal distribution of the sensitive attribute. From this vantage point, the optimal solution $\bm x$ may treat the minority groups unfairly as it emphasizes groups of higher weight. This observation motivates us to study fair stochastic programming.
Since many pertinent decision-making problems with sensitive attributes are concerned with a finite population, we assume that the entire support set $\Xi$ is finite (i.e., $\Xi=\{\bm{\xi}_i\}_{i\in [m]}$), and we assume that each group $a\in A$ consists of $m_a$ individuals represented by the set $C_a$, i.e., $\Pr_{a}\{\tilde{\bm\xi}_a=\bm{\xi}_i\}=1/m_a$ for any $i\in C_a$. Evidently, the following identities hold in view of our assumption: $\Xi=\cup_{a\in A}\{\bm \xi_i\}_{i\in C_a}$ and $C_{a}\cap C_{\bar a}=\emptyset$ for any $a<\bar a\in A$. Under this setting, the stochastic program \eqref{eq_sp_v2} further simplifies to
\begin{align}\label{eq_sp_v3}
 V^*=\min_{\bfx \in \mathcal X} \sum_{a\in A} \frac{m_a}{m} \sum_{i\in C_a}\frac{1}{m_a} Q(\bfx,\bm\xi_i).
\end{align}
Many deterministic optimization problems involving multiple groups of individuals can be viewed as a special case of \eqref{eq_sp_v3} by treating each individual as an equiprobable sample.

To measure fairness, given a decision $\bfx\in \mathcal X$, for an individual realization $\bm\xi_a$ in each $a\in A$, 
we suppose that the function $f(\bfx,\bm \xi_a)$ 
denotes its utility value, which may not be monotonic (see, e.g., \citealt{kliegr2009uta}). Our goal is to match the probability distributions of the random utility values $\{f(\bfx,\tilde{\bm\xi}_a)\}_{a\in A}$ among different groups to attain fairness, where we quantify the utility distributional disparities 
using a statistical distance metric. 
Since the random utilities may have different support sets, we employ the Wasserstein distance and propose the following Distributionally Fair Stochastic Optimization (DFSO):
\begin{align}
v^*(q)=\min_{\bfx \in \mathcal X}\left\{\WD_q^q(\bfx):=\max_{a<\bar a\in A}W_q^q\left(\Pr_{f(\bfx,\tilde{\bm\xi}_a)},\Pr_{f(\bfx,\tilde{\bm\xi}_{\bar a})}\right): \E_{\Pr}[Q(\bfx,\tilde{\bm\xi})]\leq V^*+\epsilon|V^*|\right\}.
 \tag{DFSO} \label{eq_sp_fair}
\end{align}
Here, the objective function represents the $q$th power of type-$q$ Wasserstein fairness measure. Particularly, the type-$q$ Wasserstein distance $W_q(\cdot,\cdot)$ is defined as 
\[W_{q}(\Pr_1,\Pr_2)=\inf_{\Qe} \left\{\sqrt[^q]{\int_{\Xi\times\Xi}\left\|\bm{\zeta}_{1}-\bm{\zeta}_{2}\right\|^{q} \Qe(d\bm{\zeta}_{1},d\bm{\zeta}_{2})}: 
\begin{array}{c}\displaystyle
\text{$\Qe$ is a joint distribution of $\tilde{\bm{\zeta}}_{1}$ and $\tilde{\bm{\zeta}}_{2}$}\\\displaystyle 
\text{with marginals $\Pr_1$ and $\Pr_2$, respectively}
\end{array} \right\},\]
where $\|\cdot\|$ is a norm and $q\in[1,\infty]$. In \ref{eq_sp_fair}, the goal is to minimize the maximum distributional disparities of utilities quantified by the Wasserstein distance among all pairs of groups (i.e., the Wasserstein fairness) while maintaining the cost efficiency around a near-optimal region, where $\epsilon\geq 0$ denotes the inefficiency level prescribed by the decision-maker. 
In practice, the utility function $f(\bfx,{\bm\xi})$ can be quite general. If it is equal to the recourse function $Q(\bfx,{\bm\xi})$, then the decision-maker, in this case, tries to achieve the distributional fairness of random cost among different groups.

\subsection{Literature Review}
Optimization has served as an essential tool in decision-making over the past decades. Throughout the years, the issues of fairness in optimization have been recognized and studied in the fields of resource allocation, facility location, and communication networks \citep{ogryczak2014fair, karsu2015inequity}. The commonly adopted definitions of fairness pertain to the utilities of all individuals in the population, e.g., max-min fairness, proportional fairness, and alpha fairness, or to some particular characteristics of the distribution of the utilities, e.g., spread, deviation, Jain's index, and Gini coefficient. 
Contrary to traditional definitions that consider the entire population, this paper concentrates on fairness among different groups of individuals. These fairness measures at the population level can be simply generalized to the group level by applying them to each group instead of the entire population. 
For example, \cite{samorani2022overbooked} studied the max-min fairness at the group level. They addressed the racial disparity in medical appointment scheduling by minimizing the maximum waiting time among the racial groups. \cite{cohen2022price} discussed price discrimination against protected groups and attempted to enforce nearly equal prices for different groups. 
\cite{patel2020group} considered group fairness for the knapsack problem when each item belongs to a particular group. They defined three fair knapsack notions, i.e., to bound the number of items from each group, to bound the total weight of items from each group, and to bound the total value of items from each group.
Since the traditional fairness measures are often based on summary statistics, they might be inadequate for quantifying group 
disparities in a comprehensive way. 
Our distributional fairness notion overcomes this limitation by using the Wasserstein distance to quantify the distributional disparities among different groups. The Wasserstein distance has also been used in a variety of optimization problems such as Wasserstein distributional robust optimization \citep{mohajerin2018data,blanchet2019quantifying,gao2023distributionally,hanasusanto2018conic,chen2022data,xie2021distributionally}.

Recent studies in the growing field of fair machine learning have proposed various methods for a number of tasks \citep{caton2020fairness}.
The majority of the literature has focused on group fairness, which seeks to treat different groups equally. 
Group fairness in binary classification has been extensively studied \citep{kamishima2012fairness,feldman2015certifying,barocas2016big,hardt2016equality,zafar2017fairnessb,donini2018empirical,aghaei2019learning,kallus2022assessing,taskesen2020distributionally,ye2020unbiased,wang2021wasserstein,lowy2021stochastic}. 
However, the number of works on group fairness in regression with continuous outcomes is rather limited. 
\cite{berk2017convex} introduced a family of convex fairness regularizers such that each group should have similar predicted outcomes weighted by the nearness of the true outcomes on average. \cite{agarwal2019fair, chzhen2020fair, rychener2022metrizing} used the Kolmogorov–Smirnov distance to achieve 
demographic parity. Additionally, \cite{rychener2022metrizing} summarized the common integral probability metrics for quantifying fairness, including the Kolmogorov–Smirnov distance and the Wasserstein distance. To achieve fairness, \cite{agarwal2019fair} designed a reduction-based algorithm, while \cite{chzhen2020fair} developed a post-processing algorithm for fair regression. \cite{rychener2022metrizing} proposed to solve fair regression via a stochastic gradient descent algorithm. 
According to the definition in fair machine learning literature, the demographic parity-based fairness notion ensures the probability distribution of outcomes is independent of the sensitive attribute groups. Our distributional fairness notion coincides with demographic parity when applied to machine learning problems. 
Furthermore, the proposed \ref{eq_sp_fair} formulation is a general stochastic optimization problem where fairness is integrated with efficiency. Thus, it provides flexibility to model various decision-making problems, including classification with binary utilities and regression with continuous utilities. More importantly, different from existing results in the literature, we thoroughly investigate the optimization properties of the Wasserstein fairness measure and exploit them to systematically design efficient solution algorithms with provable guarantees.

\subsection{Summary of Contributions}
The main contributions of this paper are summarized as follows:
\begin{itemize}
\item From a fresh scope, this paper establishes the fundamental result that the Wasserstein fairness measure is essentially equivalent to matching the probability distributions of distinct groups comonotonically and computing the distance of the comonotonic distributions. Using this equivalence, we show that the Wasserstein fairness measure recovers the well-known demographic parity fairness from the binary classification literature, and we reveal that the Wasserstein fairness measure is relatively close to the Kolmogorov–Smirnov one.
\item  We prove that the \ref{eq_sp_fair} under the Wasserstein fairness measure, in general, is NP-hard. However, different from other biconvex programs, we show that the epigraph of the Wasserstein fairness measure is, in general, Mixed-Integer Convex Programming Representable (MICP-R), and we provide four different representations. These are the first known MICP-R results for Wasserstein distance-based distributional fairness models. 
\item We derive two different lower bounds for the Wasserstein fairness measure: the Jensen bound for the general Wasserstein fairness measure and the Gelbrich bound for the type-2 Wasserstein fairness measure. We prove a broader condition than the well-known elliptical distributions under which the Gelbrich bound is asymptotically tight, and we provide a theoretical gap between the Wasserstein fairness measure and the Gelbrich bound. We also prove that computing the Gelbrich bound is NP-hard. 
\item Inspired by the theoretical properties of the Wasserstein fairness measure, we design effective solutions algorithms to solve the \ref{eq_sp_fair} to near-optimality. Our numerical study confirms the effectiveness of the proposed algorithms.
\end{itemize}

The remainder of the paper is organized as follows. Section \ref{sec_poperties} presents properties of the Wasserstein fairness measure. Section \ref{sec_micp_preliminaries} formalizes definitions  
and develops two exact mixed-integer convex programming representations of the epigraph of the Wasserstein fairness measure. Section \ref{sec_lower_bounds} studies two lower bounds of the Wasserstein fairness measure. Section \ref{sec_numerical_study} reports the numerical study, and Section \ref{sec_conclusion} concludes the paper. Proofs and additional results are relegated to the appendix.

\noindent\textbf{Notation.} 
Bold lowercase letters (e.g., $\vect{x}$) denote vectors, bold uppercase letters (e.g., $\vect{Z}$) denote matrices, and the corresponding regular letters (e.g., $x_i, Z_{ij}$) denote their components. 
For any $n\in\Ze_+$, we let $[n]:=\{1,2,\ldots,n\}$ and use $\Re_+^n:=\{\vect{x}\in \Re^n:x_i\geq0, \forall i\in [n]\}$. For any $n_1<n_2\in\Ze_+$, we let $[n_1, n_2]:=\{n_1,n_1+1,\ldots,n_2\}$.
For a set $A$, we let  $a<\bar a\in A$ denote $a,\bar{a}\in A$ such that $a<\bar{a}$. 
The indicator function $\I(B)$ takes value $1$ if $B$ is true and $0$ otherwise. Additional notation will be introduced as needed.

\section{Properties of the Wasserstein Fairness Measure}\label{sec_poperties}

This section presents various notable properties of the Wasserstein fairness measure. To begin with, let us define the cumulative distribution functions of the random functions $\{f(\bfx,\tilde{\bm\xi}_a)\}_{a\in A}$ as $F_a(t\mid\bfx)=\Pr_{a}\{f(\bfx,\tilde{\bm\xi}_a)\leq t\}$ 
for all $a\in A$. Correspondingly, we define the inverse distribution functions $F^{-1}_{a}(y\mid\bfx)=\inf \{t: F_{a}(t\mid\bfx)\geq y\}$ for all $a\in A$. 

\subsection{Comonotonicity and Complexity}

This subsection investigates the comonotonicity property of the Wasserstein fairness measure and the complexity of \ref{eq_sp_fair}, which motivate us to develop strong mixed-integer convex programming formulations for \ref{eq_sp_fair}.

One property of the Wasserstein fairness measure in \ref{eq_sp_fair} is that it can be simplified as the integral of the difference of inverse cumulative distributions. 
\begin{lemma}[Proposition 2.17 in \cite{santambrogio2015optimal}]\label{lem_icf_WD}
For any $a<\bar a\in A$ and a fixed decision $\bfx$, the Wasserstein distance $W_q\left(\Pr_{f(\bfx,\tilde{\bm\xi}_a)},\Pr_{f(\bfx,\tilde{\bm\xi}_{\bar a})}\right)$ can be expressed as
\begin{align}\label{eq_lem_inv_cdf_WD}
W_q\left(\Pr_{f(\bfx,\tilde{\bm\xi}_a)},\Pr_{f(\bfx,\tilde{\bm\xi}_{\bar a})}\right)= &\sqrt[^q]{\int_{0}^{1}\left|F^{-1}_{a}(y\mid\bfx)-F_{\bar a}^{-1}(y\mid\bfx)\right|^{q}dy},
\end{align}
where $F^{-1}_a$ is the inverse distribution function of the random function $f(\bfx,\tilde{\bm\xi}_a)$ for each $a\in A$. 
When $q=1$, the type-1 Wasserstein distance $W_1\left(\Pr_{f(\bfx,\tilde{\bm\xi}_a)},\Pr_{f(\bfx,\tilde{\bm\xi}_{\bar a})}\right)$ coincides with the $L_1$ distance between the cumulative distribution functions
\begin{align}W_1\left(\Pr_{f(\bfx,\tilde{\bm\xi}_a)},\Pr_{f(\bfx,\tilde{\bm\xi}_{\bar a})}\right)= &\int_{\Re} \left|F_{a}(t\mid\bfx)-F_{\bar a}(t\mid\bfx)\right|dt,\label{eq_lem_inv_cdf_WD2}
\end{align}
 where $F_a$ is the cumulative distribution function of the random function $f(\bfx,\tilde{\bm\xi}_a)$ for each $a\in A$. 
\end{lemma}

\Cref{lem_icf_WD} shows that the Wasserstein fairness measure $\WD_q(\bfx):=\max_{a<\bar a\in A}W_q(\Pr_{f(\bfx,\tilde{\bm\xi}_a)},\Pr_{f(\bfx,\tilde{\bm\xi}_{\bar a})})$ can be viewed as the largest $L_q$-norm of the difference between inverse distribution functions. 
Hence, in \ref{eq_sp_fair}, minimizing $\WD_q(\bfx)$ implies attempting to match the distributions of utilities between any two 
groups $a<\bar a\in A$. 
Remarkably, as established in the existing literature, type-1 Wasserstein fairness measure $\WD_1(\bfx)$ is equivalent to the maximum $L_1$ distance between the cumulative distribution functions. For each pair of groups $a<\bar{a}\in A$, under our assumption of discrete distributions, the integrals in \eqref{eq_lem_inv_cdf_WD} and \eqref{eq_lem_inv_cdf_WD2} can be simplified to be summations. These properties motivate us to study the exact MICP formulations of the Wasserstein fairness measure.

Another interesting byproduct of \Cref{lem_icf_WD} is that when achieving the infimum of the Wasserstein distance, the two distributions must be aligned comonotonically, where the comonotonicity of two random variables is formally defined as follows.
\begin{definition}\label{def_comontonic}
A pair of random variables $X = (X_1,X_2)$ is comonotonic if and only if it can be represented as
$\displaystyle (X_{1},X_{2})\stackrel{\text{d}}{=}(F_{X_{1}}^{-1}(U),F_{X_{2}}^{-1}(U))$, where $U$ is the standard uniform random variable, and $F_{X_{1}}^{-1}, F_{X_{2}}^{-1}$ are the inverse distribution functions of $X_1, X_2$.
\end{definition}
This 
gives rise to an interesting result for the following Wasserstein fairness measure. 
\begin{restatable}{proposition}{propcomono}\label{prop_comono}
For a given decision $\bfx \in\mathcal X$, when computing the Wasserstein fairness measure 
in \ref{eq_sp_fair}, the optimal joint distribution is comonotonic for any pair $a<\bar a\in A$.
\end{restatable}
 \noindent\textit{Proof. } See Appendix~\ref{proof_prop_comono}.\QEDA

\Cref{prop_comono} shows that the Wasserstein fairness measure, in fact, aligns the two distinct groups' utility function values comonotonically, computes the $L_q$ norm of the difference of their inverse distribution functions, and then takes the maximum value among all the pairs of groups. It helps us study the new exactness conditions of the well-known lower bound (i.e., the Gelbrich bound) of the type-2 Wasserstein fairness measures. This result also motivates us to study its relation with another popular distributional fairness notion: the Kolmogorov–Smirnov fairness measure.

We conclude the subsection by proving the NP-hardness of  \ref{eq_sp_fair} via a reduction from the well-known chance-constrained stochastic program \citep{charnes1959chance,ahmed2018relaxations}.
\begin{restatable}{theorem}{thmnpharddfso}\label{thm_np_hard_dfso}
Solving \ref{eq_sp_fair} is, in general, strongly NP-hard, even when $\mathcal X$ is a polytope, $\epsilon=\infty$, and $f(\bfx,\bm\xi)$ is a linear function.
\end{restatable}
\noindent\textit{Proof. } See Appendix~\ref{proof_thm_np_hard_dfso}.\QEDA

\subsection{Recovering the Demographic Parity Fairness Measure of Binary Outcomes}

Demographic parity of binary outcomes, defined as $\DP(\bfx)$, requires the probability of beneficial or detrimental outcomes to be independent of the sensitive attribute. In the following, we show that the proposed $\WD_q(\bfx)$ recovers $\DP(\bfx)$ if the utility function is Bernoulli. Let us first define $\DP(\bfx)$.
\begin{definition}
Suppose that $\Pr\{f(\bfx,\tilde{\bm\xi})\in\{0,1\}\}=1$. The binary demographic parity fairness measure is defined as 
$$\DP(\bfx)=\max_{a<\bar{a} \in A} \left|\Pr\{f(\bfx,\tilde{\bm\xi}_{a})=0\}-\Pr\{f(\bfx,\tilde{\bm\xi}_{\bar a})=0\} \right|=\max_{a<\bar{a} \in A} \left|\Pr_a\{f(\bfx,\tilde{\bm\xi}_{a})=1\}-\Pr_{\bar{a}}\{f(\bfx,\tilde{\bm\xi}_{\bar a})=1\} \right|.$$
\end{definition}
We next show that the Wasserstein fairness measure $\WD_q(\bfx)$ is equivalent to $\DP(\bfx)$ in view of \Cref{lem_icf_WD}.
\begin{restatable}{proposition}{propbernoulli}\label{prop_bernoulli}
For a Bernoulli utility function $f(\bfx,\tilde{\bm\xi})\in\{0,1\}$, $\WD_q(\bfx)$ is equivalent to $\DP(\bfx)$.
\end{restatable}
 \noindent\textit{Proof. }See Appendix \ref{proof_prop_bernoulli}.\QEDA

The result in \Cref{prop_bernoulli} reveals that the proposed Wasserstein fairness measure constitutes a generalization of the binary demographic parity fairness measure.

\subsection{Comparison with the Kolmogorov–Smirnov Fairness Measure
}
Instead of the sum of differences, we can use the supremum of differences to measure the demographic parity fairness as defined below. 
\begin{definition}[Kolmogorov–Smirnov Fairness Measure, \citealt{agarwal2019fair}]
The distributional fairness of a decision $\bfx$ can be measured using the Kolmogorov–Smirnov distance: 
\begin{align}\label{dfn_KSD}
\KSD(\bfx)=\max_{a<\bar a\in A} &\sup_{t} \left|F_{a}(t\mid\bfx)-F_{\bar a}(t\mid\bfx)\right|.
\end{align}
\end{definition}
The Kolmogorov–Smirnov distance measures the largest difference of cumulative distribution functions between any two distinct groups. To compute $\KSD(\bfx)$, one needs to discretize $t$, which is easily done in view of our assumption of finite populations. 
Specifically, the assumption implies that the cumulative distributions $\{F_{a}(y\mid\bfx)\}_{a\in A}$ and their inverse counterparts $\{F_{a}^{-1}(y\mid\bfx)\}_{a\in A}$ are of finitely many values, defined formally as follows. 
\begin{definition}[Breaking Points]\label{def_break_point}
For any $\bfx \in \mathcal X$ and any pair  $a<\bar a\in A$, the breaking points of $F_{a}^{-1}(y\mid\bfx)-F_{\bar a}^{-1}(y\mid\bfx)$ are denoted by $\bfb_{a\bar a}(\bfx)=(b_{ja\bar a}(\bfx))_{j\in{J_{a\bar a}}}$ with index set $J_{a\bar a}:=\{1,2,\cdots, |J_{a\bar a}|\}$. We further define the widths $w_{ja\bar a}(\bfx)=b_{ja\bar a}(\bfx)-b_{(j-1)a\bar a}(\bfx)$ for $j \in J_{a\bar a}\setminus \{1\}$ and calculate the largest width as $\eta(\bfx)=\max_{a<\bar a\in A, j\in J_{a\bar a}\setminus \{1\}}w_{ja\bar a}(\bfx)$.
\end{definition}

Based on \Cref{def_break_point}, we propose the following lower and upper bounds on $\KSD(\bfx)$ in terms of $\WD_q(\bfx)$, which shows that the two measures are close to each other within a constant factor. 
The key idea is to recast the type-$q$ Wasserstein fairness measure as
\begin{align*}
\WD_q(\bfx)=&\max_{a<\bar a\in A} \sqrt[^q]{\int_{0}^{1}\left|F_{a}^{-1}(y\mid\bfx)-F_{\bar a}^{-1}(y\mid\bfx)\right|^{q}dy}\\
=&\max_{a<\bar a\in A} \sqrt[^q]{\sum_{j\in J_{a\bar a}\setminus \{1\}} w_{ja\bar a}(\bfx) \left|F_{a}^{-1}\left(b_{ja\bar a}(\bfx)\mid\bfx\right)-F_{\bar a}^{-1}\left(b_{ja\bar a}(\bfx)\mid\bfx\right)\right|^{q}},
\end{align*}
in the spirit of \Cref{def_break_point}. Then we bound the difference between $\WD_1(\bfx)$ and $\KSD(\bfx)$. Next, we use the relationship between $\WD_1(\bfx)$ and $\WD_q(\bfx)$ to finally bound $\WD_q(\bfx)$ and $\KSD(\bfx)$. 

\begin{restatable}{proposition}{propboundksdWF}
\label{prop_bound_ksd_WF}
For any feasible $\bfx\in \mathcal X$ and $q\in[1,\infty]$, the following inequalities hold:
\[\frac{1}{\max_{a<\bar a\in A}\eta(\bfx)^{\frac{1-q}{q}}(t_{2a\bar a}(\bfx)-t_{1a\bar a}(\bfx))}\WD_q(\bfx) \leq \KSD(\bfx) \leq \frac{1}{\min_{a<\bar a\in A}\mu(\Delta_{a\bar a}(\bfx))}\WD_q(\bfx).\]
Here, 
$t_{1a\bar a}(\bfx)=\min\{\min_{t}\{t: F_{a}(t\mid\bfx)>0\}, \min_{t} \{t: F_{\bar{a}}(t\mid\bfx)>0\}\}$, $t_{2a\bar a}(\bfx)=\max\{\sup_{t} \{t: F_{a}(t\mid\bfx)<1\}, \sup_{t} \{t: F_{\bar a}(t\mid\bfx)<1\}\}$, and $\Delta_{a\bar a}(\bfx)=\{\bar{t}: |F_{a}(\bar{t}\mid\bfx)-F_{\bar a}(\bar{t}\mid\bfx)|=\sup_{t} |F_{a}(t\mid\bfx)-F_{\bar a}(t\mid\bfx)|\}$ with its 
Lebesgue measure $\mu(\Delta_{a\bar a}(\bfx))$.
\end{restatable}
\noindent\textit{Proof. } See Appendix~\ref{proof_prop_bound_ksd_WF}. \QEDA

\Cref{prop_bound_ksd_WF} theoretically establishes that the Wasserstein and Kolmogorov–Smirnov fairness measures are rather similar to each other. The bounds can be independent of the decision variables $\bfx$ by finding the least-favorable coefficients. It is worth mentioning that the existing literature (see, e.g., \citealt{ross2011fundamentals}) only bound Kolmogorov–Smirnov fairness measure from above by type-$1$ Wasserstein fairness measure when the underlying random variables are continuous. In our following derivation, the Wasserstein fairness measure shows amenable optimization properties. Our numerical study demonstrates the advantage of the proposed methods for the Wasserstein fairness measure compared to the existing ones for the Kolmogorov–Smirnov fairness measure.

\section{Mixed-Integer Convex Programming Formulations of \ref{eq_sp_fair}}\label{sec_micp_preliminaries}
This section focuses on deriving exact Mixed-Integer Convex Programming (MICP) 
formulations of \ref{eq_sp_fair}. To begin with, we observe that
under the discrete-distribution assumption, using epigraphical variable $\nu$, the proposed \ref{eq_sp_fair} can be formulated as the mathematical program
\begin{subequations}\label{eq_sp_fair_simple}
\begin{align}
 v^*(q)=\min_{(\bfx ,\nu)\in \F_q}\quad&\nu,\\
 \text{s.t.} \quad&\sum_{i\in [m]}\frac{1}{m}Q(\bfx,\bm{\xi}_i)\leq V^*+\epsilon|V^*|,\label{eq_tol}
\end{align}
\end{subequations}
where we introduce the set $\F_{q}$ to denote the  epigraph of the  Wasserstein fairness measure, as follows:
\begin{equation}\label{eq_F_q}
\F_{q}=\left\{(\bfx,\nu)\in\mathcal X\times \Re_+:W_q^q\left(\Pr_{f(\bfx,\tilde{\bm\xi}_a)},\Pr_{f(\bfx,\tilde{\bm\xi}_{\bar a})}\right) \leq \nu,\forall a<\bar a\in A \right\}.
\end{equation}
For the formulations, we will also utilize the following set that corresponds to the graph of the function $f(\cdot,\bm{\xi}_i)$ for each realization $\bm{\xi}_i\in \Xi$:
\[X_i=\left\{(\bfx,\bar w_i)\in \mathcal X\times \Re: f(\bfx,\bm{\xi}_i)=\bar w_i\right\}.\]

\Cref{sec_micpr_x_i} discusses the concept of MICP representability and presents MICP formulations for the graphs of utility functions $\{X_i\}_{i\in [m]}$. \Cref{sec_M3} and \Cref{sec_M4} explore two different ways of representing the epigraph of Wasserstein fairness measure $\F_q$ in \eqref{eq_F_q}. The first formulation uses \Cref{lem_icf_WD} to represent quantiles using mixed-integer programming formulations. The second formulation is a variation of the first one, using aggregate rather than individual quantiles. We have two additional formulations presented in Appendix \ref{sec_reform_add}, where the \hyperref[WD_q_trans_micpr]{Discretized Formulation} (see Appendix \ref{sec_M1}) is based on the discretization of the transportation decisions by observing that the inflated transportation decision variables can be restricted to integers, and the \hyperref[model2_WD1]{Complementary Formulation} (see Appendix \ref{sec_M2}) is to recast the set $\F_q$ using linear programming with complementary slackness constraints and linearize the complementary slackness constraints. Besides, we derive an equivalent MICP formulation for the Kolmogorov–Smirnov fairness measure $\KSD(\bfx)$, which can be found in Appendix \ref{sec_micpr_KSD}. 

\subsection{Mixed-Integer Convex Programming Representability}\label{sec_micpr_x_i} 
To begin with, we introduce the notion of MICP representability and develop formulations for various families of utility functions, depending on whether the sets $\{X_i\}_{i\in [m]}$ are MICP representable (MICP-R) or not MICP-R. 
The MICP-R sets are defined as follows.
\begin{definition}[Theorem 4.1 in \citealt{lubin2022mixed}]
A set $S\subseteq\Re^n$ is MICP-R if and only if there exists $d,p\in\Ze_{+}$, a convex set $C\subseteq\Re^d$, and a closed convex family $(B_z)_{z\in C}\subseteq \Re^{n+p}$ such that $S=\bigcup_{z\in C\cap\Ze^d}\proj_{x}(B_z)$. 
\end{definition}
In addition, \cite{lubin2022mixed} also provided the following sufficient condition for not MICP-R. 
\begin{lemma}[Lemma 4.1 in \citealt{lubin2022mixed}]
\label{lem_not-micpr}
A set $S\subseteq\Re^n$ is not MICP-R if there exists $R\subseteq S, |R|=\infty$ such that $(\bfx+\bfx')/2\notin S$ for all $\bfx,\bfx'\in R, \bfx\neq\bfx'$.
\end{lemma}
Our MICP-R results rely on the McCormick representation. 
\begin{definition}[McCormick Representation, \citealt{mccormick1976computability}]
Consider a bilinear set $\{(\psi,\kappa,\nu)\in\Re\times\{\kappa^L, \kappa^U\}\times[\nu^L, \nu^U]: \psi=\kappa \nu\}$ with given lower bounds $\kappa^L, \nu^L$ and upper bounds $\kappa^U, \nu^U$. 
Its McCormick representation is
\begin{equation}\nonumber
\MC(\kappa^L, \kappa^U, \nu^L, \nu^U) = \left\{(\psi,\kappa,\nu): \begin{aligned}
& \psi\in\Re, \kappa\in\{\kappa^L, \kappa^U\}, \nu^L\leq \nu\leq \nu^U,\\
& \psi\geq \kappa^L \nu+\kappa \nu^L-\kappa^L \nu^L, \psi\geq \kappa^U \nu+\kappa \nu^U-\kappa^U \nu^U,\\
& \psi\leq \kappa^U \nu+\kappa \nu^L-\kappa^U \nu^L, \psi\leq \kappa \nu^U+\kappa^L \nu-\kappa^L \nu^U
\end{aligned}\right\}.
\end{equation}
\end{definition}

Next, we discuss three special cases when the sets $\{X_i\}_{i\in [m]}$ are MICP-R. 
\begin{proposition}
Suppose that $f(\bfx,\bm\xi)=\bm\xi^\top \bm{r}(\bfx)+ s(\bfx)$, where $\bm{r}(\bfx)$ and $s(\bfx)$ are linear functions. Then the sets $\{X_i\}_{i\in [m]}$ are MICP-R.
\end{proposition}
 \noindent\textit{Proof. }
 We have 
\begin{equation}\nonumber
\begin{aligned}
X_{i}&=\left\{(\bfx, \bar{w}_{i})\in \mathcal X\times\Re: \bm\xi_i^\top \bm{r}(\bfx)+ s(\bfx)= \bar{w}_{i}\right\},
\end{aligned}
\end{equation}
which is an MICP-R set. 
\QEDA

\begin{proposition}
Suppose that $f(\bfx, {\bm\xi})=\max_{\tau\in T}\left\{\bm\xi^\top \bm{r}_\tau(\bfx)+ s_\tau(\bfx)\right\}$, where $\{\bm{r}_\tau(\bfx)\}_{\tau \in T}$ and $\{s_\tau(\bfx)\}_{\tau\in T}$ are linear functions. Then the sets $\{X_i\}_{i\in [m]}$ are MICP-R.
\end{proposition}
 \noindent\textit{Proof. }
Suppose that $M_{i}\geq \max_{\bm{x}\in \mathcal X,\eqref{eq_tol}} \left|\max_{\tau\in T}\left\{\bm\xi_i^\top \bm{r}_\tau(\bfx)+ s_\tau(\bfx)\right\}\right|$ for each $i\in [m]$. Then, we have
\begin{equation}\nonumber
\begin{aligned}
X_{i}&=\left\{(\bfx, \bar{w}_{i})\in \mathcal X\times\Re: \max_{\tau\in T}\left\{\bm\xi_i^\top \bm{r}_\tau(\bfx)+ s_\tau(\bfx)\right\}=\bar{w}_{i}\right\}\\
&=\left\{(\bfx, \bar{w}_{i})\in \mathcal X\times\Re: \begin{aligned}
&\bar{w}_{i}\geq \bm\xi_i^\top \bm{r}_\tau(\bfx)+ s_\tau(\bfx), \forall \tau\in T, \\
&\bar{w}_{i}\leq \bm\xi_i^\top \bm{r}_\tau(\bfx)+ s_\tau(\bfx)+M_{i}(1-z_{i\tau}), \forall \tau\in T, \\ 
&\sum_{\tau\in T} z_{i\tau} = 1, z_{i\tau}\in\{0,1\}, \forall \tau\in T
\end{aligned}\right\},
\end{aligned}
\end{equation}
which is an MICP-R set. 
\QEDA

\begin{proposition}
Suppose that $f(\bfx, {\bm\xi})=\min_{\tau\in T}\left\{\bm\xi^\top \bm{r}_\tau(\bfx)+ s_\tau(\bfx)\right\}$, where $\{\bm{r}_\tau(\bfx)\}_{\tau \in T}$ and $\{s_\tau(\bfx)\}_{\tau\in T}$ are linear functions. Then the sets $\{X_i\}_{i\in [m]}$ are MICP-R.
\end{proposition}
 \noindent\textit{Proof. }
Recall that $M_{i}\geq \max_{\bm{x}\in \mathcal X,\eqref{eq_tol}} \left|\min_{\tau\in T}\left\{\bm\xi_i^\top \bm{r}_\tau(\bfx)+ s_\tau(\bfx)\right\}\right|$ for each $i\in [m]$. Thus, we have
\begin{equation}\nonumber
\begin{aligned}
X_{i}&=\left\{(\bfx, \bar{w}_{i})\in \mathcal X\times\Re: \min_{\tau\in T}\left\{\bm\xi_i^\top \bm{r}_\tau(\bfx)+ s_\tau(\bfx)\right\}=\bar{w}_{i}\right\}\\
&=\left\{(\bfx, \bar{w}_{i})\in \mathcal X\times\Re: \begin{aligned}
&\bar{w}_{i}\leq \bm\xi_i^\top \bm{r}_\tau(\bfx)+ s_\tau(\bfx), \forall \tau\in T, \\
&\bar{w}_{i}\geq \bm\xi_i^\top \bm{r}_\tau(\bfx)+ s_\tau(\bfx)-M_{i}(1-z_{i\tau}), \forall \tau\in T, \\ 
&\sum_{\tau\in T} z_{i\tau} = 1, z_{i\tau}\in\{0,1\}, \forall \tau\in T
\end{aligned}\right\},
\end{aligned}
\end{equation}
which is an MICP-R set. 
\QEDA

When the utility functions are exponential or logarithmic, their corresponding $\F_q$ sets are typically not MICP-R according to \Cref{lem_not-micpr}. Hence we propose to approximate them using piecewise linear functions (see Appendix~\ref{sec_not_micpr}).

\subsection{Quantile Formulation}\label{sec_M3}
In this subsection, we propose a quantile-based formulation to represent the set $\F_q$ motivated by \Cref{lem_icf_WD}. That is, we first equivalently rewrite set $\F_q$ as
\begin{equation}\label{eq_F_q2}
\F_{q}=\left\{(\bfx,\nu)\in\mathcal X\times \Re_+:\int_{0}^{1}\left|F_{a}^{-1}(y\mid\bfx)-F_{\bar a}^{-1}(y\mid\bfx)\right|^{q}dy \leq \nu,\forall a<\bar a\in A \right\}.
\end{equation}
Since all the random parameters have finite support, let us sort the distinct elements of the set 
\[\{0\}\cup \left\{\frac{i}{m_a}\right\}_{i\in [m_a]} \cup\left\{\frac{i}{m_{\bar a}}\right\}_{i\in [m_{\bar a}]}:=\left\{\hat{b}_{ia\bar a}\right\}_{i\in[\hat{m}_{a\bar a}]}\]
in the ascending order as 
$0:=\hat{b}_{1a\bar a}<\cdots< \hat{b}_{(\hat{m}_{a\bar a})a\bar a}:=1$ for each $a<\bar a\in A$. 
Observe that in the equation \eqref{eq_lem_inv_cdf_WD}, the value $F^{-1}_{a}(y\mid\bfx)-F^{-1}_{\bar a}(y\mid\bfx)$ is a constant whenever $y\in (\hat{b}_{ia\bar a},\hat{b}_{(i+1)a\bar a}]$ for $i\in [\hat{m}_{a\bar a}-1]$. Thus, the set $\{\hat{b}_{ia\bar a}\}_{i\in[\hat{m}_{a\bar a}]}$ helps simplify the Wasserstein fairness measure as 
\begin{equation}\label{eq_lem_inv_cdf_WD1}
\begin{aligned}
\WD_q^q(\bfx)=&\max_{a<\bar a\in A} {\sum_{i\in [\hat{m}_{a\bar a}-1]}\int_{\hat{b}_{ia\bar a}}^{\hat{b}_{(i+1)a\bar a}}\left|F^{-1}_{a}(y\mid\bfx)-F^{-1}_{\bar a}(y\mid\bfx)\right|^{q}dy},\\
=&\max_{a<\bar a\in A} {\sum_{i\in [\hat{m}_{a\bar a}-1]}(\hat{b}_{(i+1)a\bar a}-\hat{b}_{ia\bar a})\left|F^{-1}_{a}(\hat{b}_{(i+1)a\bar a}\mid\bfx)-F^{-1}_{\bar a}(\hat{b}_{(i+1)a\bar a}\mid\bfx)\right|^{q}}.
\end{aligned}
\end{equation}

Next, we define the quantile set 
$\Omega_a(k)=\{(\bfx,t_{ka})\in \mathcal X\times \Re: F^{-1}_{a}(k/m_a\mid\bfx)=t_{ka}\}$ for each $k\in [m_a]$ and $a\in A$. Using the graph representation  $(\bfx,\bar{w}_i)\in X_i$ for each $i\in[m]$, we propose the following equivalent formulation of the quantile set $\Omega_a(k)$.
\begin{restatable}{proposition}{propquantset}\label{prop_quant_set}
Suppose that $M_{i}\geq \max_{\bm{x}\in \mathcal X,\eqref{eq_tol}} |f(\bfx,{\bm\xi}_i)|$ for each $i\in [m]$. For each $k\in [m_a]$ and $a\in A$, the quantile set $\Omega_a(k)$ is equivalent to
\begin{equation}\label{set_Qa}
\Omega_a(k)=\left\{(\bfx,t_{ka})\in \mathcal X\times \Re: \begin{aligned}
& \pi_{ika}\in\{0,1\}, z_{ika}\in\{0,1\},\pi_{ika}\leq z_{ika},(\bfx,\bar{w}_{i})\in X_i,\forall i\in C_a,\\
& \sum_{i\in C_a}z_{ika}=k, \sum_{i\in C_a}\pi_{ika}=1,t_{ka}=\sum_{i\in C_a}\hat{t}_{ika},\\
& t_{ka}\geq \bar{w}_{i}-(M_i+M_{(k)})(1-z_{ika}), 
t_{ka}\leq \bar{w}_{i}+(M_i+M_{(k)})z_{ika},\\ 
& (\hat{t}_{ika},\pi_{ika},\bar{w}_{i})\in\MC(0,1,-M_i,M_i), \forall i\in C_a
\end{aligned}\right\},
\end{equation}
where $M_{(i)}$ is the $i$th smallest value of the vector $\bm M$.
\end{restatable}

\noindent\textit{Proof. } See Appendix~\ref{proof_prop_quant_set}.
\QEDA

To reformulate the set $\F_{q}$ defined in \eqref{eq_F_q2}, we use the quantile-based representation \eqref{eq_lem_inv_cdf_WD1} of the Wasserstein fairness measure by plugging in the MICP-R quantile sets $\{\Omega_a(k)\}_{k\in [m_a],a\in A}$ defined in \eqref{set_Qa}.
\begin{restatable}{theorem}{thmMthird}\textbf{(Quantile Formulation)}\label{thm_M3}
Suppose that the set $X_{i}=\{(\bfx, \bar{w}_{i})\in \mathcal X\times\Re: f(\bfx,{\bm\xi}_i)=\bar{w}_{i}\}$ is MICP-R and $M_{i}\geq \max_{\bm{x}\in \mathcal X,\eqref{eq_tol}} |f(\bfx,{\bm\xi}_i)|$ for each $i\in [m]$. We further define the quantile set $\Omega_a(k)=\{(\bfx,t_{ka}) \in \mathcal X\times \Re: F^{-1}_{a}(k/m_a\mid\bfx)=t_{ka}\}$, which admits a MICP-R form \eqref{set_Qa}. Then $\F_{q}$ can be represented as 
\begin{equation}\label{model3_WD}
\F_{q}=\left\{(\bm{x},\nu)\in \mathcal X\times \Re_+: \begin{aligned}
&\sum_{i\in [\hat{m}_{a\bar a}-1]}\left(\hat{b}_{(i+1)a\bar a}-\hat{b}_{ia\bar a}\right) \eta_{ia\bar a}^{q}\leq\nu,\forall a<\bar a\in A,\\
& \left|\sum_{j\in[m_a]}\delta_{ija\bar a 1}t_{ja}-\sum_{j\in [m_{\bar a}]}\delta_{ija\bar a 2}t_{j\bar a}\right|\leq \eta_{ia\bar a}, \forall i\in [\hat{m}_{a\bar a}-1], a<\bar a\in A,\\
& (\bfx, t_{ja})\in \Omega_a(j), \forall j\in [m_a], a\in A
\end{aligned}\right\},
\end{equation}
where 
\[\delta_{ija\bar a 1}=\I\left(\left(\hat{b}_{ia\bar a},\hat{b}_{(i+1)a\bar a}\right]\subseteq \left(\frac{j-1}{m_a}, \frac{j}{m_a}\right]\right),\forall i\in [\hat{m}_{a\bar a}-1], j\in [m_a], a<\bar a\in A,\]
and
\[\delta_{ija\bar a 2}=\I\left(\left(\hat{b}_{ia\bar a},\hat{b}_{(i+1)a\bar a}\right]\subseteq \left(\frac{j-1}{m_{\bar a}}, \frac{j}{m_{\bar a}}\right]\right), \forall i\in [\hat{m}_{a\bar a}-1], j\in [m_{\bar a}], a<\bar a \in A.\]
\end{restatable}

\noindent\textit{Proof. } See Appendix~\ref{proof_thm_M3}.
\QEDA

According to \eqref{set_Qa}, we see that the continuous variable $t_{ja} = F^{-1}_{a}(j/m_a\mid\bfx)$ represents the $j/m_a$th smallest quantile value, and the binary variable $z_{ija}$ indicates whether up to the $i$th smallest quantile value is selected or not for each $i\in C_a, j\in [m_a],$ and $a\in A$. Therefore, we obtain the following monotonicity-based valid inequalities.
\begin{proposition}\label{prop_M3_ineq}
The following inequalities are valid for the \hyperref[model3_WD]{Quantile Formulation}
\begin{equation}\label{M3_ineq}
t_{ja}\leq t_{(j+1)a},z_{ija}\geq z_{i(j+1)a},\forall i\in C_a, j\in [m_a], a\in A.
\end{equation}
\end{proposition}

\subsection{Aggregate Quantile Formulation}\label{sec_M4}
Motivated by the quantile set $\Omega_a(k)$ in \eqref{set_Qa}, we develop another formulation using the aggregate quantiles in this subsection. We also show that this formulation can be quite strong compared to others. To begin with, let us define the aggregate quantile variable $\bar{t}$ and the aggregate quantile set $\bar{\Omega}_{a}(k)=\{(\bfx,\bar{t}_{ka})\in \mathcal X\times \Re: \sum_{i=1}^{k}F^{-1}_{a}(i/m_a\mid\bfx)=\bar{t}_{ka}\}$ for each $k\in [m_a]$ and $ a\in A$. 
Letting $(\bfx,\bar{w}_i)\in X_i$ for each $i\in[m]$, we present the following representation of the set $\bar{\Omega}_{a}(k)$ 
\begin{equation*}
\begin{aligned}
\bar{\Omega}_{a}(k)&=\left\{(\bfx,\bar{t}_{ka})\in \mathcal X\times \Re: \begin{aligned}
&(\bfx,\bar{w}_{i})\in X_i, \forall i\in C_a,  \\
&\bar{t}_{ka}\leq\min_{\bar{\bfz}}\left\{\sum_{i\in C_a}\bar{z}_{ika} \bar{w}_i:\bar{z}_{ika}\in\{0,1\},\forall i\in C_a, \sum_{i\in C_a}\bar{z}_{ika}=k \right\},\\
&  \bar{t}_{ka}\geq\min_{\bfz}\left\{\sum_{i\in C_a}z_{ika} \bar{w}_i:z_{ika}\in\{0,1\},\forall i\in C_a, \sum_{i\in C_a}z_{ika}=k \right\}
\end{aligned}\right\},
\end{aligned}
\end{equation*}
where similar to \eqref{set_Qa}, we let the binary variables $z_{ika},\bar{z}_{ika}$ indicate whether up to the $i$th smallest quantile values is selected or not for each $i\in C_a, k\in [m_a],$ and $a\in A$. By dualizing the first minimization problem and linearizing the bilinear terms in the second minimization problem, we arrive at the following MICP-R set.
\begin{restatable}{proposition}{propsetbaromega}\label{prop_set_bar_omega}
Suppose that $M_{i}\geq \max_{\bm{x}\in \mathcal X,\eqref{eq_tol}} |f(\bfx,{\bm\xi}_i)|$ for each $i\in [m]$. For each $k\in [m_a]$ and $a\in A$, the aggregate quantile set $\bar{\Omega}_{a}(k)$ is equivalent to
\begin{equation}\label{set_bar_Qa}
\begin{aligned}
\bar{\Omega}_{a}(k)&=\left\{(\bfx,\bar{t}_{ka})\in \mathcal X\times \Re: \begin{aligned}
& z_{ika}\in\{0,1\},(\bfx,\bar{w}_{i})\in X_i, \forall i\in C_a, \sum_{i\in C_a}z_{ika}=k, \\
& \bar{t}_{ka}\leq k \pi_{ka} -\sum_{i\in C_a}\rho_{ika},\pi_{ka}-\rho_{ika}\leq \bar{w}_i, \rho_{ika}\geq0, \forall i\in C_a,\\
& \bar{t}_{ka}\geq\sum_{i\in C_a} s_{ika},(s_{ika},z_{ika},\bar{w}_{i})\in\MC(0,1,-M_i,M_i), \forall i\in C_a
\end{aligned}\right\}.
\end{aligned}
\end{equation}
\end{restatable}

To represent the set $\F_{q}$ defined in \eqref{eq_F_q2}, we simply plug in the representation of the aggregate quantile sets $\{\bar{\Omega}_a(k)\}_{k\in [m_a],a\in A}$ into the representation \eqref{eq_F_q2}, which motivates the following formulation.
\begin{restatable}{theorem}{thmM4}\textbf{(Aggregate Quantile Formulation)}\label{thm_M4}
Suppose that the set $X_{i}=\{(\bfx, \bar{w}_{i})\in \mathcal X\times\Re: f(\bfx,{\bm\xi}_i)=\bar{w}_{i}\}$ is MICP-R and $M_{i}\geq \max_{\bm{x}\in \mathcal X,\eqref{eq_tol}} |f(\bfx,{\bm\xi}_i)|$ for each $i\in [m]$. We further define the aggregate quantile set $\bar{\Omega}_{a}(k)=\{(\bfx,\bar{t}_{ka})\in \mathcal X\times \Re: \sum_{i=1}^{k}F^{-1}_{a}(i/m_a\mid\bfx)=\bar{t}_{ka}\}$, which admits a MICP-R form \eqref{set_bar_Qa}. Then $\F_{q}$ can be represented as 
\begin{equation}\label{model4_WD}
\F_{q}=\left\{(\bm{x},\nu)\in \mathcal X\times \Re_+: \begin{aligned}
&\sum_{i\in [\hat{m}_{a\bar a}-1]}\left(\hat{b}_{(i+1)a\bar a}-\hat{b}_{ia\bar a}\right) \eta_{ia\bar a}^{q}\leq\nu,\forall a<\bar a\in A,\\
& \left|\sum_{j\in[m_a]}\delta_{ija\bar a 1}t_{ja}-\sum_{j\in [m_{\bar a}]}\delta_{ija\bar a 2}t_{j\bar a}\right|\leq \eta_{ia\bar a}, \forall i\in [\hat{m}_{a\bar a}-1], a<\bar a\in A,\\
& t_{ja}=\bar{t}_{ja}-\bar{t}_{(j-1)a}, (\bfx, \bar{t}_{ja})\in \bar{\Omega}_{a}(j), \bar{t}_{0a}=0,\forall j\in[m_a], a\in A
\end{aligned}\right\},
\end{equation}
where the parameters $\bm{\delta}$ are defined in \Cref{thm_M3}.
\end{restatable}
 \noindent\textit{Proof. }
The proof follows \Cref{thm_M3} with the fact that $t_{ja}=\bar{t}_{ja}-\bar{t}_{(j-1)a}$ for all $j\in[m_a],a\in A$.
\QEDA
 
\
We remark that the inequalities \eqref{M3_ineq} are also valid for the \hyperref[model4_WD]{Aggregate Quantile Formulation}.

\subsection{Summary of the Different Formulations}
The different formulations have their own strengths from the derivations according to their developments. Their formulation complexities are summarized in \Cref{table1}, where we suppress the term $O(|A|^2)$ for simplicity. In our numerical study, we observe that the \hyperref[model4_WD]{Aggregate Quantile Formulation} consistently outperforms the others in terms of computational time, which might be because it has the least amount of binary variables and the smallest big-M coefficients.

\begin{table}[htbp]
\centering
\small
\caption{Formulation Complexity Comparisons}\label{table1}
\renewcommand{\arraystretch}{1.2}
\begin{tabular}{|c|c|c|c|c|}
\hline
Formulation& \textbf{\# of Constraints} & \textbf{\begin{tabular}[c]{@{}c@{}}\# of Binary \\ Variables\end{tabular}} & \textbf{\begin{tabular}[c]{@{}c@{}}\# of Continuous \\ Variables\end{tabular}} & \textbf{\begin{tabular}[c]{@{}c@{}}Largest Big-M \\ Coefficient\end{tabular}} \\ \hline
\hyperref[WD_q_trans_micpr]{\textbf{Discretized}} & $O(m^2\log(m))$ & $O(m^2\log(m))$ & $O(m^2\log(m))$ & $\max_{i\in [m]}M_i$ \\ \hline
\hyperref[model2_WD1]{\textbf{Complementary}} & $O(m^2)$ & $O(m^2)$ & $O(m^2)$ & \begin{tabular}[c]{@{}c@{}}$\max_{a<\bar a\in A}\sum_{(i,j)\in C_a\times C_{\bar a}}$\\ $(M_i+M_j)^q$\end{tabular} \\ \hline
\hyperref[model3_WD]{\textbf{Quantile}} & $O(m^2)$ & $O(m^2)$ & $O(m^2)$ & $2\max_{i\in [m]}M_i$ \\ \hline
\hyperref[model4_WD]{\textbf{Aggregate Quantile}} & $O(m^2)$ & $O(m^2)$ & $O(m^2)$ & $\max_{i\in [m]}M_i$ \\ \hline
\end{tabular}
\end{table}

In the following, we show that the \hyperref[model3_WD]{Quantile Formulation} and the \hyperref[model4_WD]{Aggregate Quantile Formulation} can be stronger than the other two under some assumptions. 
\begin{restatable}{proposition}{propMoneMtwoconti}\label{prop_M1_M2_conti}
Suppose that the big-M coefficients $\bm{M},\hat{\bm{M}}$ are large enough as specified in the proof. Then, by relaxing the binary variables, 
\begin{enumerate}[(i)]
\item the continuous relaxation value of the \hyperref[WD_q_trans_micpr]{Discretized Formulation} is zero;
\item the continuous relaxation value of the \hyperref[model2_WD1]{Complementary Formulation} is zero;
\item 
the continuous relaxation value of the \hyperref[model3_WD]{Quantile Formulation} is 
\[\min_{\bfx\in \X}\max_{a<\bar{a}\in A}\left(\hat{b}_{(\hat{m}_{a\bar a})a\bar a}-\hat{b}_{(\hat{m}_{a\bar a}-1)a\bar a}\right)\left|F_{a}^{-1}(1\mid\bfx)-F_{\bar a}^{-1}(1\mid\bfx)\right|^q;\]
\item 
the continuous relaxation value of the \hyperref[model4_WD]{Aggregate Quantile Formulation} is at least
\[\min_{\bfx\in \X}\max_{a<\bar{a}\in A}\left|\frac{1}{m_a}\sum_{i\in C_{a}}F_{a}^{-1}\left(\frac{i}{m_a}\Bigm|\bfx\right)-\frac{1}{m_{\bar a}}\sum_{i\in C_{\bar a}}F_{\bar a}^{-1}\left(\frac{i}{m_{\bar a}}\Bigm|\bfx\right)\right|^q.\] 
\end{enumerate}
\end{restatable}
\noindent\textit{Proof. } See Appendix~\ref{proof_prop_M1_M2_conti}\QEDA

As a side product of \Cref{prop_M1_M2_conti}, we see that 
\begin{corollary}\label{cor_prop_M1_M2_conti}
For any $q\geq 1$, the continuous relaxation value of the \hyperref[model4_WD]{Aggregate Quantile Formulation} is at least as good as the Jensen bound presented in \Cref{sec_jensen}.
\end{corollary}

The continuous relaxations of all the formulations can have nonzero objective values if one optimizes the big-M coefficients or adds valid inequalities. We numerically test each formulation in \Cref{ex_regression} and observe that the \hyperref[model4_WD]{Aggregate Quantile Formulation} performs best overall.

\subsection{An Alternating Minimization (AM) Algorithm}\label{sec_AM}
When solving large instances with thousands of populations, the exact formulations in the previous subsections 
may suffer from slow convergence to find an optimal solution. Therefore, in this subsection, motivated by the representation in \Cref{lem_icf_WD}, we design a fast AM algorithm that can effectively solve \ref{eq_sp_fair} instances to near optimality. 

To this end, according to \eqref{eq_lem_inv_cdf_WD1}, we can recast \ref{eq_sp_fair} as
\begin{subequations}\label{eq_invese_dfso}
\begin{align}
v^*(q) = \min_{\bfx\in\mathcal X ,\nu} \quad & \nu,\\
\text{s.t.}\quad & 
\sum_{i\in [\hat{m}_{a\bar a}-1]}\left(\hat{b}_{(i+1)a\bar a}-\hat{b}_{ia\bar a}\right)\left|F^{-1}_{a}\left(\hat{b}_{(i+1)a\bar a}\mid\bfx\right)-F^{-1}_{\bar a}\left(\hat{b}_{(i+1)a\bar a}\mid\bfx\right)\right|^{q} \leq \nu, \forall a<\bar a\in A,\\
& \eqref{eq_tol}\nonumber,
\end{align}
\end{subequations}
which has been used to derive the \hyperref[model3_WD]{Quantile Formulation} and the \hyperref[model4_WD]{Aggregate Quantile Formulation} of \ref{eq_sp_fair}. This formulation is also valuable for deriving the AM algorithm. Specifically, we can run the AM algorithm as follows: (i) First, we pick a feasible solution $\bfx_0$ (e.g., an optimal solution that minimizes the total cost \eqref{eq_sp}); (ii) At iteration $t\geq 0$, we find the inverse distribution functions $\{F_a^{-1}(\cdot\mid\bfx_t)\}_{a\in A}$, which can be done via sorting with time complexity $O(m\log m)$; (iii) For each $i\in [\hat{m}_{a\bar a}-1]$ and $a<\bar a\in A$, let $f(\bfx_t, \bm\xi_{\hat{s}_a(i)}):=F^{-1}_{a}(\hat{b}_{(i+1)a\bar a}\mid\bfx_t)$ and  $f(\bfx_t, \bm\xi_{\hat{s}_{\bar a}(i)}):=F^{-1}_{\bar a}(\hat{b}_{(i+1)a\bar a}\mid\bfx_t)$; (iv) Next, we solve the following program by fixing the inverse distribution functions in the DFSO \eqref{eq_invese_dfso}:
\begin{subequations}\nonumber
\begin{align}
v_{t+1}(q) = \min_{\bfx\in\mathcal X ,\nu} \quad & \nu,\\
\text{s.t.}\quad & 
\sum_{i\in [\hat{m}_{a\bar a}-1]}\left(\hat{b}_{(i+1)a\bar a}-\hat{b}_{ia\bar a}\right)\left|f\left(\bfx,\bm\xi_{\hat{s}_a(i)}\right)-f\left(\bfx,\bm\xi_{\hat{s}_{\bar a}(i)}\right)\right|^{q} \leq \nu, \forall a<\bar a\in A,\\
& \eqref{eq_tol}\nonumber,
\end{align}
\end{subequations}
with an optimal solution $\bfx_{t+1}$; and (v) Let $t:=t+1$ and repeat Step (ii) to Step (iv) until the stopping criterion is invoked (e.g., $|v_t-v_{t+1}|<\bar{\epsilon}$ for some small threshold $\bar{\epsilon}$). The benefit of the proposed AM algorithm is that it completely eliminates the necessity of auxiliary binary variables introduced by the exact MICP-R formulations. In addition to its computational advantage, our numerical study shows that the proposed AM algorithm can successfully find optimal solutions in many instances.

\section{Two Lower Bounds for the Wasserstein Fairness Measure}\label{sec_lower_bounds}
In this section, we study a compact Jensen lower bound for type-$q$ Wasserstein fairness measure (i.e., $\WD_{q}^q(\bfx)$) and the well-known Gelbrich lower bound for type-$2$ Wasserstein fairness measure (i.e., $\WD_{2}^2(\bfx)$). In particular, we derive new conditions under which the Gelbrich bound is tight. To obtain the equivalent MICP-R formulations, we assume that $f(\bfx,\bm\xi)=\bm\xi^\top \bm{r}(\bfx)+ s(\bfx)$, where $\bm{r}(\bfx)=\bm{A}\bm{x}+\hat{\bm{a}}_0$ and $s(\bfx) = \hat{\bm{a}}_1^\top \bfx+ \hat{a}_2 $ are linear functions. For notational convenience, we define the mean and covariance matrix for each group $a\in A$ as $\bm\mu_a =\E_{\Pr}[\tilde{\bm\xi}_a]$ and $\bm\Sigma_a =\text{Cov}_{\Pr}[\tilde{\bm\xi}_a]$, respectively.
\subsection{The Jensen Bound for the $q$th Power of Type-$q$ Wasserstein Fairness Measure $\WD_{q}^q($\boldmath{$x$}\unboldmath{$)$}}\label{sec_jensen}
We first introduce the Jensen bound for $\WD_{q}^q(\bfx)$, which enables us to ascertain that the semidefinite relaxation of the Gelbrich bound is relatively weak. The following theorem establishes the relation between the Wasserstein fairness measure and the Jensen bound. 
\begin{theorem}[The Jensen Bound]
For any $q\geq 1$, $\WD_q^q(\bfx)$ is bounded by
\[\WD_q^q(\bfx)\geq \max_{a<\bar{a}\in A}\left|\bm\mu_{a}^{\top}\bm{r}(\bfx) -\bm\mu_{\bar{a}}^{\top}\bm{r}(\bfx) \right|^q := v_J(q).\]
\end{theorem}
\noindent\textit{Proof. }
For any $q$, $a<\bar{a}\in A$, and joint distribution $\Qe_{a,\bar{a}}$ of $f(\bfx,\tilde{\bm\xi}_{a})$ and $f(\bfx,\tilde{\bm\xi}_{\bar{a}})$ with marginals $\Pr_a,\Pr_{\bar a}$, we have
\begin{align*}
\E_{\Qe_{a,\bar{a}}}[|f(\bfx,{\bm\xi}_{a})-f(\bfx,{\bm\xi}_{\bar{a}})|^q]
& \geq \left|\E_{\Qe_{a,\bar{a}}}[f(\bfx,{\bm\xi}_{a})]-\E_{\Qe_{a,\bar{a}}}[f(\bfx,{\bm\xi}_{\bar{a}})]\right|^q, \\
& = \left|\E_{\Pr_{a}}[f(\bfx,{\bm\xi}_{a})]-\E_{\Pr_{\bar{a}}}[f(\bfx,{\bm\xi}_{\bar{a}})]\right|^q = \left|\bm\mu_{a}^{\top}\bm{r}(\bfx) -\bm\mu_{\bar{a}}^{\top}\bm{r}(\bfx) \right|^q.
\end{align*}
Here, the inequality is due to Jensen's inequality, and the first equality is because the random vectors ${\bm\xi}_{a},{\bm\xi}_{\bar{a}}$ are governed by the marginal distributions $\Pr_a,\Pr_{\bar a}$, respectively. Then, we obtain
\[\WD^q_q(\bfx)=\max_{a<\bar{a}\in A} W_{q}^q\left(\Pr_{a},\Pr_{\bar{a}}\right) \geq \max_{a<\bar{a}\in A}\left|\bm\mu_{a}^{\top}\bm{r}(\bfx) -\bm\mu_{\bar{a}}^{\top}\bm{r}(\bfx) \right|^q,\]
which completes the proof.\QEDA
 
This result gives rise to the following model for computing the Jensen bound for $\WD_q^q(\bfx)$: 
\begin{align}\label{lower_bound_Wq}
v_J(q) = \min_{\bfx\in\mathcal X,\nu} \left\{\nu:\left|\bm\mu_{a}^{\top}\bm{r}(\bfx) -\bm\mu_{\bar{a}}^{\top}\bm{r}(\bfx) \right|^q \leq \nu, \forall a<\bar{a}\in A,\eqref{eq_tol}\right\}.
\end{align}

\subsection{The Gelbrich Bound for the Squared Type-$2$ Wasserstein Fairness Measure $\WD_{2}^{2}($\boldmath{$x$}\unboldmath{$)$}}\label{sec_gelbrich}

When $q=2$, there is a popular Gelbrich bound for $ W_{2}^2\left(\Pr_{a},\Pr_{\bar{a}}\right)$ for any $a<\bar{a}\in A$, which has been studied in many optimal transport works (see, e.g., \citealt{kuhn2019wasserstein}). Formally, the Gelbrich bound for $W_{2}^2\left(\Pr_{a},\Pr_{\bar{a}}\right)$ is defined as follows.
\begin{definition}[The Gelbrich Bound, Theorem 2.1 in \citealt{gelbrich1990formula}]\label{def_gelbrich} For any $a<\bar{a}\in A$, the squared type-2 Wasserstein distance $W_{2}^2\left(\Pr_{a},\Pr_{\bar{a}}\right)$ is bounded by
\[W_{2}^2\left(\Pr_{a},\Pr_{\bar{a}}\right)\geq \left(\bm\mu_{a}^{\top}\bm{r}(\bfx) -\bm\mu_{\bar{a}}^{\top}\bm{r}(\bfx) \right)^2 + \left(\sqrt{\bm{r}(\bfx)^{\top}\bm\Sigma_{a}\bm{r}(\bfx)} - \sqrt{\bm{r}(\bfx)^{\top}\bm\Sigma_{\bar a}\bm{r}(\bfx)}\right)^2.
\]
\end{definition}

According to \Cref{def_gelbrich}, the Gelbrich bound can be computed via the following nonconvex program:
\begin{subequations}\label{gelbrich_bound}
\begin{align}
v_{G} = \min_{\bfx\in\mathcal X,\nu} \quad &\nu,\\
\text{s.t.} \quad &
 \left(\bm\mu_{a}^{\top}\bm{r}(\bfx) -\bm\mu_{\bar{a}}^{\top}\bm{r}(\bfx) \right)^2 + \left(\sqrt{\bm{r}(\bfx)^{\top}\bm\Sigma_{a}\bm{r}(\bfx)} - \sqrt{\bm{r}(\bfx)^{\top}\bm\Sigma_{\bar a}\bm{r}(\bfx)}\right)^2 \leq \nu, \forall a<\bar{a}\in A,\\
&  \eqref{eq_tol}.\notag
\end{align}
\end{subequations}
where $\bm\mu_{a}$ and $\bm\Sigma_{a}$ are the mean and covariance of $\tilde{\bm\xi}_{a}$ for all $a$.

Using the Cholesky decomposition $\bm\Sigma_a=\bm L_a\bm L_a^{\top}$ for each $a\in A$, we can recast \eqref{gelbrich_bound} as
\begin{subequations}\label{gelbrich_Cholesky}
\begin{align}
v_{G} = \min_{\bfx\in\mathcal X,\bfz,\nu} \quad &\nu,\\
\text{s.t.} \quad
& \bfz_{a} =\bm  L_{a}^{\top}\bm{r}(\bfx), \forall a\in A,\label{eq_z_Chol}\\
& \left(\bm\mu_{a}^{\top}\bm{r}(\bfx) - \bm\mu_{\bar{a}}^{\top}\bm{r}(\bfx) \right)^2 + \left( \|\bfz_{a}\|_2 - \|\bfz_{\bar{a}}\|_2 \right)^2 \leq \nu, \forall a<\bar{a}\in A,\\
& \eqref{eq_tol}\nonumber.
\end{align}
\end{subequations}

Our numerical study shows that the Gelbrich bound \eqref{gelbrich_Cholesky} can be very close to the true optimal value $v^*(2)$. 
Unfortunately, computing the Gelbrich bound \eqref{gelbrich_Cholesky} constitutes an intractable nonconvex program, which we formally prove to be generically NP-hard.
\begin{restatable}{theorem}{thmgelNP}\label{thm_gel_NP}
Computing the Gelbrich bound is strongly NP-hard even when $\epsilon=\infty$ and $|A|=2$.
\end{restatable}
 \noindent\textit{Proof. } See Appendix~\ref{proof_thm_gel_NP}.\QEDA

We remark that, in practice, one can compute the Gelbrich bound~\eqref{gelbrich_Cholesky} by employing off-the-shelf solvers, which are based on the spatial branch and bound algorithm. To further expedite the solution process, we can tighten the bounds of decision variables and auxiliary variables in formulation \eqref{gelbrich_Cholesky}, which significantly decreases the number of branch and bound nodes and thus accelerates the computation.

\noindent \textbf{AM Algorithm.} The complexity result motivates us to solve \eqref{gelbrich_Cholesky} using a highly effective AM method. We first rewrite the formulation \eqref{gelbrich_Cholesky} as
\begin{align*}
v_{G} = \min_{\bfx\in\mathcal X,\bfz, \nu} \quad &\nu,\\
\text{s.t.} \quad
&\left(\bm\mu_{a}^{\top}\bm{r}(\bfx) - \bm\mu_{\bar{a}}^{\top}\bm{r}(\bfx) \right)^2 + 2\|\bfz_{a}\|_2^2 + 2\|\bfz_{\bar{a}}\|_2^2-\left( \|\bfz_{a}\|_2 + \|\bfz_{\bar{a}}\|_2 \right)^2  \leq \nu, \forall a<\bar{a}\in A,
\\
&\eqref{eq_tol}, \eqref{eq_z_Chol}.\nonumber
\end{align*}
Using the convex conjugate representation, we have 
$$-\left( \|\bfz_{a}\|_2 + \|\bfz_{\bar{a}}\|_2 \right)^2 =\min_{w_{a\bar a}\geq 0,\bm{\alpha}_a,\bm{\alpha}_{\bar a}}\left\{-2 \bm{\alpha}_a^\top \bfz_{a} -2 \bm{\alpha}_a^\top \bfz_{a}+w_{a\bar a}^2:\|\bm{\alpha}_a\|_2\leq w_{a\bar a}, \|\bm{\alpha}_{\bar a}\|_2\leq w_{a\bar a} \right\}.$$
Thus, we can equivalently restate the formulation \eqref{gelbrich_Cholesky} as
\begin{subequations}\label{gelbrich_AM}
\begin{align}
v_{G} = \min_{\bfx\in\mathcal X,\bfz,\bm{w},\bm{\alpha}, \nu}\quad&\nu,\\
\text{s.t.} \quad
&\left(\bm\mu_{a}^{\top}\bm{r}(\bfx) - \bm\mu_{\bar{a}}^{\top}\bm{r}(\bfx) \right)^2 + 2\|\bfz_{a}\|_2^2 + 2\|\bfz_{\bar{a}}\|_2^2-2 \bm{\alpha}_a^\top \bfz_{a} \notag\\
&-2\bm{\alpha}_a^\top \bfz_{a}+w_{a\bar a}^2  \leq \nu, |\bm{\alpha}_a\|_2\leq w_{a\bar a}, \|\bm{\alpha}_{\bar a}\|_2\leq w_{a\bar a}, \forall a<\bar{a}\in A, \\
&\eqref{eq_tol}, \eqref{eq_z_Chol}.\nonumber
\end{align}
\end{subequations}
In the AM method, at each iteration $t$, given a solution $(\bfx_t, \bfz_{t},\nu_t)$, we  compute the solution $(\bm{w}_t, \bm{\alpha}_{t},\bar{\nu}_t)$ in closed-form, as follows:
\begin{align*}
&\bar{\nu}_t=
\max_{a<\bar a\in A}\left\{\left(\bm\mu_{a}^{\top}\bm{r}(\bfx_t) - \bm\mu_{\bar{a}}^{\top}\bm{r}(\bfx_t) \right)^2 +\left( \|\bfz_{at}\|_2 - \|\bfz_{\bar{a}t}\|_2 \right)^2\right\},\\
& w_{a\bar a t}=\|\bfz_{at}\|_2 + \|\bfz_{\bar{a}t}\|_2, \bm{\alpha}_{at}=\frac{w_{a\bar a t}}{\|\bfz_{at}\|_2}\bfz_{at}, \bm{\alpha}_{\bar a t}=\frac{w_{a\bar a t}}{\|\bfz_{\bar at}\|_2}\bfz_{\bar{a}t}.
\end{align*}
Then we fix the values of $(\bm{w}_t, \bm{\alpha}_{t})$ and resolve \eqref{gelbrich_AM} with respect to the variables $(\bfx, \bfz,\nu)$. The procedure is repeated until we reach a prescribed tolerance. Our numerical study finds that the AM approach works extremely well in quickly finding near-optimal solutions.

\noindent\textbf{Semidefinite Programming Relaxation.} Alternatively, in the Gelbrich bound formulation \eqref{gelbrich_Cholesky}, let us introduce a new variable $\sigma_a=\|\bfz_{a}\|_2$ for each $a\in A$. For each pair $a<\bar a\in A$, let us denote $\bm s_{a\bar{a}} = \begin{bmatrix} \sigma_{a} & \bfz_{a} & \sigma_{\bar{a}}  & \bfz_{\bar{a}} \end{bmatrix}^{\top}$ and $ \bm{Z}_{a\bar{a}} =  \bm s_{a\bar{a}}\cdot \bm s_{a\bar{a}}^{\top}$.
Then one can show that the Gelbrich bound \eqref{gelbrich_Cholesky} can be converted to a semidefinite programming formulation with rank-one constraint, as follows:
\begin{subequations}\label{gelbrich_SDP_rank1}
\begin{align}
v_{G} = \min_{\bfx\in\mathcal X,\bm s, \bm Z,\nu} \quad&\nu,\\
\text{s.t.} \quad
& \left(\bm\mu_{a}^{\top}\bm{r}(\bfx) - \bm\mu_{\bar{a}}^{\top}\bm{r}(\bfx) \right)^2 + \notag\\
&\left(Z_{a\bar{a}11} - 2Z_{a\bar{a}1(n+2)} +  Z_{a\bar{a}(n+2) (n+2)}\right)\leq \nu, \forall a< \bar{a}\in A, \label{gelbrich_SDP_rank1_eq1}\\
& \bfz_{a} = \bm  L_{a}^{\top}\bm{r}(\bfx)\label{gelbrich_SDP_rank1_eq2}, \sigma_a\geq 0,\forall a\in A,\\
& Z_{a\bar{a}11} = \sum_{i=2}^{n+1}Z_{a\bar{a}ii}, Z_{a\bar{a}(n+2) (n+2)} = \sum_{i=n+3}^{2n+2}Z_{a\bar{a}ii}, \forall a< \bar{a}\in A, \label{gelbrich_SDP_rank1_eq5}\\
& \bm s_{a\bar{a}} = \begin{bmatrix} \sigma_{a} & \bfz_{a} & \sigma_{\bar{a}}  & \bfz_{\bar{a}} \end{bmatrix}^{\top}, \forall a< \bar{a}\in A, \label{gelbrich_SDP_rank1_eq6}\\
& \bm{Z}_{a\bar{a}} =  \bm s_{a\bar{a}}\cdot \bm s_{a\bar{a}}^{\top}, \forall a< \bar{a}\in A, \label{gelbrich_SDP_Zvar}\\
& \eqref{eq_tol}\nonumber.
\end{align}
\end{subequations}
The rank one constraints in \eqref{gelbrich_SDP_Zvar} are difficult to handle in practice. A simple way is to
relax \eqref{gelbrich_SDP_Zvar} as the semidefinite inequalities
\[\bm{Z}_{a\bar{a}} \succeq \bm s_{a\bar{a}}\cdot \bm s_{a\bar{a}}^{\top}, \forall a< \bar{a}\in A.\]
Using the Schur complement, we obtain the semidefinite relaxation of the Gelbrich bound \eqref{gelbrich_Cholesky} as
\begin{subequations}\label{gelbrich_SDP_relax}
\begin{align}
v_{\underline G} = \min_{\bfx\in\mathcal X,\bm s,\bm Z,\nu} \quad &\nu,\\
\text{s.t.} \quad
& \begin{bmatrix} 1 & \bm s_{a\bar{a}}^{\top}\\ \bm s_{a\bar{a}} &  \bm{Z}_{a\bar{a}} \end{bmatrix}\succeq 0,\forall a< \bar{a}\in A, \label{gelbrich_SDP_rank1_schur}\\
& \eqref{eq_tol},\eqref{gelbrich_SDP_rank1_eq1}-\eqref{gelbrich_SDP_rank1_eq6}.\nonumber
\end{align}
\end{subequations}

We see that the semidefinite relaxation \eqref{gelbrich_SDP_relax} is stronger than the type-2 Jensen bound $v_J(2)$ in \eqref{lower_bound_Wq} since $\bm{Z}_{a\bar{a}}$ is positive semidefinite and $\left(Z_{a\bar{a}11} - 2Z_{a\bar{a}1(n+2)} +  Z_{a\bar{a}(n+2) (n+2)}\right)\geq0$ for every pair $a<\bar{a}\in A$. On the other hand, if we allow the relative tolerance of the semidefinite constraints (see \citealt{aps2019mosek}), then for up to any prescribed tolerance, we can show that $v_{\underline G}\leq  v_J(2)$. This result is summarized in the following proposition.
\begin{restatable}{proposition}{thmgelSDP}\label{thm_gel_SDP}
The semidefinite relaxation \eqref{gelbrich_SDP_relax} of the Gelbrich bound model satisfies $v_{\underline G}\geq  v_J(2)$. On the other hand, for any relative tolerance  $\beta>0$ of the semidefinite constraints in \eqref{gelbrich_SDP_rank1_schur} such that $\bm{Z}_{a\bar{a}} -\bm s_{a\bar{a}}\cdot \bm s_{a\bar{a}}^{\top}\succeq -\beta \lambda_{min}^+(\bm{Z}_{a\bar{a}}) \bm I_{2n+2}$, where $\lambda_{min}^+(\cdot )$ denotes the smallest nonzero eigenvalue, we have $v_{\underline G}\leq  v_J(2)$. 
\end{restatable}
 \noindent\textit{Proof. }
See Appendix~\ref{proof_thm_gel_SDP}.\QEDA

\subsection{Tightness of the Gelbrich Bound}

\noindent \textbf{The Same Univariate Marginal Distribution Condition:} In the literature, it is known 
that the Gelbrich bound is tight when the random parameters $\{\tilde{\bm\xi}_a\}_{a\in A}$ are asymptotically elliptical as $m_a\rightarrow \infty$ for all $a\in A$. We generalize this result by establishing a weaker condition that achieves the tightness of the Gelbrich bound. Our result shows that when the random utility functions of different groups can be linearly transformed to the same univariate random variable, then the Gelbrich bound is asymptotically tight.
\begin{restatable}{theorem}{thmexactgelbrich}\label{thm_exact_gelbrich}Suppose that for any pair $a<\bar a\in A$, the optimal comonotonic random variables $(f(\bfx,\tilde {\bm\xi}_{a})-\bm\mu_{a}^{\top}\bm{r}(\bfx)-s(\bfx), f(\bfx, \tilde{\bm\xi}_{\bar a})-\bm\mu_{\bar a}^{\top}\bm{r}(\bfx)-s(\bfx))\xrightarrow{m_{a}\rightarrow \infty,m_{\bar a}\rightarrow \infty}(\sqrt{\bm{r}(\bfx)^{\top}\bm\Sigma_{a}\bm{r}(\bfx)}\tilde u, \sqrt{\bm{r}(\bfx)^{\top}\bm\Sigma_{\bar a}\bm{r}(\bfx)}\tilde u)$ for a univariate random variable $\tilde u$ with zero mean and unit variance. Then the Gelbrich bound is asymptotically tight.
\end{restatable}
 \noindent\textit{Proof. } See Appendix~\ref{proof_thm_exact_gelbrich}. \QEDA

\Cref{thm_exact_gelbrich} shows that the tightness of the Gelbrich bound applies to a much broader family of distributions than elliptical. In fact, from the proof, we can see that the Gelbrich bound is derived using the Cauchy-Schwarz inequality, i.e.,
\[\E_{\Qe_{a,\bar{a}}}\left[\left(f(\bfx,\tilde{\bm\xi}_{a})-\bm\mu_{a}^{\top}\bm{r}(\bfx)-s(\bfx)\right)\left(f(\bfx,\tilde{\bm\xi}_{\bar{a}})-\bm\mu_{\bar a}^{\top}\bm{r}(\bfx)-s(\bfx)\right)\right]\leq \sqrt{\bm{r}(\bfx)^{\top}\bm\Sigma_{a}\bm{r}(\bfx)} \cdot \sqrt{\bm{r}(\bfx)^{\top}\bm\Sigma_{\bar a}\bm{r}(\bfx)}.\]
Thus, the tightness result holds whenever there exists a joint distribution such that the Cauchy-Schwarz inequality becomes equality.

More importantly, we can theoretically bound the gap between the optimal Gelbrich bound $v_G$ and the optimal value of \ref{eq_sp_fair} $v^*(2)$ under type $q=2$ Wasserstein distance.
\begin{restatable}{theorem}{thmexactgelbrichgap}\label{thm_exact_gelbrich_gap}Suppose that for any group $a\in A$, the individual samples $\{\bm{\xi}_i\}_{i\in C_a}$ satisfy $f(\bfx,\bm\xi_{i})-\bm\mu_{a}^{\top}\bm{r}(\bfx)-s(\bfx)\stackrel{\text{d}}{=}\sqrt{\bm{r}(\bfx)^{\top}\bm\Sigma_{a}\bm{r}(\bfx)}u_i$ for each $i\in C_a$, where $\{u_i\}_{i\in C_a}$ are i.i.d. samples of a univariate sub-Gaussian random variable $\tilde u_a$ with zero mean and unit variance, and $\{\tilde u_a\}_{a\in A}$ obey the same distribution. Then with probability at most $1-\hat{\eta}$ such that $\hat{\eta}>0$ is small, we have $$v^*(2)-\bar{C}_1(\hat{\eta}\min_{a\in A}\sqrt{m_{a}})^{-1}\leq v_G\leq v^*(2)$$ for some positive constant $\bar{C}_1$. 
\end{restatable}
 \noindent\textit{Proof. }See Appendix~\ref{proof_thm_exact_gelbrich_gap}.\QEDA

\noindent\textbf{Different Groups with Proportional Covariances:} 
The result in \Cref{thm_exact_gelbrich} necessitates the same marginal distributions. We relax this assumption by establishing another tightness condition, such that the marginal distributions of different groups can be distinct.
\begin{restatable}{theorem}{thmexactgelbrichvtwo}\label{thm_exact_gelbrich_v2}Suppose that for any pair $a<\bar a\in A$, the optimal comonotonic random variables $(f(\bfx,\tilde{\bm\xi}_{a})-\bm\mu_{a}^{\top}\bm{r}(\bfx)-s(\bfx), f(\bfx,\tilde{\bm\xi}_{\bar a})-\bm\mu_{\bar a}^{\top}\bm{r}(\bfx)-s(\bfx))\xrightarrow{m_{a}\rightarrow \infty,m_{\bar a}\rightarrow \infty}(\hat{\bm\xi}_{ a}^{\top}\bm{r}(\bfx), \hat{\bm\xi}_{\bar a}^{\top}\bm{r}(\bfx))$, where the random vectors $\hat{c}_a^{-1}\hat{\bm\xi}_{ a},\hat{c}_{\bar{a}}^{-1}\hat{\bm\xi}_{\bar a}$ obey the same distribution with zero mean and covariance matrix $\bm{\Sigma}_{a\bar a}$ for some positive parameters $\hat{c}_a, \hat{c}_{\bar{a}}$. Then the Gelbrich bound is asymptotically tight.
\end{restatable}
 \noindent\textit{Proof. }See Appendix \ref{proof_thm_exact_gelbrich_v2}.\QEDA

Similar to \Cref{thm_exact_gelbrich_gap}, we can theoretically bound the gap between the optimal Gelbrich bound $v_G$ and the optimal value of \ref{eq_sp_fair} $v^*(2)$ under type $q=2$ Wasserstein distance.
\begin{theorem}\label{thm_exact_gelbrich_gap_v2}Suppose that for any group $a\in A$, the individual samples $\{\bm{\xi}_i\}_{i\in C_a}$ satisfy $f(\bfx,\bm\xi_{i})-\bm\mu_{a}^{\top}\bm{r}(\bfx)-s(\bfx):=\bm\xi_{i}^{\top}\bm{r}(\bfx)$ for each $i\in C_a$, where $\{\bm\xi_i\}_{i\in C_a}$ are i.i.d. and sampling from $\hat{\bm\xi}_{ a}$ and 
the random vectors $\hat{c}_a^{-1}\hat{\bm\xi}_{ a},\hat{c}_{\bar{a}}^{-1}\hat{\bm\xi}_{\bar a}$ obey the same sub-Gaussian distribution with zero mean and covariance matrix $\bm{\Sigma}_{a\bar a}$. Then with probability at most $1-\hat{\eta}$ such that $\hat{\eta}>0$ is small, we have $$v^*(2)-\bar{C}_2(\hat{\eta}\min_{a\in A}\sqrt{m_{a}})^{-1}\leq v_G\leq v^*(2)$$ for some positive constant $\bar{C}_2$. 
\end{theorem}
 
\noindent\textit{Proof. }The proof is similar to that of \Cref{thm_exact_gelbrich_gap} and is thus omitted.
\QEDA

\section{Numerical Study}\label{sec_numerical_study}
In this section, we apply our framework to several fair optimization problems. We consider the fair regression problem and the fair allocation problem of scarce medical resources. An additional numerical study on the fair knapsack problem can be found in Appendix \ref{sec_fair_knapsack}. All the instances in this section are executed in Python 3.7 with calls to Gurobi 10.0.0 on a PC with an Apple M2 Pro processor and 16GB of memory.

\subsection{Fair Regression}\label{ex_regression}
Consider the regression problem aiming to predict the response vector $\bfy\in\Re^m$ using features $\bm{\xi}\in\Re^{m\times n}$, where the loss function is given by the mean squared error (MSE) $Q(\bfx,\bm{\xi}_i) = |\bm{\xi}_i^{\top}\bfx - y_i|^2$   or the mean absolute error (MAE) $Q(\bfx,\bm{\xi}_i) = |\bm{\xi}_i^{\top}\bfx - y_i|$. 
In terms of demographic parity fairness, we choose the utility function of the fair regression problem to be $f(\bfx,\bm{\xi}) = \bm{\xi}^{\top}\bfx$.
We conduct two experiments to test the proposed methods: (i) using hypothetical data to evaluate the performance of the exact formulations, AM algorithm, and two lower bounds, and (ii) using real data to compare \ref{eq_sp_fair} against two state-of-the-art methods.

\subsubsection{Formulation comparisons}
We compare the proposed methods for solving fair regression with MAE, where we (i) test the exact formulations on small populations and (ii) test the AM algorithm and lower bounds on large populations. In this experiment, we choose $|A|=2$ (i.e., we study the fairness among two groups) and $q=2$ (i.e., we consider type-2 Wasserstein distance), and we set $\epsilon=10\%$ as the inefficiency level. The hypothetical data is generated in the following manner. 
The response $\tilde{y}$ is generated from $\tilde{y}=\tilde{\bm{\xi}}^\top(\bfx^{0}) +\text{noise}$. The first $\lfloor n/2\rfloor$ components of the vector $\bfx^{0}$ are randomly sampled i.i.d. from the uniform distribution $\text{Unif}(-1, 0)$, the next $\lfloor n/2\rfloor-1$   components are sampled from $\text{Unif}(0, 10)$, and the last component is set to zero. The last component $\tilde{\xi}_{\kappa}\in\{-1,1\}$ of the vector $\tilde{\bm{\xi}}$ corresponds to the sensitive attribute. In the generated dataset, the first $\lceil m/2\rceil$ data points are assigned with the sensitive attribute 
$\xi_{\kappa}=-1$, where their features $(\tilde{\xi}_{j})_{j\in [\kappa-1]}$ 
are independently drawn from $\{\text{Unif}(0, j)\}_{j\in [\kappa-1]}$. The remaining data points are assigned $\xi_{\kappa}=1$, where their features $(\tilde{\xi}_{j})_{j\in [\kappa-1]}$ are independently drawn from $\{\text{Unif}(0, j+2)\}_{j\in [\kappa-1]}$. 
The $\text{noise}$ follows the uniform distribution $\text{Unif}(-0.1, 0.1)\times \E[\tilde{\bm{\xi}}]^{\top}(\bfx^{0})$.

In the first comparison, we generate data sets of a small population with sizes $m\in\{15,20,\dots,100\}$ and feature dimension $\kappa=10$ to compare the exact formulations against the AM algorithm of \ref{eq_sp_fair} and the two lower bounds. In the second comparison, we generate data sets of a large population with sizes $m\in\{100,200,\dots,3{,}000\}$  to illustrate the solution quality of the AM algorithm and two lower bounds. We solve the \hyperref[set_F_model1]{Vanilla Formulation}, the four exact MICP-R formulations, the AM algorithm in Section \ref{sec_AM}, the \hyperref[lower_bound_Wq]{Jensen bound}, and the \hyperref[gelbrich_Cholesky]{Gelbrich bound} in the first comparison. We test the AM algorithm, the \hyperref[lower_bound_Wq]{Jensen bound}, and the \hyperref[gelbrich_Cholesky]{Gelbrich bound} in the second comparison. 
Particularly, the AM algorithm of \ref{eq_sp_fair} in Section \ref{sec_AM} is initialized with the Gelbrich bound solution obtained by executing its corresponding AM algorithm described in Section \ref{sec_gelbrich}. 

In the first comparison, we report each instance's objective value, lower bound, optimality gap, and running time. Let ``Obj.Val'' denote the objective value and ``LB'' denote the lower bound. We use the dashed line ``--'' if ``Obj.Val'' is not available. The optimality gap denoted by ``Gap'' is computed by $\text{(UB-LB)/UB}\times100\%$, where we use the optimal objective value as UB if available. 
We define the best upper bound as the smallest ``Obj.Val'' of the exact formulations and the AM algorithm of \ref{eq_sp_fair}, and the best lower bound as the largest  ``LB'' of the exact formulations 
and ``Obj.Val'' of \hyperref[gelbrich_Cholesky]{Gelbrich bound}. For some instances, ``Obj.Val'' may not be available for the exact formulations. In this case, we use the best upper bound to compute their optimality gaps, use the best lower bound to compute the AM's optimality gap, and use the AM's objective value to compute the \hyperref[lower_bound_Wq]{Jensen bound}'s and \hyperref[gelbrich_Cholesky]{Gelbrich bound}'s optimality gaps. The running time in seconds is denoted as ``Time''. We set the time limit to 3,600 seconds. 
In the second comparison, we plot the gaps between AM and the two lower bounds over $10$ replications, where the gap is computed by $\text{(UB-LB)/UB}\times100\%$. We report the mean and standard deviation of the gaps and also illustrate the average running time of each method.

The first comparison results are displayed in Tables \ref{table_MIP_012}-\ref{table_MIP_Conti}. The \hyperref[set_F_model1]{Vanilla Formulation} cannot solve the small population instances within the time limit. The upper bounds of the \hyperref[set_F_model1]{Vanilla Formulation} tend to be close to the optimal value, while the lower bounds are nearly zero. In fact, the gap of the \hyperref[set_F_model1]{Vanilla Formulation} is $100\%$ when the population size of the instance is $m\geq30$. 
The \hyperref[WD_q_trans_micpr]{Discretized Formulation} can solve the instance with $m=15$. The quality of the incumbent solution at the time limit then deteriorates rapidly as $m$ increases. 
The \hyperref[model2_WD1]{Complementary Formulation} performs similarly to the \hyperref[WD_q_trans_micpr]{Discretized Formulation}. Its upper bounds are often worse than other formulations, and the lower bounds are always zero for all the instances. This demonstrates the weakness of the \hyperref[WD_q_trans_micpr]{Discretized Formulation} and the \hyperref[model2_WD1]{Complementary Formulation}, which is consistent with \Cref{prop_M1_M2_conti}. The performances of the \hyperref[model3_WD]{Quantile Formulation} and the \hyperref[model4_WD]{Aggregate Quantile Formulation} in Table \ref{table_MIP_34} are significantly better. The \hyperref[model3_WD]{Quantile Formulation} is able to solve instances up to $m\leq40$ to optimality, and it returns nonzero lower bounds except for the last instance. The optimality gap of the \hyperref[model3_WD]{Quantile Formulation} becomes larger as $m$ increases. Remarkably, the \hyperref[model4_WD]{Aggregate Quantile Formulation} can solve instances with $m\leq60$ and $m=70$ to optimality. The running time for each instance is less than $10$ seconds when the population size is $m\leq40$.  The \hyperref[model4_WD]{Aggregate Quantile Formulation} cannot be solved optimally for larger instances; however, it still consistently provides high quality lower bounds with small gaps. We observe that the upper bounds of the \hyperref[model3_WD]{Quantile Formulation} and the \hyperref[model4_WD]{Aggregate Quantile Formulation} may not be available for instances with large $m$. This is potentially due to these two MICP-R formulations having many variables and constraints, which causes the solver to have difficulty finding a feasible solution for large instances. Therefore, we instead use the AM algorithm to solve instances for which the \hyperref[model4_WD]{Aggregate Quantile Formulation} cannot provide an optimal solution within the time limit.

In fact, as shown in Table \ref{table_AM_Jen_Gel}, the AM algorithm provides very near-optimal solutions to instances of a small population using less than one second. It has a zero gap for most instances when the optimal solution is available, that is, $m\leq60$ and $m=70$. Its solution is close to the best lower bound when the optimal solution is unavailable, where the gap is less than $7\%$. In particular, the AM algorithm has better objective values than the \hyperref[set_F_model1]{Vanilla Formulation} for all instances in this experiment.
On the other hand, the \hyperref[lower_bound_Wq]{Jensen bound} has a gap of around $30\%$ for each instance, and its running time is short due to the simplicity of its model formulation. On the contrary, the \hyperref[gelbrich_Cholesky]{Gelbrich bound}'s gap decreases and running time slightly increases when the population size $m$ increases. The gap of the \hyperref[gelbrich_Cholesky]{Gelbrich bound} 
is around $15\%$ when the population size is $m\geq50$. Since the number of features $\kappa=10$ is small, the \hyperref[gelbrich_Cholesky]{Gelbrich bound} model can solve all instances to optimality, where each instance's running time is less than $3$ seconds. Besides, we also compute the continuous relaxation values of the exact MICP formulations. 
In Table \ref{table_MIP_Conti}, the continuous relaxation values of the first three formulations are zero for most instances. The \hyperref[model4_WD]{Aggregate Quantile Formulation} always has a nonzero continuous relaxation value, and it is greater than the objective value of the \hyperref[lower_bound_Wq]{Jensen bound} as shown in \Cref{cor_prop_M1_M2_conti}. The continuous relaxation gap of the \hyperref[model4_WD]{Aggregate Quantile Formulation} is around $30\%$ overall, which numerically verifies this formulation's strength. 

\begin{figure}[htbp]
\begin{subfigure}[Gap between AM and the Jensen Bound]{{
\centering\includegraphics[width=0.3\textwidth]{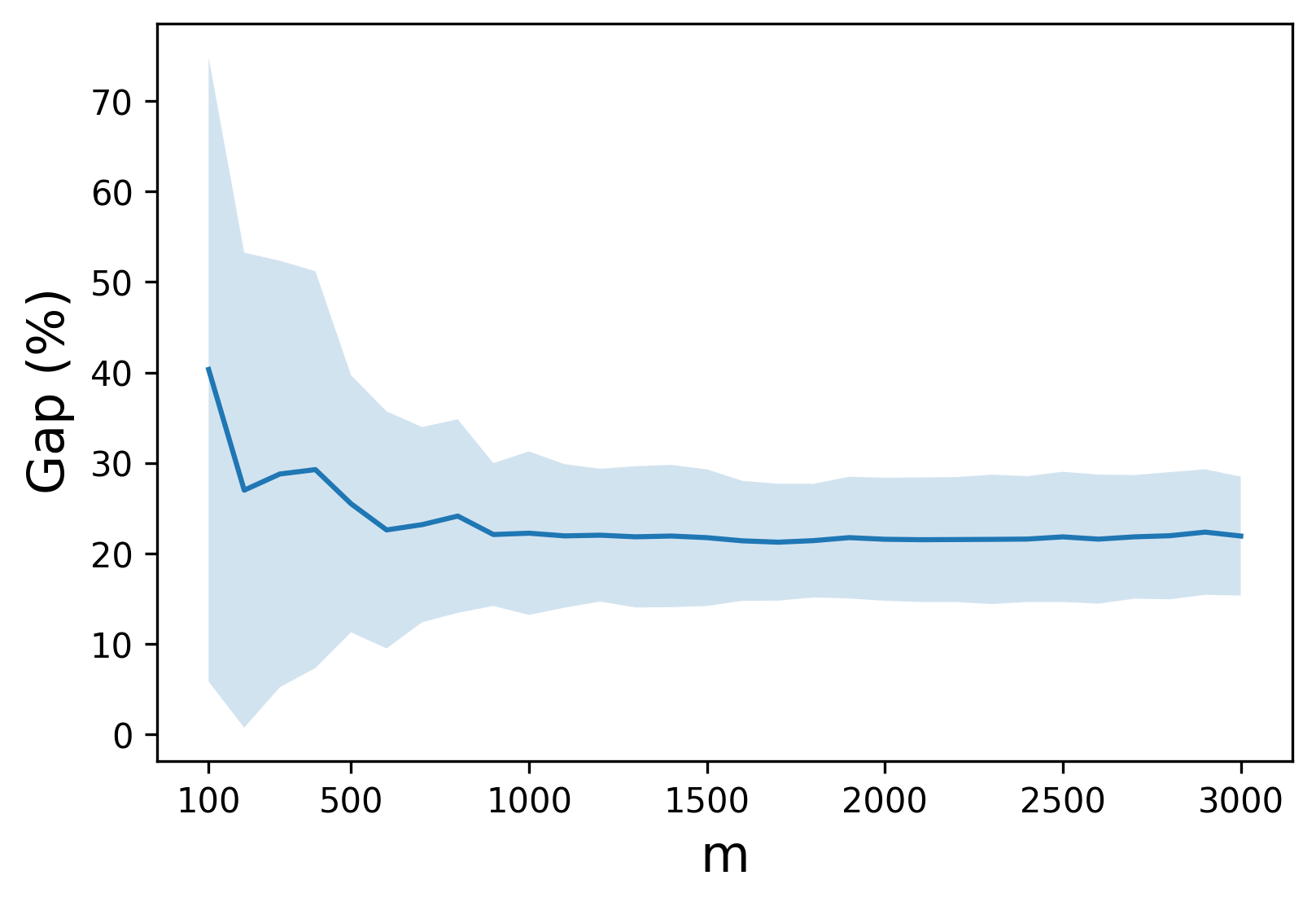}}}
\end{subfigure}
\hfill
\begin{subfigure}[Gap between AM and the Gelbrich Bound]{{
\centering\includegraphics[width=0.3\textwidth]{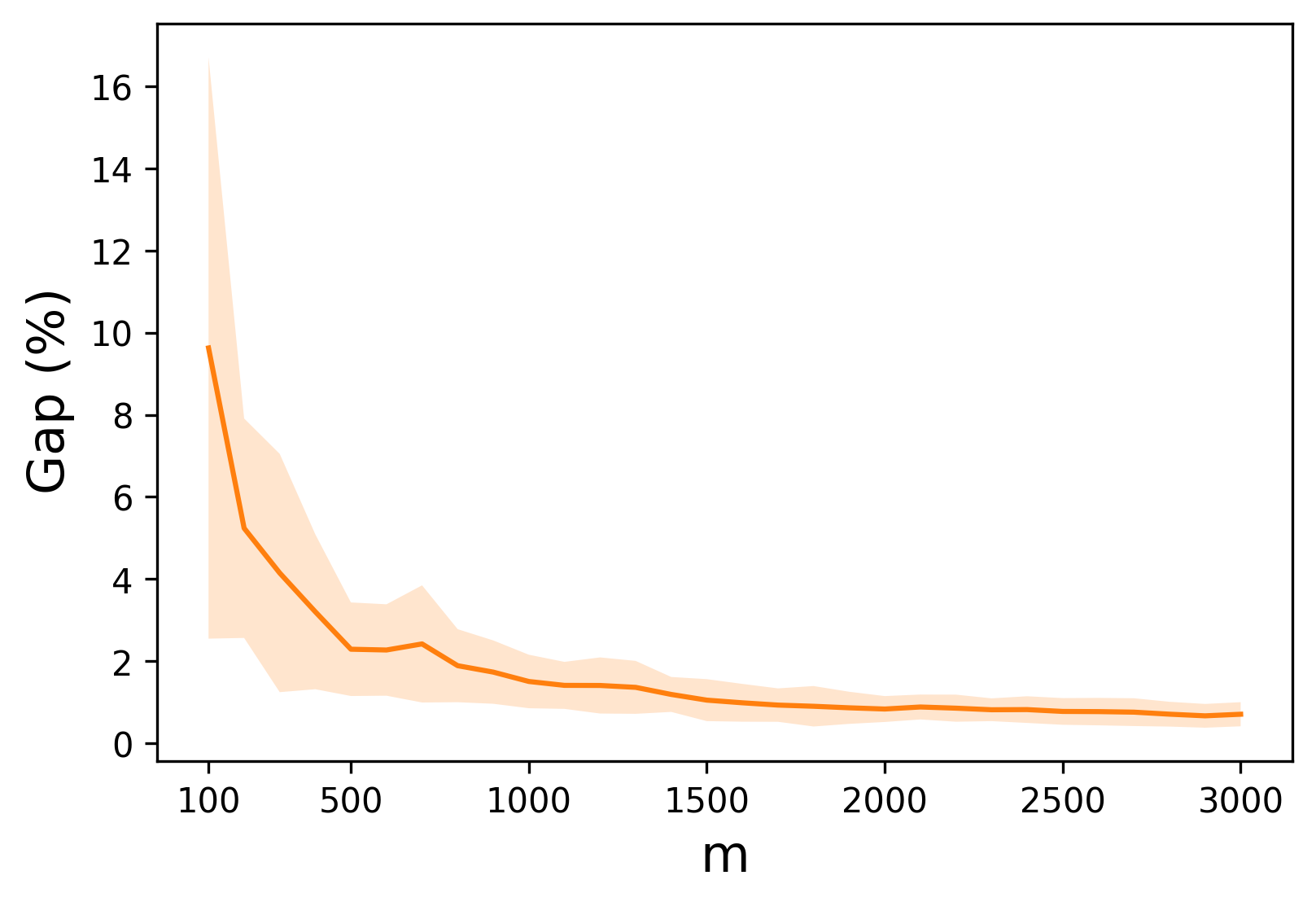}}}
\end{subfigure}
\hfill
\begin{subfigure}[Running time]{{
\centering\includegraphics[width=0.3\textwidth]{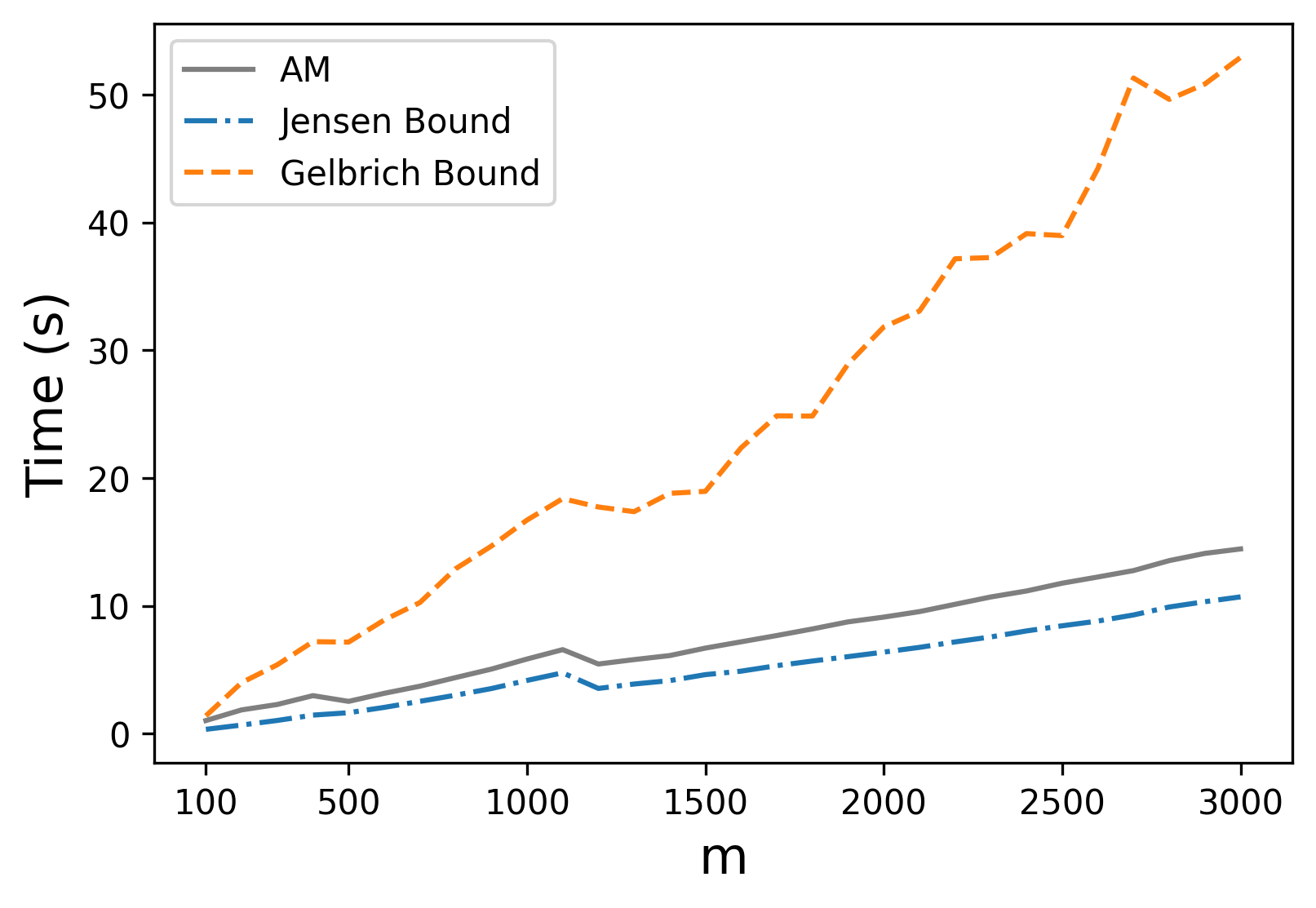}}\label{fig_time}}
\end{subfigure}
\caption{Gap between AM and the two lower bounds for a large population size. The mean and standard deviation of the gap over 10 replications, as well as the average running time of each method, are illustrated.}\label{fig_convergence}
\end{figure}

The second comparison is presented in Figure \ref{fig_convergence}. We see that both gaps stabilize when the population size is large enough. The gap between AM and the \hyperref[lower_bound_Wq]{Jensen bound} decreases from $40\%$ to $21\%$ when the population size $m$ grows from $100$ to $1{,}000$. This gap is around $21\%$ when the population size $m\geq 1{,}000$. The gap between AM and the \hyperref[gelbrich_Cholesky]{Gelbrich bound} drops from $10\%$ to $1\%$ when the population size $m$ grows from $100$ to $1{,}500$. This gap decreases to $0.8\%$ after $m=1{,}500$. The small gap between AM and the \hyperref[gelbrich_Cholesky]{Gelbrich bound} verifies that the solution of AM is near optimal and demonstrates the strength of the \hyperref[gelbrich_Cholesky]{Gelbrich bound}. Meanwhile, the running time of these methods grows slowly. Since the \hyperref[gelbrich_Cholesky]{Gelbrich bound} formulation is nonconvex, it requires a longer time to solve large population instances. Figure \ref{fig_time} shows that AM is much faster than the \hyperref[gelbrich_Cholesky]{Gelbrich bound}, and the \hyperref[lower_bound_Wq]{Jensen bound} is slightly faster than AM. When $m=3{,}000$, AM, the \hyperref[lower_bound_Wq]{Jensen bound}, and the \hyperref[gelbrich_Cholesky]{Gelbrich bound} take $15$, $11$, and $53$ seconds on average, respectively. The stable and efficiently solvable lower bound solutions are useful to initialize AM and verify its solution quality.

The two comparisons in this experiment confirm the effectiveness of the proposed methods in solving \ref{eq_sp_fair}. In practice, we suggest choosing the \hyperref[model4_WD]{Aggregate Quantile Formulation} to solve fair decision-making problems with a small population and switch to the AM method if the population size is large, where we can use the \hyperref[lower_bound_Wq]{Jensen bound} or the \hyperref[gelbrich_Cholesky]{Gelbrich bound} to initialize and establish the quality of the AM method.

\begin{table}[htbp]
\centering
\caption{Results of Exact MICP Formulations}\label{table_MIP_012}
\tiny
\begin{tabular}{c|rrrr|rrrr|rrrr}
\hline
\multirow{2}{*}{m} & \multicolumn{4}{c|}{\hyperref[set_F_model1]{Vanilla Formulation}} & \multicolumn{4}{c|}{\hyperref[WD_q_trans_micpr]{Discretized Formulation}} & \multicolumn{4}{c}{\hyperref[model2_WD1]{Complementary Formulation}} \\ \cline{2-13}
& \multicolumn{1}{c}{Obj.Val} & \multicolumn{1}{c}{LB} & \multicolumn{1}{c}{Gap (\%)} & \multicolumn{1}{c|}{Time} & \multicolumn{1}{c}{Obj.Val} & \multicolumn{1}{c}{LB} & \multicolumn{1}{c}{Gap (\%)} & \multicolumn{1}{c|}{Time} & \multicolumn{1}{c}{Obj.Val} & \multicolumn{1}{c}{LB} & \multicolumn{1}{c}{Gap (\%)} & \multicolumn{1}{c}{Time} \\ \hline
15  & 342.43 & 99.92 & 70.82  & 3600.00 & 342.43 & 342.40 & 0.01   & 339.36  & 342.44  & 0.00 & 100.00 & 3600.00 \\
20  & 230.62 & 8.23  & 96.43  & 3600.00 & 230.62 & 184.95 & 19.80  & 3600.00 & 231.30  & 0.00 & 100.00 & 3600.00 \\
25  & 136.81 & 0.62  & 99.55  & 3600.00 & 135.03 & 84.28  & 37.58  & 3600.00 & 140.37  & 0.00 & 100.00 & 3600.00 \\
30  & 174.38 & 0.00  & 100.00 & 3600.00 & 172.21 & 52.11  & 69.74  & 3600.00 & 199.55  & 0.00 & 100.00 & 3600.00 \\
35  & 136.94 & 0.00  & 100.00 & 3600.00 & 133.81 & 12.34  & 90.78  & 3600.00 & 219.28  & 0.00 & 100.00 & 3600.00 \\
40  & 256.27 & 0.00  & 100.00 & 3600.00 & 257.70 & 53.99  & 79.05  & 3600.00 & 1249.10 & 0.00 & 100.00 & 3600.00 \\
45  & 226.81 & 0.01  & 100.00 & 3600.00 & 228.73 & 20.61  & 90.99  & 3600.00 & 613.29  & 0.00 & 100.00 & 3600.00 \\
50  & 170.45 & 0.00  & 100.00 & 3600.00 & 177.92 & 21.70  & 87.80  & 3600.00 & 708.46  & 0.00 & 100.00 & 3600.00 \\
55  & 205.50 & 0.00  & 100.00 & 3600.00 & 230.96 & 11.95  & 94.83  & 3600.00 & 684.75  & 0.00 & 100.00 & 3600.00 \\
60  & 134.96 & 0.00  & 100.00 & 3600.00 & 772.27 & 1.35   & 99.82  & 3600.00 & 1142.74 & 0.00 & 100.00 & 3600.00 \\
65  & 150.43 & 0.00  & 100.00 & 3600.00 & 176.77 & 0.03   & 99.98  & 3600.00 & 658.52  & 0.00 & 100.00 & 3600.00 \\
70  & 138.49 & 0.00  & 100.00 & 3600.00 & 144.42 & 0.00   & 100.00 & 3600.00 & 596.50  & 0.00 & 100.00 & 3600.00 \\
75  & 140.58 & 0.00  & 100.00 & 3600.00 & 242.28 & 0.00   & 100.00 & 3600.00 & 1021.06 & 0.00 & 100.00 & 3600.00 \\
80  & 169.10 & 0.00  & 100.00 & 3600.00 & 253.35 & 0.00   & 100.00 & 3600.00 & 793.09  & 0.00 & 100.00 & 3600.00 \\
85  & 148.41 & 0.00  & 100.00 & 3600.00 & 320.75 & 0.00   & 100.00 & 3600.00 & 773.68  & 0.00 & 100.00 & 3600.00 \\
90  & 174.45 & 0.00  & 100.00 & 3600.00 & 400.22 & 0.00   & 100.00 & 3600.00 & 835.96  & 0.00 & 100.00 & 3600.00 \\
95  & 177.90 & 0.00  & 100.00 & 3600.00 & 587.42 & 0.00   & 100.00 & 3600.00 & 898.22  & 0.00 & 100.00 & 3600.00 \\
100 & 161.34 & 0.00  & 100.00 & 3600.00 & 819.40 & 0.00   & 100.00 & 3600.00 & 727.02  & 0.00 & 100.00 & 3600.00 \\ \hline
\end{tabular}
\end{table}

\begin{table}[htbp]
\centering
\caption{Results of Exact MICP Formulations}\label{table_MIP_34}
\tiny
\begin{tabular}{c|rrrr|rrrr}
\hline
\multirow{2}{*}{m} & \multicolumn{4}{c|}{\hyperref[model3_WD]{Quantile Formulation}} & \multicolumn{4}{c}{\hyperref[model4_WD]{Aggregate Quantile Formulation}} \\ \cline{2-9}
& \multicolumn{1}{c}{Obj.Val} & \multicolumn{1}{c}{LB} & \multicolumn{1}{c}{Gap (\%)} & \multicolumn{1}{c|}{Time} & \multicolumn{1}{c}{Obj.Val} & \multicolumn{1}{c}{LB} & \multicolumn{1}{c}{Gap (\%)} & \multicolumn{1}{c}{Time} \\ \hline
15  & 342.43 & 342.43 & 0.00   & 0.51    & 342.43 & 342.43 & 0.00  & 0.21    \\
20  & 230.62 & 230.62 & 0.00   & 6.94    & 230.62 & 230.62 & 0.00  & 0.50    \\
25  & 135.03 & 135.03 & 0.00   & 40.80   & 135.03 & 135.03 & 0.00  & 1.50    \\
30  & 172.21 & 172.21 & 0.00   & 356.06  & 172.21 & 172.21 & 0.00  & 3.14    \\
35  & 133.54 & 133.54 & 0.00   & 271.32  & 133.54 & 133.54 & 0.00  & 5.77    \\
40  & 252.42 & 252.42 & 0.00   & 2379.63 & 252.42 & 252.42 & 0.00  & 8.61    \\
45  & 219.17 & 177.31 & 19.10  & 3600.00 & 219.17 & 219.17 & 0.00  & 46.05   \\
50  & 170.16 & 120.52 & 29.17  & 3600.00 & 169.99 & 169.99 & 0.00  & 169.07  \\
55  & 204.76 & 137.83 & 32.69  & 3600.00 & 204.76 & 204.76 & 0.00  & 230.42  \\
60  & ---    & 66.85  & 48.91  & 3600.00 & 130.84 & 130.84 & 0.00  & 1022.89 \\
65  & ---    & 47.18  & 66.90  & 3600.00 & 142.54 & 142.27 & 0.19  & 3600.00 \\
70  & ---    & 30.48  & 77.57  & 3600.00 & 135.92 & 135.91 & 0.01  & 3134.58 \\
75  & ---    & 31.39  & 77.24  & 3600.00 & ---    & 128.64 & 6.72  & 3600.00 \\
80  & ---    & 25.32  & 84.03  & 3600.00 & ---    & 154.58 & 2.51  & 3600.00 \\
85  & ---    & 27.86  & 80.89  & 3600.00 & 157.36 & 143.52 & 8.79  & 3600.00 \\
90  & ---    & 26.65  & 84.47  & 3600.00 & ---    & 165.23 & 3.71  & 3600.00 \\
95  & ---    & 4.57   & 97.37  & 3600.00 & ---    & 162.50 & 6.51  & 3600.00 \\
100 & ---    & 0.00   & 100.00 & 3600.00 & ---    & 124.59 & 12.17 & 3600.00 \\ \hline
\end{tabular}
\end{table}

\begin{table}[htbp]
\centering
\caption{Results of AM Algorithm of \ref{eq_sp_fair}, \hyperref[lower_bound_Wq]{Jensen Bound} and \hyperref[gelbrich_Cholesky]{Gelbrich Bound}}\label{table_AM_Jen_Gel}
\tiny
\begin{tabular}{c|rrr|rrr|rrr}
\hline\multirow{2}{*}{m} & \multicolumn{3}{c|}{AM} & \multicolumn{3}{c|}{\hyperref[lower_bound_Wq]{Jensen Bound}} & \multicolumn{3}{c}{\hyperref[gelbrich_Cholesky]{Gelbrich Bound}} \\ \cline{2-10} 
& \multicolumn{1}{c}{Obj.Val} & \multicolumn{1}{c}{Gap (\%)} & \multicolumn{1}{c|}{Time} & \multicolumn{1}{c}{Obj.Val} & \multicolumn{1}{c}{Gap (\%)} & \multicolumn{1}{c|}{Time} & \multicolumn{1}{c}{Obj.Val} & \multicolumn{1}{c}{Gap (\%)} & \multicolumn{1}{c}{Time} \\ \hline
15  & 342.43 & 0.00  & 0.17 & 207.12 & 39.51 & 0.06 & 222.14 & 35.13 & 0.39 \\
20  & 230.62 & 0.00  & 0.26 & 89.13  & 61.35 & 0.08 & 105.91 & 54.07 & 0.41 \\
25  & 135.03 & 0.00  & 0.37 & 46.03  & 65.91 & 0.09 & 50.41  & 62.67 & 0.48 \\
30  & 172.21 & 0.00  & 0.30 & 109.56 & 36.38 & 0.10 & 110.85 & 35.63 & 0.48 \\
35  & 133.54 & 0.00  & 0.35 & 92.95  & 30.39 & 0.11 & 93.07  & 30.31 & 0.50 \\
40  & 252.42 & 0.00  & 0.23 & 187.66 & 25.65 & 0.13 & 194.65 & 22.89 & 0.54 \\
45  & 219.17 & 0.00  & 0.63 & 138.17 & 36.96 & 0.14 & 176.41 & 19.51 & 0.59 \\
50  & 170.17 & 0.10  & 0.52 & 116.18 & 31.65 & 0.16 & 143.05 & 15.85 & 0.72 \\
55  & 204.76 & 0.00  & 0.62 & 141.53 & 30.88 & 0.16 & 178.35 & 12.90 & 0.65 \\
60  & 130.84 & 0.00  & 0.46 & 69.87  & 46.60 & 0.18 & 105.48 & 19.38 & 1.18 \\
65  & 142.54 & 0.00  & 0.90 & 83.21  & 41.63 & 0.19 & 118.92 & 16.57 & 1.40 \\
70  & 135.92 & 0.00  & 0.51 & 83.84  & 38.31 & 0.20 & 119.37 & 12.18 & 0.81 \\
75  & 137.91 & 6.72  & 0.68 & 85.22  & 38.21 & 0.22 & 119.79 & 13.14 & 0.89 \\
80  & 158.57 & 2.51  & 0.73 & 112.03 & 29.35 & 0.24 & 145.23 & 8.41  & 0.96 \\
85  & 145.77 & 1.54  & 0.81 & 104.78 & 28.12 & 0.24 & 134.34 & 7.84  & 0.99 \\
90  & 171.59 & 3.71  & 0.80 & 126.73 & 26.15 & 0.29 & 160.95 & 6.20  & 1.08 \\
95  & 173.82 & 5.36  & 0.91 & 132.55 & 23.74 & 0.31 & 164.49 & 5.36  & 2.00 \\
100 & 141.85 & 6.59  & 0.92 & 91.76  & 35.32 & 0.32 & 132.50 & 6.59  & 2.27 \\ \hline
\end{tabular}
\end{table}

\begin{table}[htbp]
\centering
\caption{Results of Continuous Relaxation Values of Exact MICP Formulations}\label{table_MIP_Conti}
\tiny
\begin{tabular}{c|rrr|rrr|rrr|rrr}
\hline\multirow{2}{*}{m} & \multicolumn{3}{c|}{\hyperref[WD_q_trans_micpr]{Discretized Formulation}} & \multicolumn{3}{c|}{\hyperref[model2_WD1]{Complementary Formulation}} & \multicolumn{3}{c|}{\hyperref[model3_WD]{Quantile Formulation}} & \multicolumn{3}{c}{\hyperref[model4_WD]{Aggregate Quantile Formulation}} \\ \cline{2-13} 
& \multicolumn{1}{c}{Obj.Val} & \multicolumn{1}{c}{Gap (\%)} & \multicolumn{1}{c|}{Time} & \multicolumn{1}{c}{Obj.Val} & \multicolumn{1}{c}{Gap (\%)} & \multicolumn{1}{c|}{Time} & \multicolumn{1}{c}{Obj.Val} & \multicolumn{1}{c}{Gap (\%)} & \multicolumn{1}{c|}{Time} & \multicolumn{1}{c}{Obj.Val} & \multicolumn{1}{c}{Gap (\%)} & \multicolumn{1}{c}{Time} \\ \hline
15  & 0.20 & 99.94  & 0.24  & 0.00 & 100.00 & 0.38  & 0.34 & 99.90  & 0.23 & 282.42 & 17.52 & 0.16 \\
20  & 0.00 & 100.00 & 0.41  & 0.00 & 100.00 & 0.60  & 0.00 & 100.00 & 0.63 & 147.18 & 36.18 & 0.64 \\
25  & 0.00 & 100.00 & 0.50  & 0.00 & 100.00 & 1.00  & 0.38 & 99.72  & 0.60 & 57.40  & 57.49 & 0.63 \\
30  & 0.00 & 100.00 & 0.77  & 0.00 & 100.00 & 1.05  & 0.00 & 100.00 & 0.49 & 111.36 & 35.33 & 1.07 \\
35  & 0.00 & 100.00 & 1.09  & 0.00 & 100.00 & 5.64  & 0.51 & 99.62  & 0.62 & 98.71  & 26.08 & 0.67 \\
40  & 0.00 & 100.00 & 1.36  & 0.00 & 100.00 & 17.50 & 0.00 & 100.00 & 0.72 & 219.32 & 13.11 & 0.82 \\
45  & 0.00 & 100.00 & 1.73  & 0.00 & 100.00 & 37.80 & 0.96 & 99.56  & 1.22 & 179.74 & 17.99 & 1.03 \\
50  & 0.00 & 100.00 & 2.11  & 0.00 & 100.00 & 5.01  & 0.00 & 100.00 & 1.09 & 134.92 & 20.63 & 1.25 \\
55  & 0.00 & 100.00 & 2.63  & 0.00 & 100.00 & 5.76  & 0.65 & 99.68  & 1.81 & 167.17 & 18.36 & 1.51 \\
60  & 0.00 & 100.00 & 3.08  & 0.00 & 100.00 & 6.60  & 0.00 & 100.00 & 1.54 & 92.30  & 29.45 & 1.77 \\
65  & 0.00 & 100.00 & 3.68  & 0.00 & 100.00 & 8.72  & 0.44 & 99.69  & 2.00 & 101.41 & 28.85 & 2.24 \\
70  & 0.00 & 100.00 & 4.76  & 0.00 & 100.00 & 10.27 & 0.00 & 100.00 & 2.08 & 103.16 & 24.10 & 2.43 \\
75  & 0.00 & 100.00 & 5.46  & 0.00 & 100.00 & 13.96 & 0.32 & 99.77  & 2.71 & 99.28  & 28.01 & 3.09 \\
80  & 0.00 & 100.00 & 6.19  & 0.00 & 100.00 & 11.26 & 0.00 & 100.00 & 2.79 & 120.75 & 23.85 & 3.40 \\
85  & 0.00 & 100.00 & 7.25  & 0.00 & 100.00 & 13.57 & 0.26 & 99.82  & 3.49 & 115.93 & 20.47 & 4.08 \\
90  & 0.00 & 100.00 & 8.06  & 0.00 & 100.00 & 14.48 & 0.00 & 100.00 & 4.12 & 136.94 & 20.19 & 4.65 \\
95  & 0.00 & 100.00 & 9.27  & 0.00 & 100.00 & 16.51 & 0.21 & 99.88  & 4.99 & 140.77 & 19.01 & 5.48 \\
100 & 0.00 & 100.00 & 10.59 & 0.00 & 100.00 & 18.72 & 0.00 & 100.00 & 4.35 & 102.23 & 27.93 & 5.66 \\ \hline
\end{tabular}
\end{table}

\subsubsection{Comparison with state-of-the-art}
In the second experiment, we compare our methods against two fair regression methods from the literature using real data. In this experiment, the cost function is set to the mean squared error (MSE). We solve \ref{eq_sp_fair} using its AM algorithm in Section \ref{sec_AM}, solve the Jensen bound by solving \eqref{lower_bound_Wq}, and solve the \hyperref[gelbrich_Cholesky]{Gelbrich bound} using its AM algorithm in Section \ref{sec_gelbrich}.
The first approach that we compare is Berk \citep{berk2017convex}. 
Their work proposed a convex optimization method to incorporate group and individual fairness for fair regression. We compare with their group fairness model for demographic parity. The second approach is Agarwal \citep{agarwal2019fair}, a reduction-based fair regression algorithm that uses the Kolmogorov–Smirnov distance to measure demographic parity.
We test the performance of different approaches using criminological and educational datasets.
The \textit{Communities and Crime}  dataset contains socio-economic, law enforcement, and crime data of different communities in the US. The goal is to predict the number of violent crimes per $100{,}000$ of the population with race (black versus non-black) as the sensitive attribute. It contains $1{,}994$ samples characterized by $127$ features. We create the sensitive attribute by thresholding the percentage of the black population following \cite{calders2013controlling}. 
The \textit{Law School} \citep{wightman1998lsac} dataset consists of student records from the Law School Admission Council (LSAC) National Longitudinal Bar Passage Study. The goal is to predict a student's GPA with race (white versus non-white) as the sensitive attribute. The dataset contains $20{,}649$ samples characterized by $12$ features.

For both datasets, we split the data into $70\%$ for training and $30\%$ for testing. We repeat this procedure $10$ times and report the average performance. We evaluate the different methods based on the trade-off between MSE and fairness scores of Wasserstein fairness and Kolmogorov–Smirnov fairness measures, where we compute $\WD_2(\bfx)$ using \eqref{eq_lem_inv_cdf_WD1} and $\KSD(\bfx)$ using \eqref{dfn_KSD} for each method. Note that we only use \ref{eq_sp_fair} to solve the Wasserstein fairness measure and plug in its solution to compute the Kolmogorov–Smirnov fairness measure.
The hyperparameters of each method are chosen as follows. For \ref{eq_sp_fair}, the \hyperref[lower_bound_Wq]{Jensen bound}, and the \hyperref[gelbrich_Cholesky]{Gelbrich bound}, we set the inefficiency level parameter $\epsilon\in\{0.05,0.1,\dots,1\}$ for the \textit{Communities and Crime} dataset and $\epsilon\in\{0.01,0.02,\dots,0.2\}$ for the \textit{Law School} dataset. For Agarwal \citep{agarwal2019fair}, we set their unfairness level parameter $\epsilon\in\{0.015,0.03,\dots,0.3\}$ for \textit{Communities and Crime} and $\epsilon\in\{0.035,0.07,\dots,0.7\}$ for \textit{Law School}. For Berk \citep{berk2017convex}, we set their unfairness penalty parameter $\lambda\in\{0.5,1,\dots,10\}$ for \textit{Communities and Crime} and $\lambda\in\{0.2,0.4,\dots,4\}$ for \textit{Law School}.

\begin{figure}[htbp]
\begin{subfigure}[Training of \textit{C \& C}]{{
\centering\includegraphics[width=0.23\textwidth]{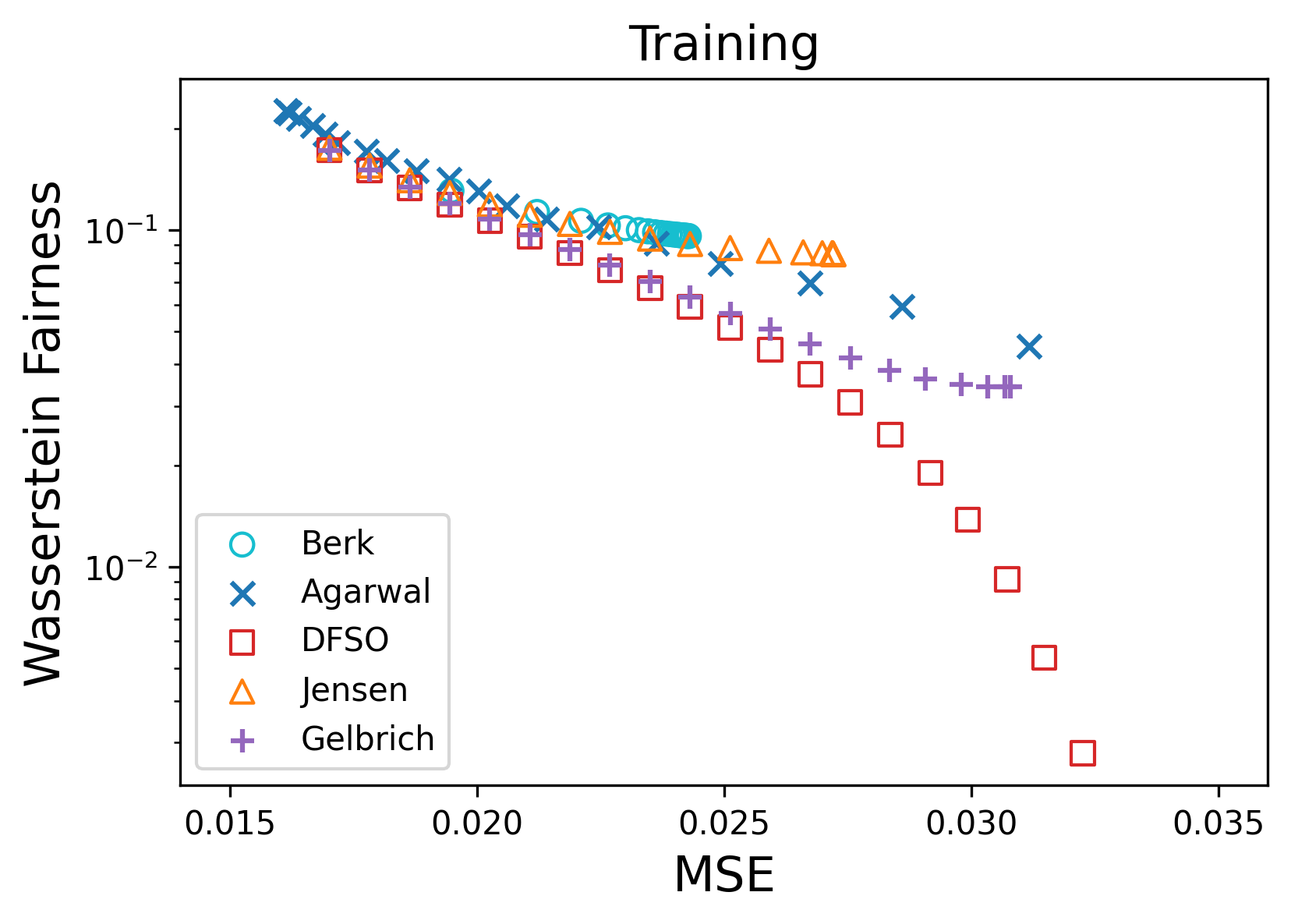}}\label{fig_crime_train_WD}}
\end{subfigure}
\hfill
\begin{subfigure}[Training of \textit{C \& C}]{{
\centering\includegraphics[width=0.23\textwidth]{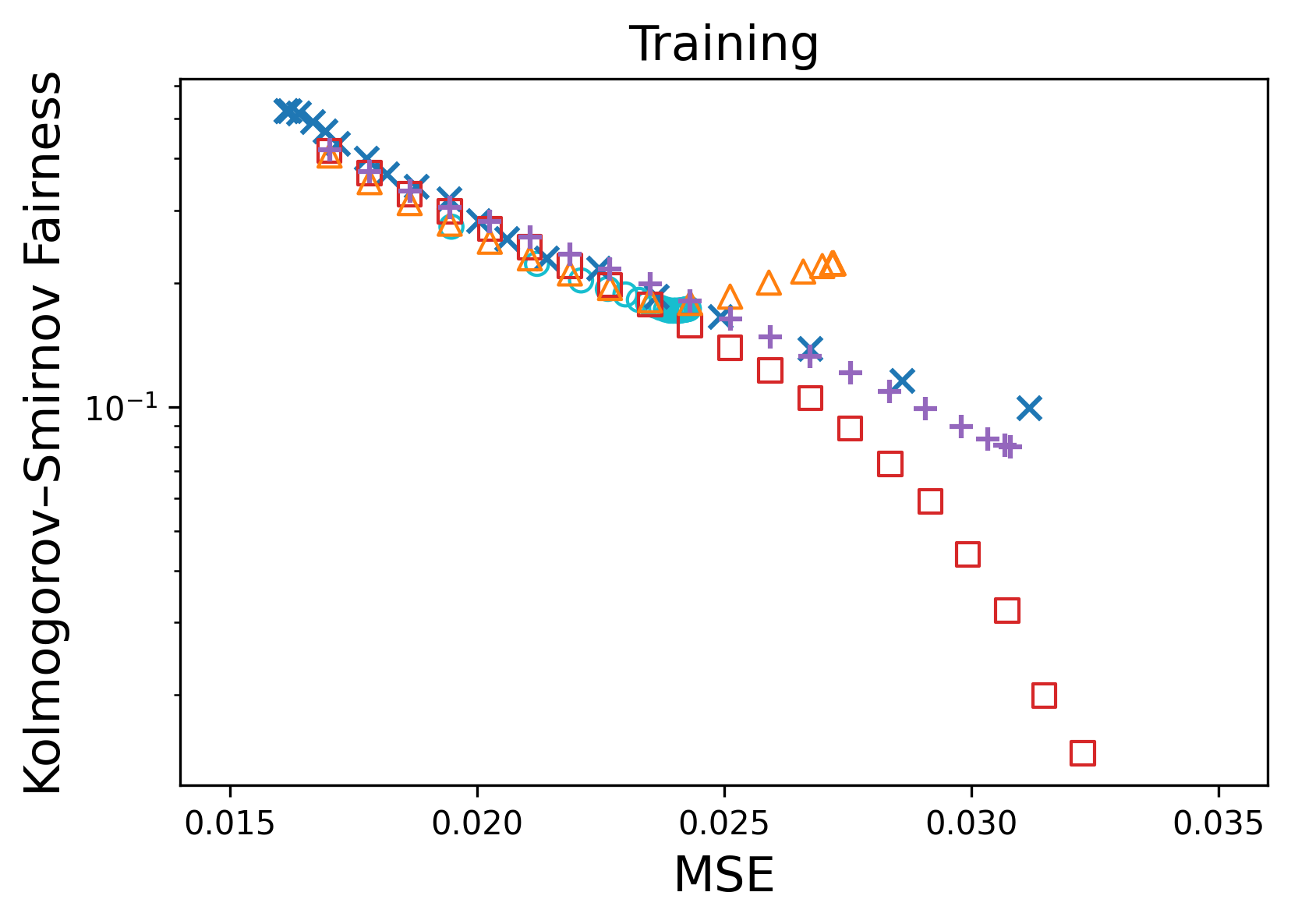}}\label{fig_crime_train_KSD}}
\end{subfigure}
\hfill
\begin{subfigure}[Training of \textit{Law School}]{{
\centering\includegraphics[width=0.23\textwidth]{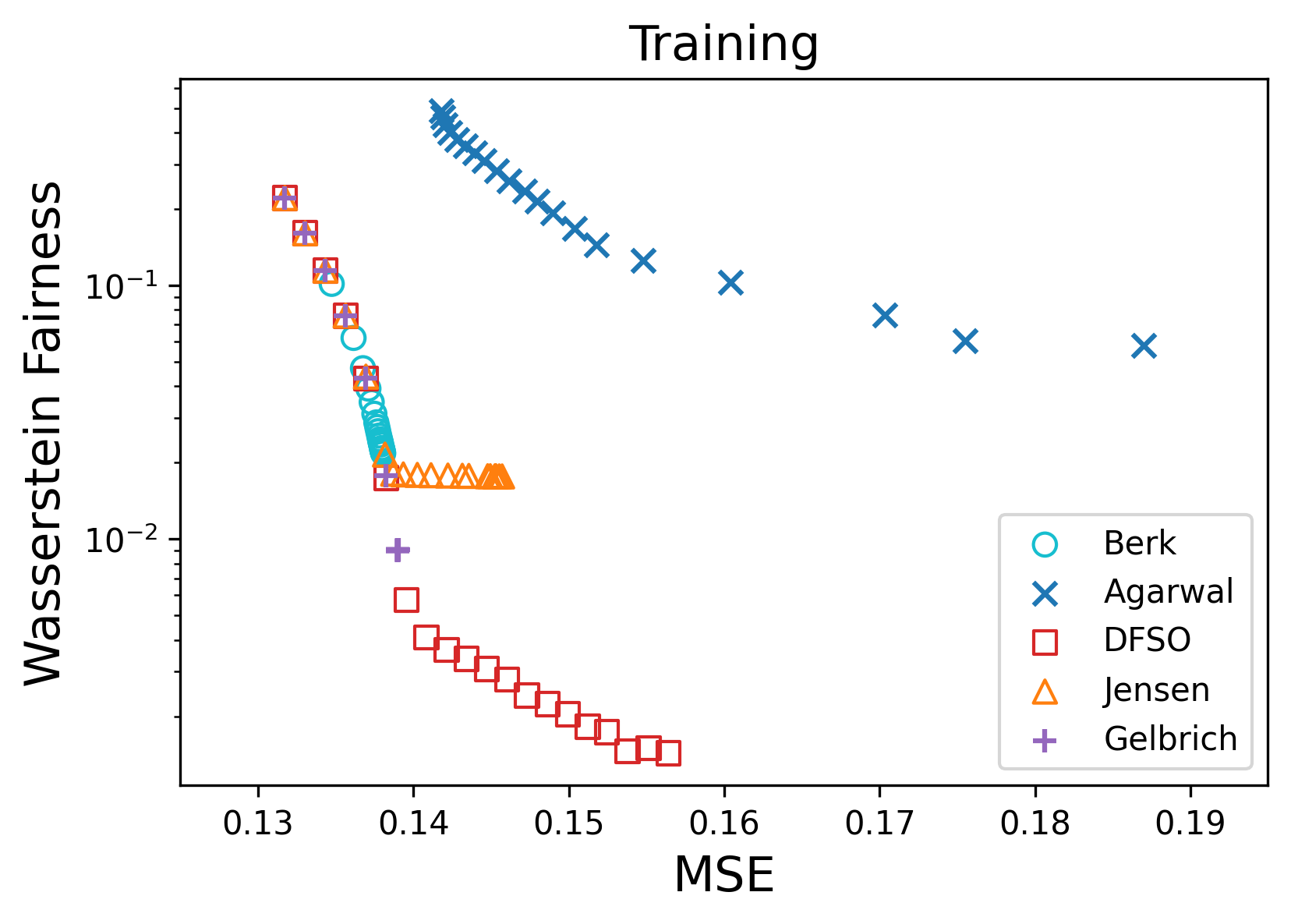}}\label{fig_law_train_WD}}
\end{subfigure}
\hfill
\begin{subfigure}[Training of \textit{Law School}]{{
\centering\includegraphics[width=0.23\textwidth]{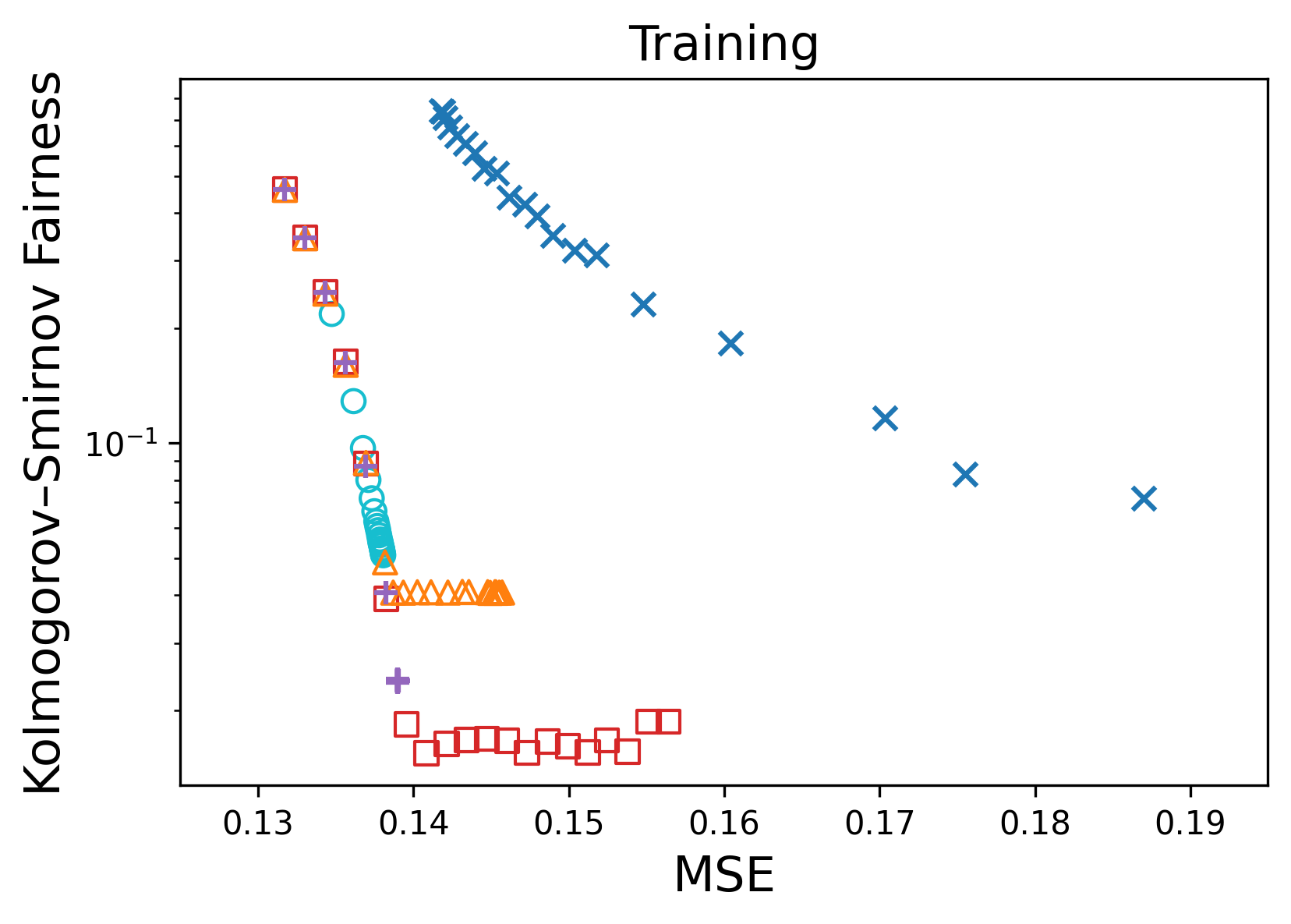}}\label{fig_law_train_KSD}}
\end{subfigure}
\hfill
\begin{subfigure}[Testing of \textit{C \& C}]{{
\centering\includegraphics[width=0.23\textwidth]{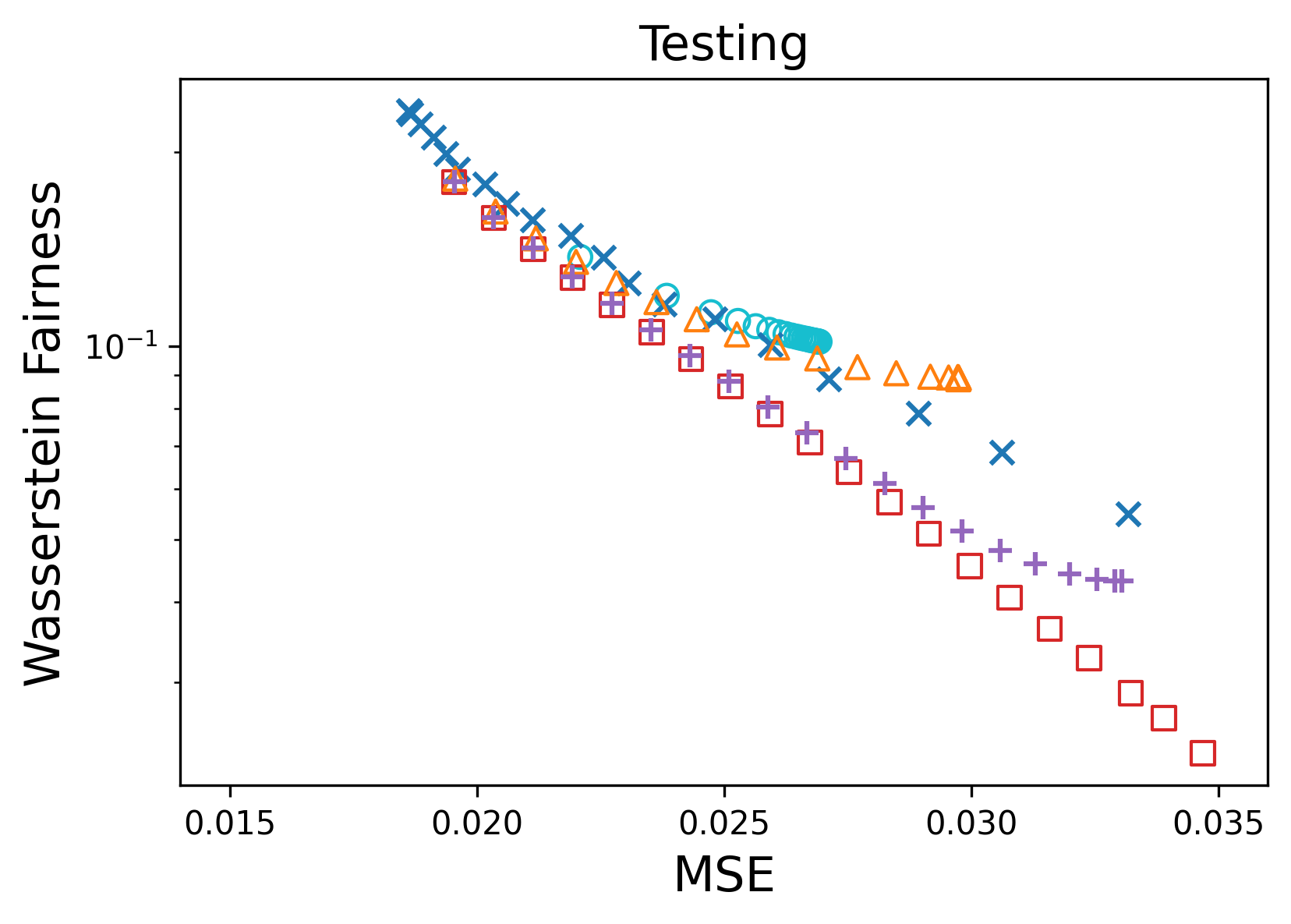}}\label{fig_crime_test_WD}}
\end{subfigure}
\hfill
\begin{subfigure}[Testing of \textit{C \& C}]{{
\centering\includegraphics[width=0.24\textwidth]{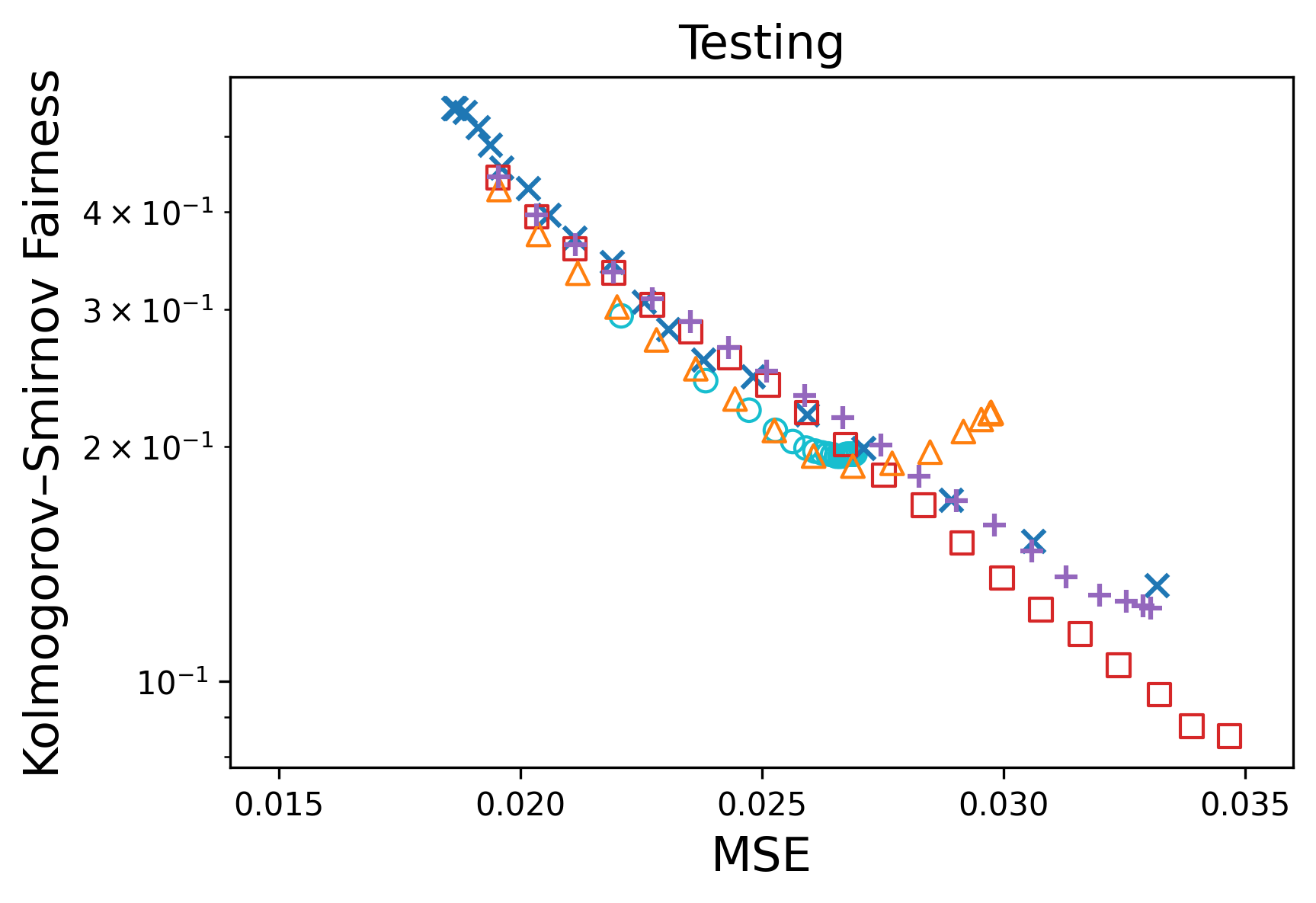}}\label{fig_crime_test_KSD}}
\end{subfigure}
\hfill
\begin{subfigure}[Testing of \textit{Law School}]{{
\centering\includegraphics[width=0.23\textwidth]{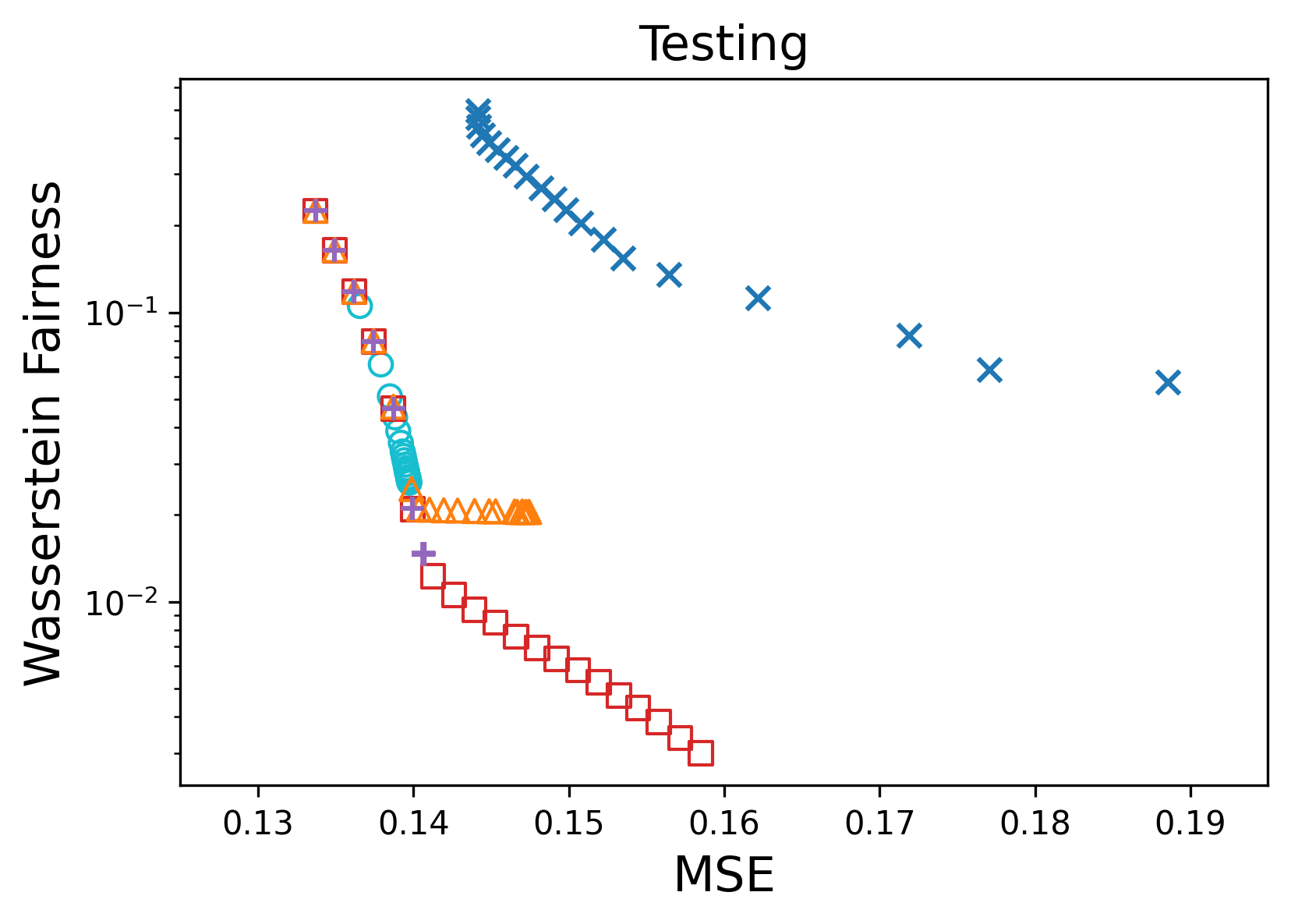}}\label{fig_law_test_WD}}
\end{subfigure}
\hfill
\begin{subfigure}[Testing of \textit{Law School}]{{
\centering\includegraphics[width=0.23\textwidth]{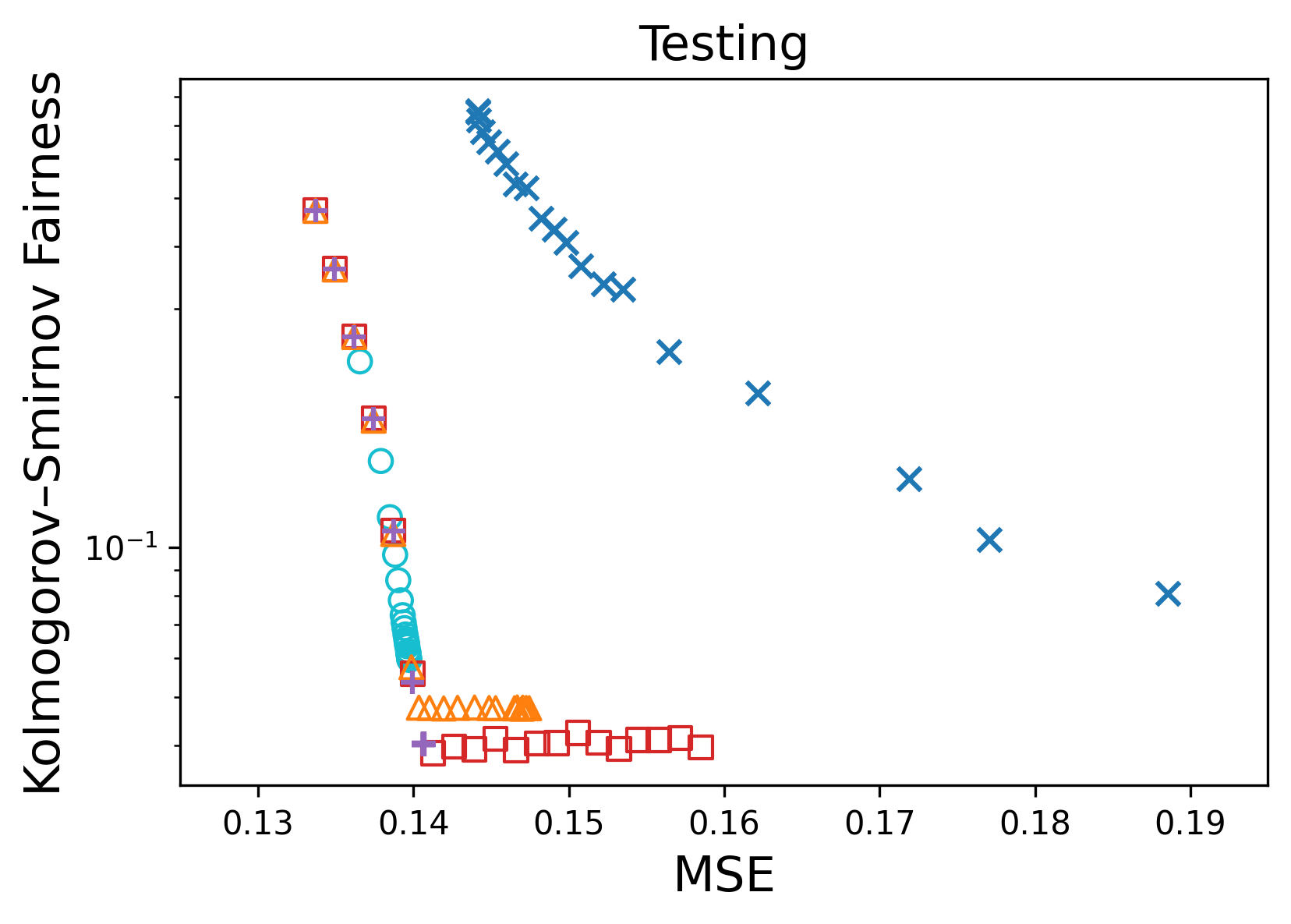}}\label{fig_law_test_KSD}}
\end{subfigure}
\caption{Fairness vs MSE for Fair Regression. The Wasserstein fairness versus MSE are shown in (a), (c), (e), (g), and Kolmogorov–Smirnov fairness versus MSE are shown in (b), (d), (f), (h). All the training and testing results are averaged over $10$ replications.}\label{fig_regression_crime_and_law}
\end{figure}

Figure \ref{fig_regression_crime_and_law} presents the trade-off between fairness scores and MSE on the two datasets described above. We observe that \ref{eq_sp_fair} consistently outperforms other methods in training and testing in view of the Wasserstein fairness measure, where it can reduce the unfairness level to significantly small values (i.e., nearly zero) with a relatively small increase in MSE. \ref{eq_sp_fair} also performs well in terms of the Kolmogorov–Smirnov fairness measure and attains small fairness scores. Two lower bounds (i.e., the \hyperref[lower_bound_Wq]{Jensen} and \hyperref[gelbrich_Cholesky]{Gelbrich} bounds) also provide good solution quality for both Wasserstein and Kolmogorov–Smirnov fairness measures. Notably, the \hyperref[lower_bound_Wq]{Jensen bound} has similar solutions as Berk \citep{berk2017convex}, and the \hyperref[gelbrich_Cholesky]{Gelbrich bound} is competitive with Agarwal \citep{agarwal2019fair}. We observe that the \hyperref[gelbrich_Cholesky]{Gelbrich bound} tends to be fairer than the \hyperref[lower_bound_Wq]{Jensen bound}. \ref{eq_sp_fair} consistently provides the best Wasserstein fairness score for the \textit{Communities and Crime} and \textit{Law School} datasets. In terms of the Kolmogorov–Smirnov fairness measure, Berk \citep{berk2017convex} and the \hyperref[lower_bound_Wq]{Jensen bound} is effective when MSE is small. Note that they cannot further improve fairness given large inefficiency level parameters or allow large MSE in the experiment. On the contrary, \ref{eq_sp_fair} and the \hyperref[gelbrich_Cholesky]{Gelbrich bound} have the capacity to improve Kolmogorov–Smirnov fairness with large MSE. It is evident that \ref{eq_sp_fair} is capable of effectively addressing the unfairness issues in fair regression problems.

\subsection{Fair Allocation of Scarce Medical Resources}

During public health emergencies such as the influenza pandemic and COVID-19, optimal allocation of scarce medical resources (e.g., therapeutics and vaccines) is a crucial yet challenging task (see, e.g., \citealt{sun2023distributionally}, \citealt{shehadeh2023equity}). \ref{eq_sp_fair} can be adapted to allocate scarce medical resources in a distributionally fair way. In this experiment, we study the fair allocation of COVID-19 vaccine across $m$ counties in Georgia. Given the total amount of vaccines $T\in\Ze_{+}$, counties' population sizes $\bm{p}\in\Ze_{+}^{m}$, and thresholds $\bm{l}\leq\bm{u}\in(0,1]^{m}$, we consider a vaccine allocation problem
\begin{align}\label{model_vaccine}
V^*=\max_{\bfx} \left\{ \sqrt[m]{\prod_{i\in[m]} x_i}:
 \sum_{i\in[m]} p_i x_i \leq T, 
l_i \leq x_i \leq u_i, \forall i\in[m]\right\},
\end{align}
where the coverage rate $x_i\in[0,1]$ is defined as the ratio of the number of allocated vaccines to the population size in each county $i\in [m]$, and the benefit function is $Q(\bfx,\bm{\xi}_i) = x_i$. We remark that the efficiency function in \eqref{model_vaccine} follows the conventional proportional fairness, which seeks to maximize the product of each individual county's utility.
The fair allocation approach that we are comparing is the max-min fairness at the group level, defined as 
\begin{equation}\label{model_vaccine_maximin}
\max_{\bfx} \min_{a\in A} \left\{\sqrt[m_a]{\prod_{i\in[m_a]} x_i}: \sqrt[m]{\prod_{i\in[m]} x_i} \geq (1-\epsilon)V^*,\sum_{i\in[m]} p_i x_i \leq T, 
l_i \leq x_i \leq u_i, \forall i\in[m] \right\}.
\end{equation}

In the experiment, we compare \ref{eq_sp_fair} against the Max-Min formulation \eqref{model_vaccine_maximin} using Georgia (GA) population data from the U.S. Census Bureau. We choose the utility function to be $f(\bfx,\bm{\xi}_i) = x_i$ for each county $i\in [m]$. The dataset includes each county's population size $p_i$ and the size of the population aged 65 years and over, denoted by $s_i$. 
We let $A=\{\text{urban, rural}\}$ in order to study the fairness among urban counties (where the population size is at least $50{,}000$) and rural counties (where the population size is less than $50{,}000$). We assume the total amount of vaccine is $20\%$ of the total population. That is, we have $T=0.2\sum_{i\in[m]}p_i$. Since older people are more vulnerable to COVID-19, we select the minimum and maximum vaccine coverage rates as $l_i=0.8 s_i T/\sum_{i\in[m]}s_i$ and $u_i=2 s_i T/\sum_{i\in[m]}s_i$, respectively. We also set the inefficiency level parameter to $\epsilon=\{0.1, 0.2, 0.267, 0.3\}$. We choose type $q=2$ Wasserstein fairness and solve \ref{eq_sp_fair} using its AM algorithm in Section \ref{sec_AM}. Since the solution of Max-Min \eqref{model_vaccine_maximin} remains unchanged when $\epsilon \geq 0.267$, we only display its results for $\epsilon=\{0.1, 0.2, 0.267\}$.

\begin{figure}[htbp]
\begin{subfigure}[Max-Min with $\epsilon=0.1$]{{
\centering\includegraphics[width=0.23\textwidth]{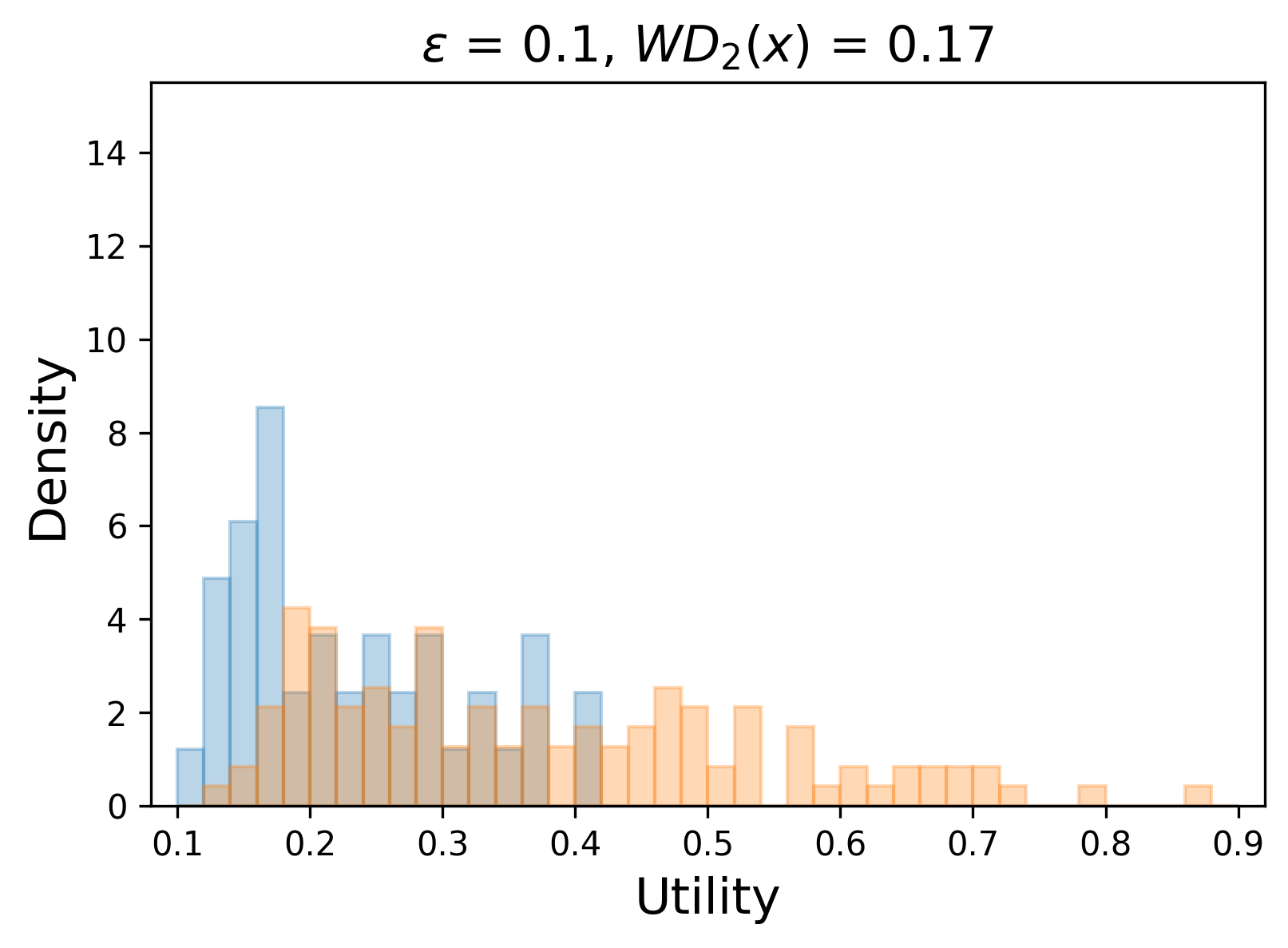}}\label{fig_vaccine_MaxMin_eps10}}
\end{subfigure}
\hfill
\begin{subfigure}[Max-Min with $\epsilon=0.2$]{{
\centering\includegraphics[width=0.23\textwidth]{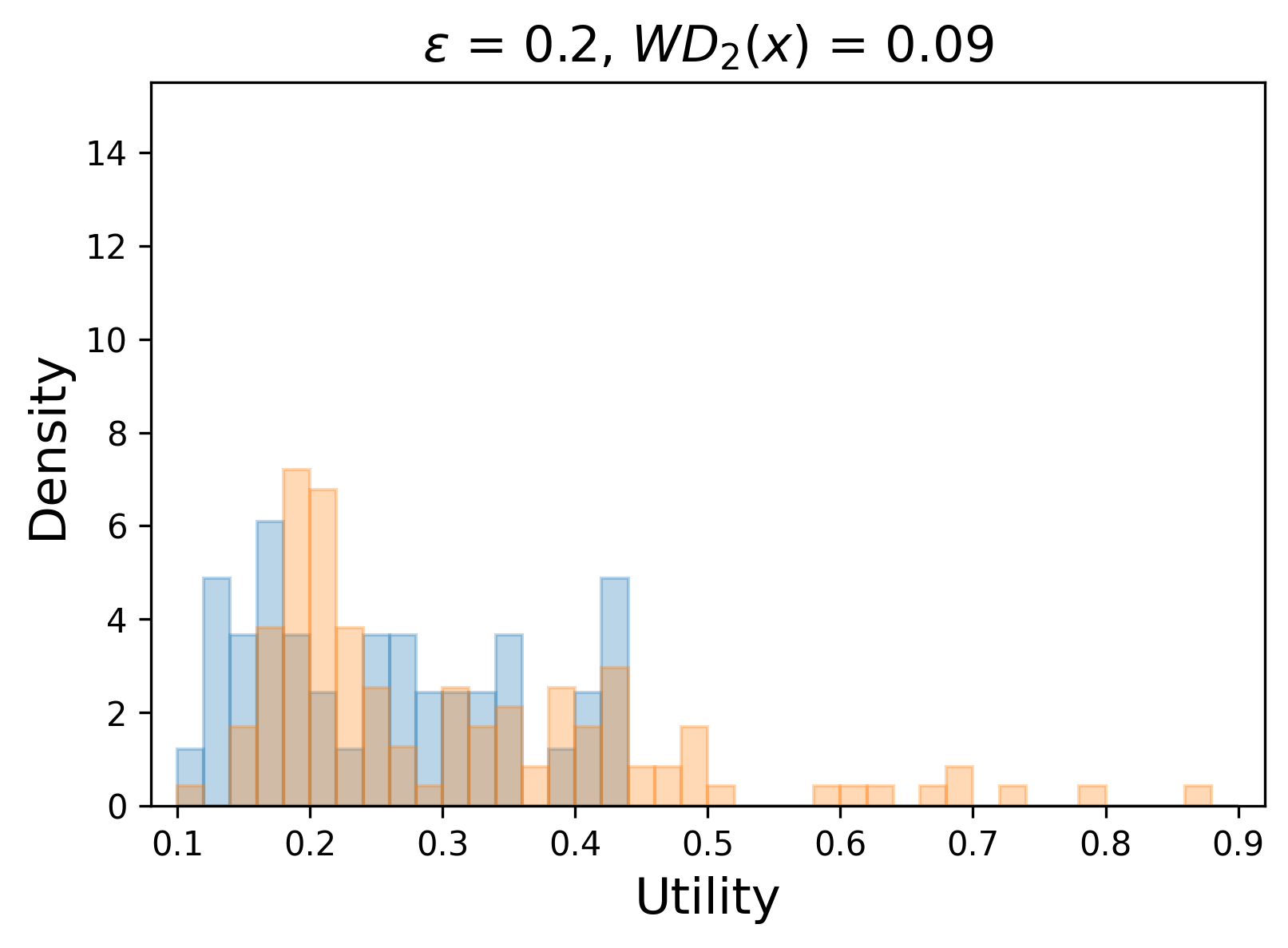}}\label{fig_vaccine_MaxMin_eps20}}
\end{subfigure}
\hfill
\begin{subfigure}[Max-Min with $\epsilon=0.267$]{{
\centering\includegraphics[width=0.23\textwidth]{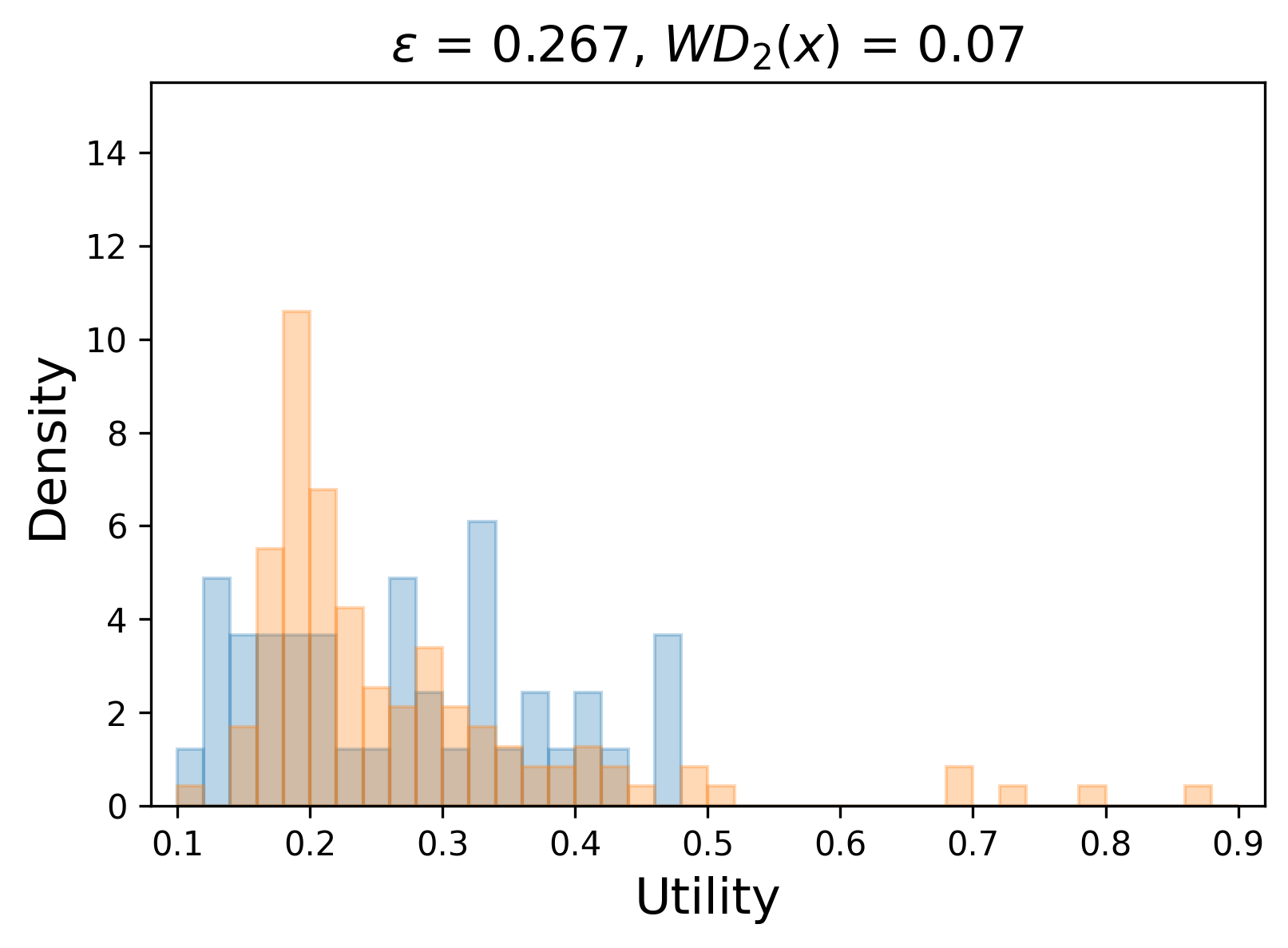}}\label{fig_vaccine_MaxMin_eps26}}
\end{subfigure}
\hfill
\begin{subfigure}[Vanilla]{{
\centering\includegraphics[width=0.24\textwidth]{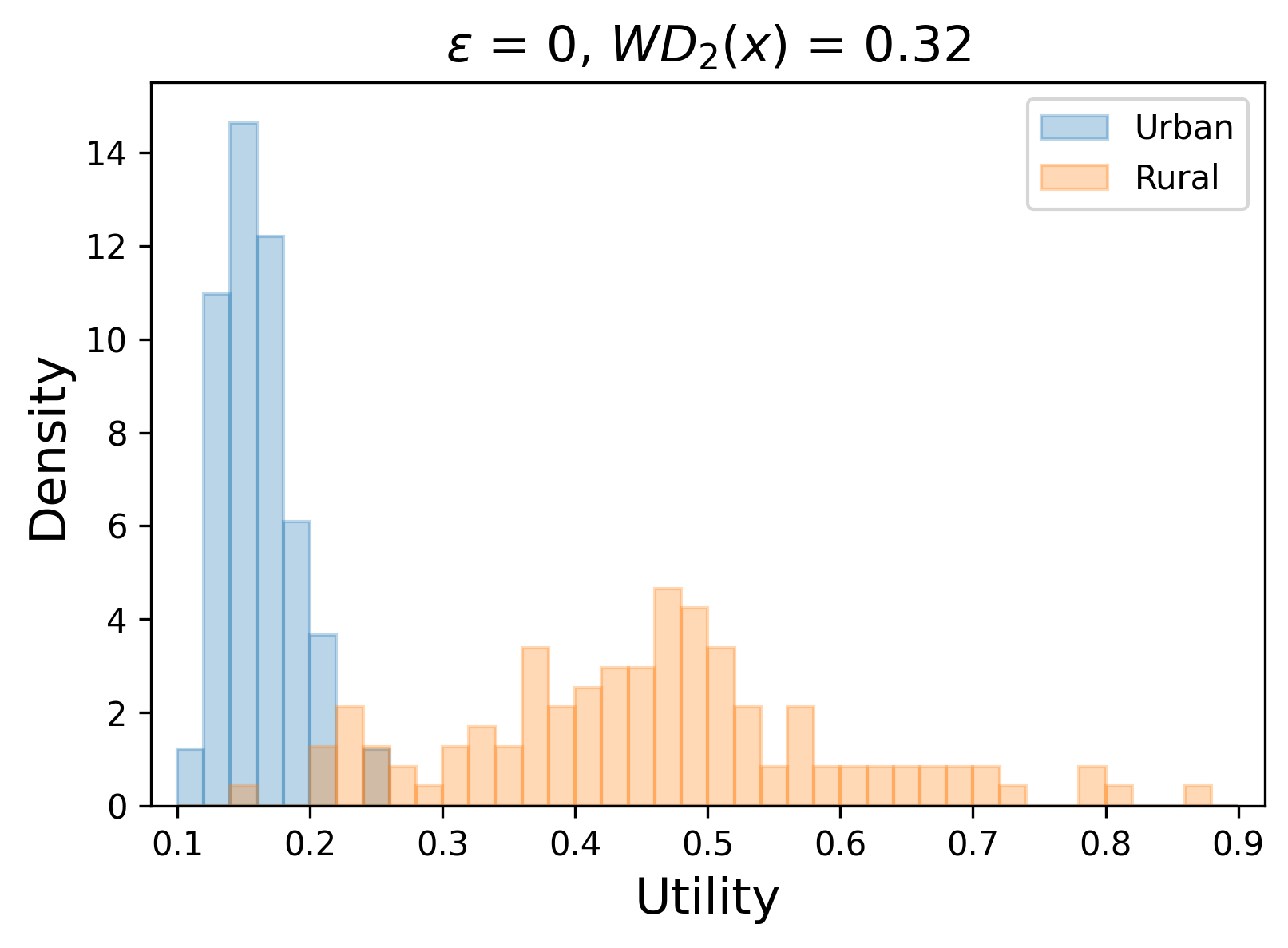}}\label{fig_vaccine_Vanilla_eps0}}
\end{subfigure}
\hfill
\begin{subfigure}[DFSO with $\epsilon=0.1$]{{
\centering\includegraphics[width=0.23\textwidth]{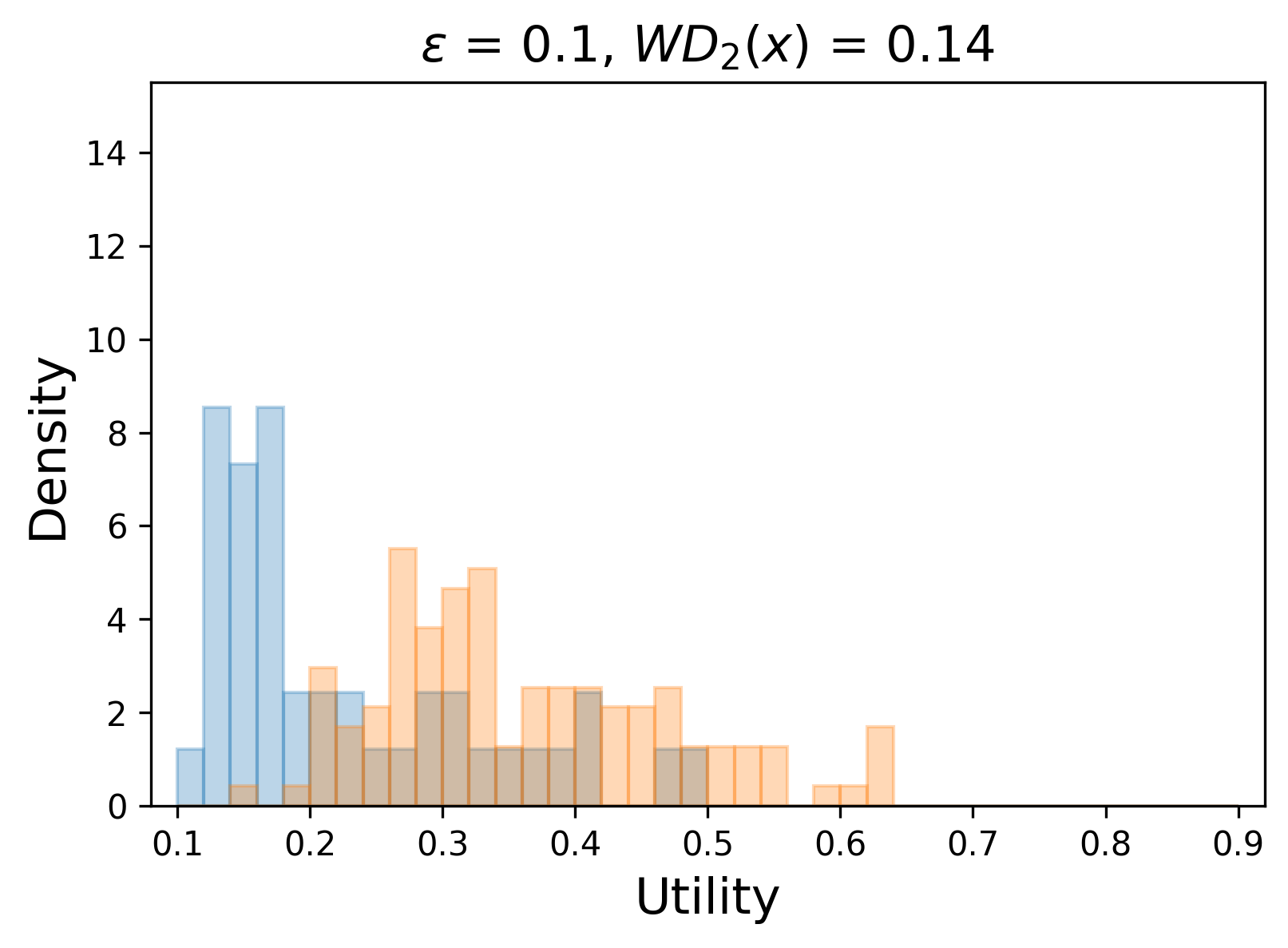}}\label{fig_vaccine_DFSO_eps10}}
\end{subfigure}
\hfill
\begin{subfigure}[DFSO with $\epsilon=0.2$]{{
\centering\includegraphics[width=0.23\textwidth]{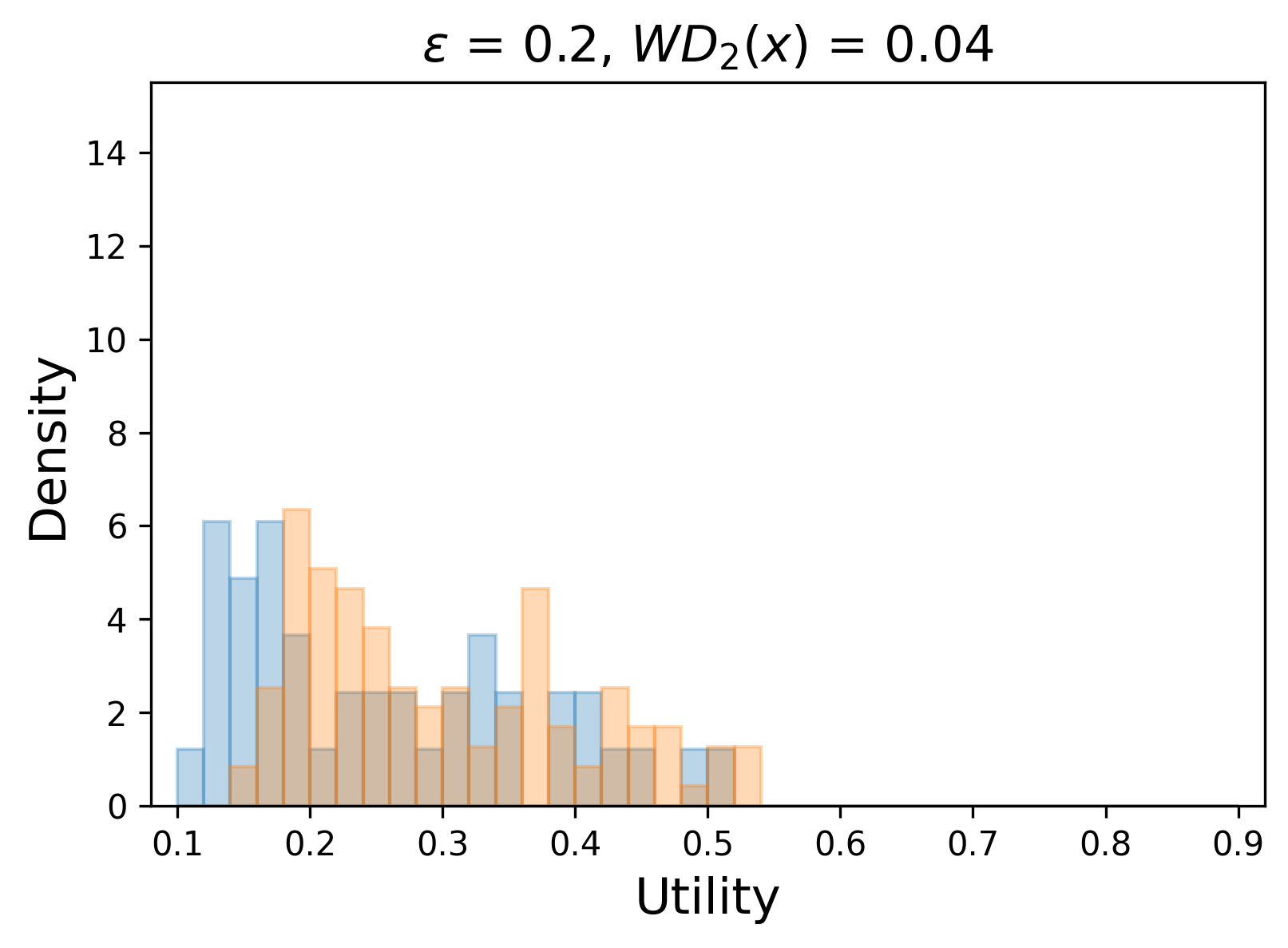}}\label{fig_vaccine_DFSO_eps20}}
\end{subfigure}
\hfill
\begin{subfigure}[DFSO with $\epsilon=0.267$]{{
\centering\includegraphics[width=0.23\textwidth]{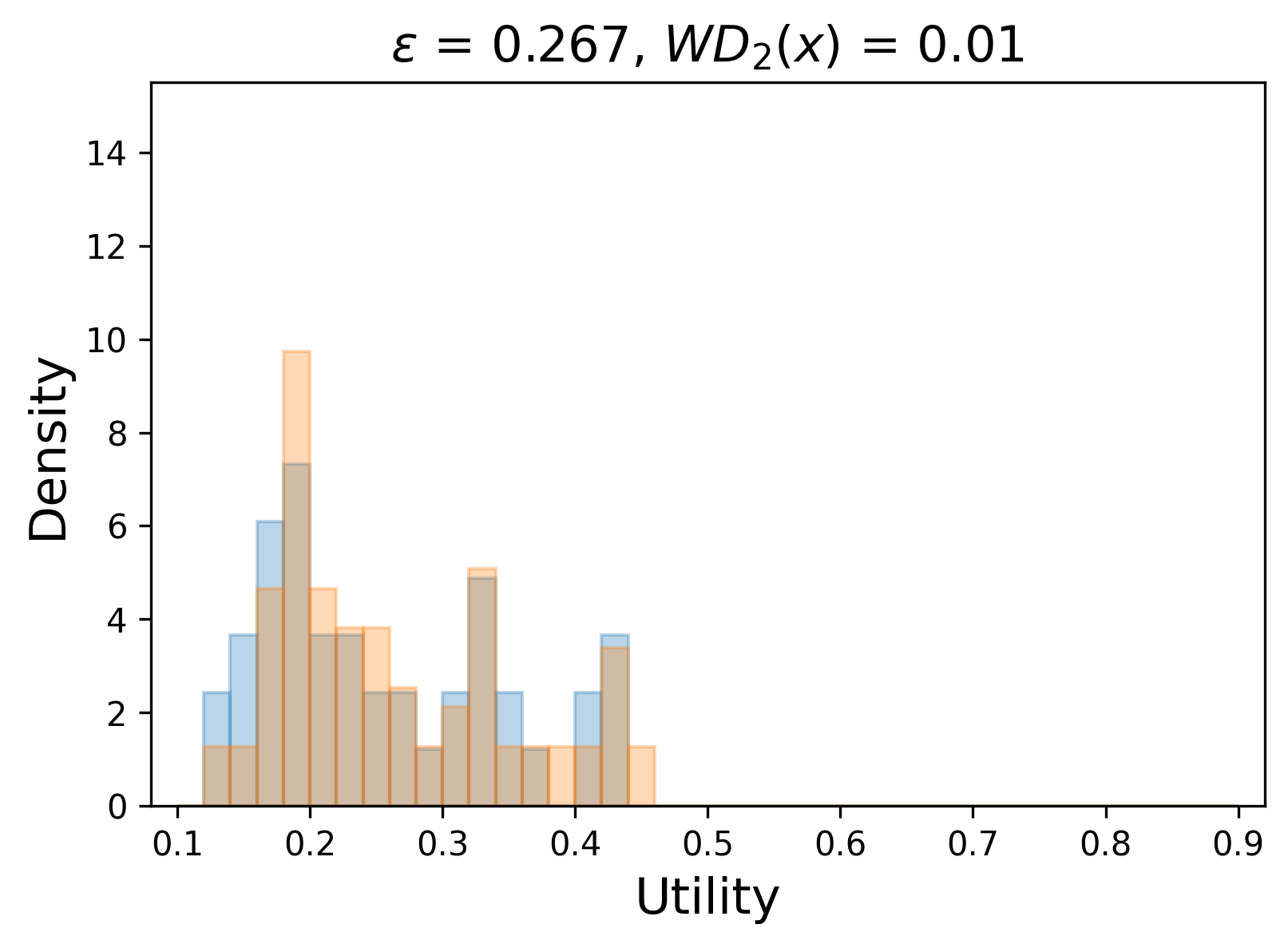}}\label{fig_vaccine_DFSO_eps26}}
\end{subfigure}
\hfill
\begin{subfigure}[DFSO with $\epsilon=0.3$]{{
\centering\includegraphics[width=0.23\textwidth]{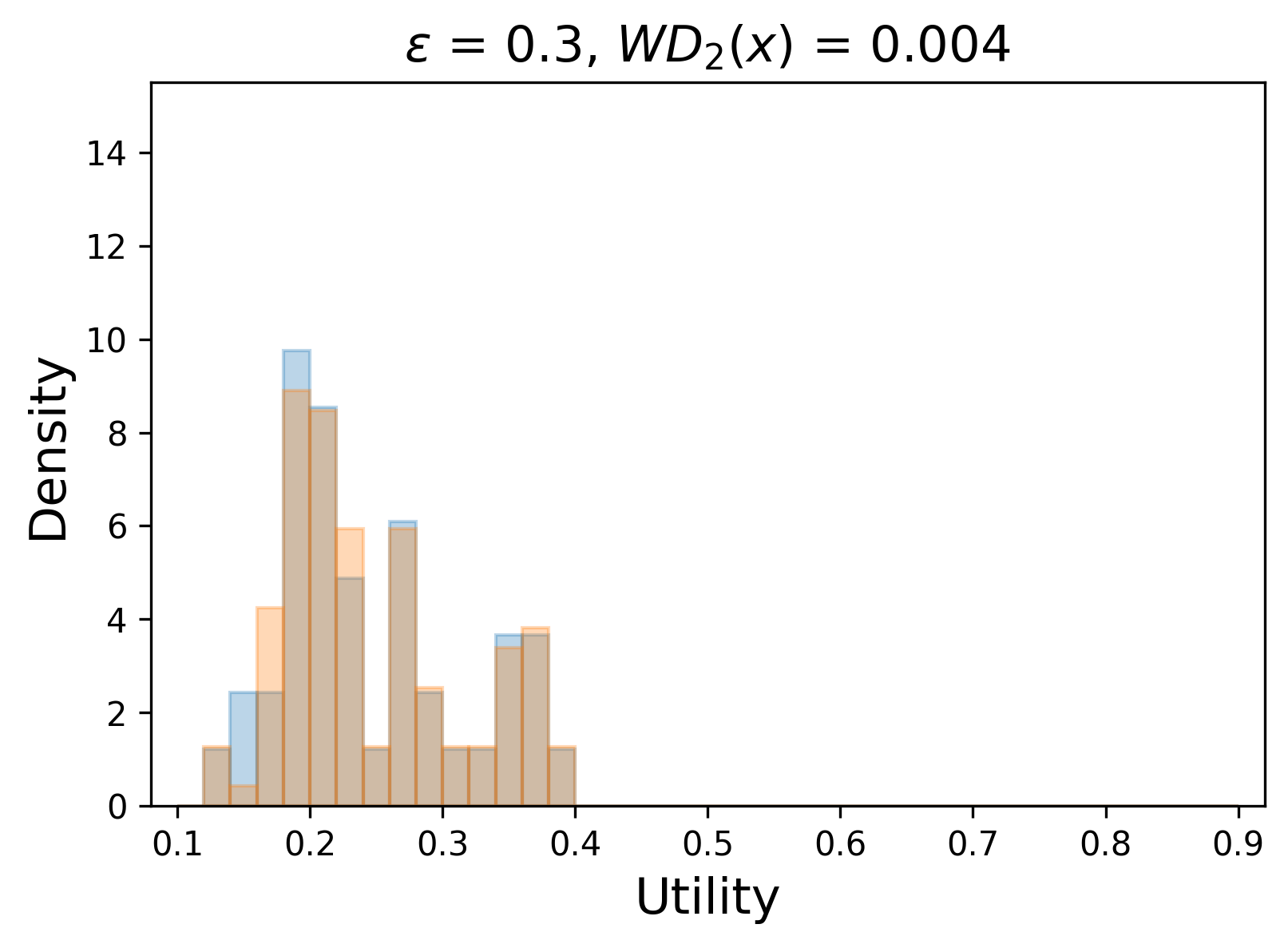}}\label{fig_vaccine_DFSO_eps30}}
\end{subfigure}
\caption{Histograms of Utility for Fair Allocation of COVID-19 Vaccine in GA}\label{fig_vaccine_maxmin_and_dfso}
\end{figure}

Figure \ref{fig_vaccine_maxmin_and_dfso} shows the histograms of utility for fair allocation of COVID-19 vaccine in GA. It can be observed that both methods can reduce the disparities of utilities among urban and rural counties, while \ref{eq_sp_fair} always has a smaller Wasserstein fairness score than Max-Min \eqref{model_vaccine_maximin} given the same inefficiency level. The solution of Max-Min \eqref{model_vaccine_maximin} remains unchanged when $\epsilon \geq 0.267$, thus Figure \ref{fig_vaccine_MaxMin_eps26} shows the fairest solution that Max-Min \eqref{model_vaccine_maximin} can provide. \ref{eq_sp_fair} can achieve a Wasserstein fairness score that is nearly zero, effectively resolving the distributional disparities. We observe that Max-Min \eqref{model_vaccine_maximin} is not sufficient to eliminate the disparity between two distributions compared to its counterpart \ref{eq_sp_fair}. This demonstrates that the proposed \ref{eq_sp_fair} can effectively address distributional fairness while achieving relatively high efficiency.

\section{Conclusion}\label{sec_conclusion}

This paper studies Distributionally Fair Stochastic Optimization (DFSO), where we employ the Wasserstein distance to measure group fairness. We propose exact mixed-integer convex programming formulations for DFSO. By exploring the properties of the Wasserstein fairness measure, we develop an efficient alternating minimization (AM) solution method and two strong lower bounds. Our numerical study shows that the proposed exact methods can solve medium-sized fair learning problems efficiently, while the proposed AM method and lower bounds work efficiently for large-scale fair optimization and learning problems. The convergence rate and solution quality of AM methods are interesting open questions. Stronger lower bounds for the general Wasserstein fairness measure are also interesting to explore in the future. Another future study is properly incorporating distributional robustness into DFSO when the individual data are noise-related.

\bibliography{reference}

\newpage
\titleformat{\section}{\large\bfseries}{\appendixname~\thesection .}{0.5em}{}
\begin{appendices}

\section{Two Additional Exact Formulations and An Equivalent MICP-R Formulation for $\KSD($\boldmath{$x$}\unboldmath{$)$}}\label{sec_reform_add}
\subsection{Discretized Formulation: Discretizing the Transportation Decisions}\label{sec_M1}
In this subsection, we develop an MICP-R formulation for the Wasserstein fairness measure set $\F_q$ by observing that the balanced transportation polytope can be integral given that the supplies and demands are both integers.
To this end, we recast the set  as
\begin{equation}\label{set_F_model1}
\F_{q}=\left\{(\bfx,\nu)\in\mathcal X\times \Re_+: \min_{{\bm{\pi}}_{a\bar a}\in {\Pi}_{a\bar a}}\left\{\sum_{i\in C_a}\sum_{j\in C_{\bar a}} \pi_{ija\bar a}\left|f(\bfx,{\bm\xi}_i)-f(\bfx,{\bm\xi}_j)\right|^q \right\}\leq \nu,\forall a<\bar a \in A\right\}
\end{equation}
where for each $a<\bar a\in A$, the transportation feasible set is given by
\begin{equation*}
\Pi_{a\bar a} = \left\{\bm{\pi}_{a\bar a}\in\Re_{+}^{m_a\times m_{\bar a}}: \sum_{i\in C_a} \pi_{ija\bar a}=\frac{1}{m_{\bar a}}, \forall j\in C_{\bar a}, \sum_{j\in C_{\bar a}} \pi_{ija\bar a}=\frac{1}{m_a}, \forall i\in C_a \right\}. 
\end{equation*}

Observe that, for each $a<\bar a\in A$, the constraint in \eqref{set_F_model1} is satisfied if and only if there exists ${\bm{\pi}}_{a\bar a}\in {\Pi}_{a\bar a}$ such that $\sum_{i\in C_a}\sum_{j\in C_{\bar a}} \pi_{ija\bar a}\left|f(\bfx,{\bm\xi}_i)-f(\bfx,{\bm\xi}_j)\right|^q \leq \nu$. Hence, the resulting set has nonconvex terms in $\bm{\pi}_{a\bar a}$ and $\bm x$ that complicate the formulation. 
\begin{restatable}{theorem}{thmmicpr}\textbf{(Discretized Formulation)}\label{thm_micpr}
Suppose that the set $X_{i}=\{(\bfx, \bar{w}_{i})\in \mathcal X\times\Re: f(\bfx,{\bm\xi}_i)=\bar{w}_{i}\}$ is MICP-R for each $i\in [m]$ and $M_{i}\geq \max_{\bm{x}\in \mathcal X,\eqref{eq_tol}} |f(\bfx,{\bm\xi}_i)|$ for each $i\in [m]$. Then $\F_q$ is equivalent to the  MICP-R set 
\begin{equation}\label{WD_q_trans_micpr}
\F_{q}=\left\{(\bfx,\nu)\in\mathcal X\times \Re_+: \begin{aligned}
&{\sum_{i\in C_a}\sum_{j\in C_{\bar a}}\sum_{k\in[\bar{\Omega}_{a\bar a}]} \frac{2^{k-1}\hat{w}_{ijka\bar a}^q}{m_am_{\bar a}}}\leq \nu, \bm{z}_{a\bar a}\in \Gamma_{a\bar a},\forall a<\bar a\in A, \\
&(\bfx, \bar{w}_{i})\in X_{i}, \forall i\in [m],
|\bar{z}_{ijka\bar a1}-\bar{z}_{ijka\bar a2}|\leq \hat{w}_{ijka\bar a}, \\
&(\bar{z}_{ijka\bar a1}, z_{ijka\bar a}, \bar{w}_{i})\in \MC(0, 1, -M_i, M_{i}), \\
&(\bar{z}_{ijka\bar a2}, z_{ijka\bar a}, \bar{w}_{j})\in \MC(0, 1, -M_j, M_{j}), \\
&\forall i\in C_a, j\in C_{\bar a}, k\in[\bar{\Omega}_{a\bar a}], a<\bar a\in A
\end{aligned}\right\},
\end{equation}
where for each $a<\bar a\in A$, we define $\bar{\Omega}_{a\bar a} = \left\lceil \log_{2}\left(\min\{m_a,m_{\bar a}\}\right) \right\rceil+1$ and
 \begin{equation*}
\Gamma_{a\bar a} = \left\{\bm{z}_{a\bar a}\in \{0,1\}^{m_a\times m_{\bar a} \times \bar{\Omega}_{a\bar a}}:
\begin{array}{l}\displaystyle
\sum_{i\in C_a}\sum_{k\in[\bar{\Omega}_{a\bar a}]} 2^{k-1} z_{ijka\bar a}=  m_a, \forall j\in C_{\bar a}, \\\displaystyle
\sum_{j\in C_{\bar a}}\sum_{k\in[\bar{\Omega}_{a\bar a}]} 2^{k-1} z_{ijka\bar a}=  m_{\bar a}, \forall i\in C_a
\end{array}
\right\}.
\end{equation*}

\end{restatable}

\noindent\textit{Proof. }
Letting $\bar{\pi}_{ija\bar a}=\pi_{ija\bar a}m_am_{\bar a}$, the set   $\F_{q}$ in \eqref{set_F_model1} can be reformulated as
\begin{equation}\nonumber
\F_{q}=\left\{(\bfx,\nu)\in\mathcal X\times \Re_+: \min_{\bar{\bm{\pi}}_{a\bar a}\in \bar{\Pi}_{a\bar a}}\left\{{\sum_{i\in C_a}\sum_{j\in C_{\bar a}} \frac{\bar{\pi}_{ija\bar a}}{m_am_{\bar a}}} \left|f(\bfx,{\bm\xi}_i)-f(\bfx,{\bm\xi}_j)\right|^{q}\right\}\leq \nu, \forall a<\bar a\in A\right\},
\end{equation}
where for each $a<\bar a\in A$, we have
\begin{equation}\nonumber
\bar{\Pi}_{a\bar a} = \left\{\bar{\bm{\pi}}_{a\bar a}\in \Re_{+}^{m_a\times m_{\bar a}}: \sum_{i\in C_a} \bar{\pi}_{ija\bar a}=m_a, \forall j\in C_{\bar a}, \sum_{j\in C_{\bar a}} \bar{\pi}_{ija\bar a}=m_{\bar a}, \forall i\in C_a \right\}. 
\end{equation}
We see that there exist integer solutions $\bar{\bm{\pi}}_{a\bar a}$ for this transportation problem according to \cite{rebman1974total}. Then, we obtain
\begin{equation}\nonumber
\F_{q}=\left\{(\bfx,\nu)\in\mathcal X\times \Re_+: \min_{\bar{\bm{\pi}}_{a\bar a}\in \bar{\Pi}_{a\bar a}\cap \Ze_{+}^{m_a\times m_{\bar a}}}\left\{{\sum_{i\in C_a}\sum_{j\in C_{\bar a}} \frac{\bar{\pi}_{ija\bar a}}{m_am_{\bar a}}} \left|f(\bfx,{\bm\xi}_i)-f(\bfx,{\bm\xi}_j)\right|^{q}\right\}\leq \nu, \forall a<\bar a\in A\right\},\\
\end{equation}
Since there exists $\bar{\bm{\pi}}_{a\bar a}\in \bar{\Pi}_{a\bar a}\cap \Ze_{+}^{m_a\times m_{\bar a}}$ for $\F_{q}$, we have 
\begin{equation}\nonumber
\F_{q}=\left\{(\bfx,\nu)\in\mathcal X\times \Re_+:\exists\bar{\bm{\pi}}_{a\bar a}\in \bar{\Pi}_{a\bar a}\cap \Ze_{+}^{m_a\times m_{\bar a}}, {\sum_{i\in C_a}\sum_{j\in C_{\bar a}} \frac{\bar{\pi}_{ija\bar a}}{m_am_{\bar a}}} \left|f(\bfx,{\bm\xi}_i)-f(\bfx,{\bm\xi}_j)\right|^{q}\leq \nu, \forall a<\bar a\in A\right\},\\
\end{equation}
Next, we binarize the
 
integer matrix variables $\bar{\bm{\pi}}_{a\bar a}$ using the expansion 
\[\bar{\pi}_{ija\bar a} = \sum_{k\in[\bar{\Omega}_{a\bar a}]} 2^{k-1} z_{ijka\bar a} \]
where $\bar{\Omega}_{a\bar a} = \left\lceil \log_{2}\left(\min\{m_a,m_{\bar a}\}\right) \right\rceil+1$ and $z_{ijka\bar a}\in\{0,1\}$ for all $i\in C_a$, $j\in C_{\bar a}$, and $k\in[\bar{\Omega}_{a\bar a}]$. 

We thus obtain
\begin{equation}\nonumber
\F_{q}=\left\{(\bfx,\nu)\in\mathcal X\times \Re_+: {\sum_{i\in C_a}\sum_{j\in C_{\bar a}}\sum_{k\in[\bar{\Omega}_{a\bar a}]} \frac{2^{k-1}z_{ijka\bar a}}{m_am_{\bar a}}} \left|f(\bfx,{\bm\xi}_i)-f(\bfx,{\bm\xi}_j)\right|^{q}\leq \nu,\bm{z}_{a\bar a}\in \Gamma_{a\bar a}, \forall a<\bar a\in A\right\}, 
\end{equation}
where for each $a<\bar a\in A$, we have
 \begin{equation}\nonumber
\Gamma_{a\bar a} = \left\{\bm{z}_{a\bar a}\in \{0,1\}^{m_a\times m_{\bar a} \times \bar{\Omega}_{a\bar a}}:
\begin{array}{l}\displaystyle
\sum_{i\in C_a}\sum_{k\in[\bar{\Omega}_{a\bar a}]} 2^{k-1} z_{ijka\bar a}=  m_a, \forall j\in C_{\bar a}, \\\displaystyle
\sum_{j\in C_{\bar a}}\sum_{k\in[\bar{\Omega}_{a\bar a}]} 2^{k-1} z_{ijka\bar a}=  m_{\bar a}, \forall i\in C_a 
\end{array}
\right\}.
\end{equation}
Then, letting $\bar{w}_{i}=f(\bfx,{\bm\xi}_i)$ and $\left|z_{ijka\bar a} \bar{w}_{i}-z_{ijka\bar a} \bar{w}_{j}\right|^{q}\leq \hat{w}_{ijka\bar a}$ 
can further linearize the set $\F_{q}$ as follows
\begin{equation}\nonumber
\F_{q}=\left\{(\bfx,\nu)\in\mathcal X\times \Re_+: \begin{aligned}
&{\sum_{i\in C_a}\sum_{j\in C_{\bar a}}\sum_{k\in[\bar{\Omega}_{a\bar a}]} \frac{2^{k-1}\hat{w}_{ijka\bar a}^q}{m_am_{\bar a}}}\leq \nu, \forall a<\bar a\in A, \\
&\left|z_{ijka\bar a} \bar{w}_i-z_{ijka\bar a} \bar{w}_j\right| \leq \hat{w}_{ijka\bar a}, \forall i\in C_a, j\in C_{\bar a}, k\in[\bar{\Omega}_{a\bar a}], a<\bar a\in A,\\
&\bm{z}_{a\bar a}\in \Gamma_{a\bar a}, \forall a<\bar a\in A, (\bfx, \bar{w}_{i})\in X_{i}, \forall i\in [m]
\end{aligned}\right\}.
\end{equation}
The conclusion follows from using McCormick representation of bilinear terms $\{z_{ijka\bar a}\bar{w}_{i}\}_{i\in C_a, j\in C_{\bar a}, k\in[\bar{\Omega}_{a\bar a}], a<\bar a\in A}$ and $\{z_{ijka\bar a}\bar{w}_{j}\}_{i\in C_a, j\in C_{\bar a}, k\in[\bar{\Omega}_{a\bar a}], a<\bar a\in A}$, and invoking the definition of the sets $\{X_i\}_{i\in [m]}$. 
\QEDA

Note that the support size of $\bar{\bm{\pi}}_{a\bar a}$ is $m_a+m_{\bar{a}}$. This motivates us to introduce new binary variables $\hat{\bm z}_{a\bar a}$ 
such that $z_{ijka\bar a}\leq \hat{z}_{ija\bar a}$ for each $i\in C_a, j\in C_{\bar a}, k\in[\bar{\Omega}_{a\bar a}], a<\bar a\in A$ and  obtain the following inequalities 
 
 valid 
 for the set $\F_q$ as 
\[\sum_{i\in C_a}\sum_{j\in C_{\bar a}}\hat{z}_{ija\bar a} \leq m_a+m_{\bar{a}},\]
for all $a<\bar a\in A$.

\subsection{Complementary Formulation: Linearizing the Complementary Slackness Constraints}\label{sec_M2}
In this subsection, we propose the second formulation of the set $\F_q$ using linear programming complementary slackness. According to the definition of the sets $\{X_i\}_{i\in [m]}$, we can represent the set $\F_q$ in \eqref{set_F_model1} as
\begin{equation}\label{set_F_equiv}
\F_{q}=\left\{(\bfx,\nu)\in\mathcal X\times \Re_+: \begin{aligned}
&  \min_{\bm{\pi}_{a\bar a}\in \Pi_{a\bar a}}\left\{\sum_{i\in C_a}\sum_{j\in C_{\bar a}} \pi_{ija\bar a}w_{ij}
\right\} \leq \nu,\forall a<\bar a\in A,\\
 & (\bfx, \bar{w}_{i})\in X_{i}, w_{ij}\geq \hat{w}_{ij}^q, \hat{w}_{ij}\geq |\bar{w}_{i}-\bar{w}_{j}|, \forall i\in [m], j\in [m]
\end{aligned} \right\},
\end{equation}
and we have $w_{ij}\leq (M_i+M_j)^q$ for each $(i,j)\in [m]\times [m]$. 
For each $a<\bar a\in A$, the dual of the left-hand side of the first constraint system in~\eqref{set_F_equiv} is
\begin{subequations}\label{set_F_equiv_dual}
\begin{align}
\max_{\bm{\mu}_{a\bar a},\bm{\lambda}_{a\bar a}} \quad& \frac{1}{m_{\bar a}}\sum_{j\in C_{\bar a}}\lambda_{ja\bar a} + \frac{1}{m_a}\sum_{i\in C_a}\mu_{ia\bar a}\leq\nu ,\\
\text{s.t.} \quad& \mu_{ia\bar a}+\lambda_{ja\bar a}\leq  w_{ij} , \forall i\in C_a, j\in C_{\bar a}.
\end{align}
\end{subequations}
According to linear programming complementary slackness, the system of linear inequalities in~\eqref{set_F_equiv_dual} is equivalent to
\begin{subequations}\label{set_F_dual}
\begin{align}
& \frac{1}{m_{\bar a}}\sum_{j\in C_{\bar a}}\lambda_{ja\bar a} + \frac{1}{m_a}\sum_{i\in C_a}\mu_{ia\bar a}\leq\nu ,\label{set_F_dual_constr1}\\
& \sum_{i\in C_a} \pi_{ija\bar a}=\frac{1}{m_{\bar a}}, \forall j\in C_{\bar a}, \\
& \sum_{j\in C_{\bar a}} \pi_{ija\bar a}=\frac{1}{m_a}, \forall i\in C_a, \\
& \pi_{ija\bar a}\geq 0, \forall i\in C_a, j\in C_{\bar a},\\
&  w_{ij} -\mu_{ia\bar a}-\lambda_{ja\bar a}\geq 0, \forall i\in C_a, j\in C_{\bar a},\label{set_F_dual_constr5}\\
& \pi_{ija\bar a}\left( w_{ij} -\mu_{ia\bar a}-\lambda_{ja\bar a}\right)=0, \forall i\in C_a, j\in C_{\bar a}.\label{set_F_dual_constr6}
\end{align}
\end{subequations}
Then, linearizing the complementary slackness constraints  \eqref{set_F_dual_constr6} allows us to derive the second MICP-R.
\begin{restatable}{theorem}{thmmodelsec}\textbf{(Complementary Formulation)}\label{thm_model_2}
Suppose that the set $X_{i}=\{(\bfx,\bar{w}_{i})\in \mathcal X\times\Re: f(\bfx,{\bm\xi}_i)=\bar{w}_{i}\}$ is MICP-R for each $i\in[m]$, $M_{i}\geq \max_{\bm{x}\in \mathcal X,\eqref{eq_tol}} |f(\bfx,{\bm\xi}_i)|$ for each $i\in[m]$, and $\hat{M}_{a\bar a}=\sum_{(i,j)\in C_a\times C_{\bar a}}(M_i+M_j)^q$ for each $a<\bar a\in A$. Then the set $\F_q$ is equivalent to
\begin{equation}\label{model2_WD1}
\F_q=\left\{\bfx: \begin{aligned}
& \frac{1}{m_{\bar a}}\sum_{j\in C_{\bar a}}\lambda_{ja\bar a} + \frac{1}{m_a}\sum_{i\in C_a}\mu_{ia\bar a}\leq\nu, \forall a<\bar a\in A,\\
& \sum_{i\in C_a} \pi_{ija\bar a}=\frac{1}{m_{\bar a}}, \forall j\in C_{\bar a}, a<\bar a\in A,
\sum_{j\in C_{\bar a}} \pi_{ija\bar a}=\frac{1}{m_a}, \forall i\in C_a,a<\bar a\in A, \\
& (\bfx, \bar{w}_{i})\in X_{i},\forall i\in [m], w_{ij}\geq \hat{w}_{ij}^q, \hat{w}_{ij}\geq |\bar{w}_{i}-\bar{w}_{j}|, \forall i\in [m], j\in [m],\\
& w_{ij}-\mu_{ia\bar a}-\lambda_{ja\bar a}\geq 0, w_{ij}-\mu_{ia\bar a}-\lambda_{ja\bar a}\leq\hat{M}_{a\bar a}\left(1-{z}_{ija\bar a}\right),\\
& \pi_{ija\bar a}\leq \min\{m_a^{-1},m_{\bar a}^{-1}\}{z}_{ija\bar a}, \pi_{ija\bar a}\geq 0, {z}_{ija\bar a}\in\{0,1\},
\forall i\in C_a, j\in C_{\bar a}, a<\bar a\in A
\end{aligned}\right\}.
\end{equation}
\end{restatable}

\noindent\textit{Proof.} 
\begin{subequations}

By introducing a large constant $\hat{M}_{a\bar a}$ for each $a<\bar a\in  A$, \eqref{set_F_dual_constr6} can be linearized as
\begin{equation}\label{set_F_dual_constr6_M}
\begin{aligned}
 & \pi_{ija\bar a}\leq \min\{m_a^{-1},m_{\bar a}^{-1}\}{z}_{ija\bar a},\; w_{ij}-\mu_{ia\bar a}-\lambda_{ja\bar a}\leq\hat{M}_{a\bar a}\left(1-{z}_{ija\bar a}\right), \\
  &{z}_{ija\bar a}\in\{0,1\}, \forall i\in C_a, j\in C_{\bar a}.  
\end{aligned}
\end{equation}
It remains to show that for each pair $a<\bar a\in A$, the big-M value $\hat{M}_{a\bar a}=\sum_{(i,j)\in C_a\times C_{\bar a}}(M_i+M_j)^q$ suffices. That is, any dual feasible solution satisfies  $w_{ij}-\mu_{ia\bar a}-\lambda_{ja\bar a}\leq \hat{M}_{a\bar a}$ for each $i\in C_a,j\in C_{\bar a}, a<\bar a\in A$.
From \eqref{set_F_dual_constr1}, we can get
\begin{equation}\label{set_F_dual_constr1_ineq}
0\leq\frac{1}{m_am_{\bar a}}\sum_{(i,j)\in C_a\times C_{\bar a}}(\mu_{ia\bar a}+\lambda_{ja\bar a})\leq\nu .
\end{equation}
According to \eqref{set_F_dual_constr5}, we also know that 
\begin{equation}\label{set_F_dual_constr5_ineq}
\mu_{ia\bar a}+\lambda_{ja\bar a}\leq w_{ij}\leq (M_i+M_j)^q.
\end{equation}
Then, we have
\[\mu_{ia\bar a}+\lambda_{ja\bar a}\geq-\sum_{(i',j')\in C_a\times C_{\bar a}\setminus\{i,j\}}(\mu_{i'a\bar a}+\lambda_{j'a\bar a})\geq-\sum_{(i',j')\in C_a\times C_{\bar a}\setminus\{i,j\}}(M_{i'}+M_{j'})^q,\]
where the first inequality is due to \eqref{set_F_dual_constr1_ineq} and the second inequality is because of \eqref{set_F_dual_constr5_ineq}.
Thus, $\hat{M}_{a\bar a}$ can be found by
\begin{equation*}
\hat{M}_{a\bar a}:= (M_i+M_j)^q+\sum_{(i',j')\in C_a\times C_{\bar a}\setminus\{i,j\}}(M_{i'}+M_{j'})^q=\sum_{(i',j')\in C_a\times C_{\bar a}}(M_{i'}+M_{j'})^q\geq w_{ij}-\mu_{ia\bar a}-\lambda_{ja\bar a},
\end{equation*}
which give us $\hat{M}_{a\bar a}=\sum_{(i,j)\in C_a\times C_{\bar a}}(M_i+M_j)^q$.
\end{subequations}
\QEDA

\subsection{A Side Product: An Equivalent MICP-R Formulation for $\KSD($\boldmath{$x$}\unboldmath{$)$}}\label{sec_micpr_KSD}

In this subsection, we propose to represent the sublevel set of the Kolmogorov–Smirnov fairness measure $\KSD(\bfx)$, denoted as $\KSD_{\hat{\delta}}(\bfx)$, using the sets $\{\Omega_a(k)\}_{a\in A,k\in [m_a]}$ defined in \eqref{set_Qa}. We note that the Kolmogorov–Smirnov fairness problem (similar to \ref{eq_sp_fair}) admits the following form:
\begin{align}
v^*=\min_{\bfx \in \mathcal X}\left\{\KSD(\bfx): \E_{\Pr}[Q(\bfx,\tilde{\bm\xi})]\leq V^*+\epsilon|V^*|\right\}. 
 \label{eq_sp_fair_KSD}
\end{align}
Hence, if we can represent the sublevel set of the function $\KSD(\bfx)$, then we can simply run a binary search to find the best objective value of problem \eqref{eq_sp_fair_KSD}. The formulation of $\KSD_{\hat{\delta}}(\bfx)$ is shown below. 
\begin{restatable}{theorem}{thmMksd}\label{thm_M_KSD}
Let the quantile set be defined as $\Omega_a(k)=\{(\bfx,t_a) \in \mathcal X\times \Re: F^{-1}_{a}(k/m_a\mid\bfx)=t_a\}$ for each $a\in A$, which admits a MICP-R form \eqref{set_Qa}. Then for a given $\hat{\delta}\in \{|i/m_a-j/m_{\bar a}|\}_{i\in [0,m_a], j\in [0,m_{\bar{a}}],a< 
\bar{a}\in A}$, the set $\KSD_{\hat{\delta}}(\bfx)$ can be expressed as
\begin{equation}\nonumber\label{model_KSD}
\KSD_{\hat{\delta}}(\bfx)=\left\{\bfx\in\mathcal X:\begin{aligned}
& (\bfx, t_{ia})\in \Omega_a(i),\forall i\in [m_a],a\in A,\\
& t_{(\lfloor m_{\bar{a}}(i/m_a-\hat{\delta})_+\rfloor)\bar{a}}\leq t_{ia}  , \forall i\in [m_a], a< \bar{a}  \in A\\
&t_{(i+1)a} \leq t_{(\lfloor m_{\bar{a}}\max(i/m_a+\hat{\delta},1)\rfloor+1)\bar{a}}, \forall i\in [0,m_a-1], a< \bar{a}  \in A
\end{aligned}\right\},
\end{equation}
where we let $t_{0a}=-\infty$ and $t_{(m_a+1)a}=+\infty$ for any $a\in A$.
\end{restatable}
 
\noindent\textit{Proof. }
Recall that for any $\bfx\in \KSD_{\hat{\delta}}(\bfx)$, we have 
\begin{subequations}
\begin{align}
&\max_{a< \bar{a}\in A} \sup_{\tau} \left|F_{a}(\tau\mid\bfx)-F_{\bar{a}}(\tau\mid\bfx)\right|\leq\hat{\delta} ,\\
\Leftrightarrow& \sup_{\tau} \left|F_{a}(\tau\mid\bfx)-F_{\bar{a}}(\tau\mid\bfx)\right|\leq\nu , \forall a< \bar{a}\in A.\label{model_KSD_eq1}
\end{align}
\end{subequations}
Let us consider the possible values of $F_a(\tau\mid\bfx)$ as $\{ 0,1/m_a, \cdots, m_a/m_a\}$. There are three cases:
\begin{enumerate}[{Case }1.]
\item If $F_a(\tau\mid\bfx)=0$, then $-\infty<\tau\leq t_{1a}$. 
For any such $\tau$, we must have $\left|F_{a}(\tau\mid\bfx)-F_{\bar{a}}(\tau\mid\bfx)\right|\leq\hat{\delta} $. 
That is, we must have
\[0\leq F_{\bar{a}}(\tau\mid\bfx) \leq \hat{\delta},\]
or equivalently $-\infty<\tau<t_{(\lfloor m_{\bar{a}}\max(\hat{\delta},1)\rfloor+1)\bar{a}}$. Therefore, the following inequalities must hold
\[\tau\leq t_{1a}\leq t_{(\lfloor m_{\bar{a}}\max(\hat{\delta},1)\rfloor+1)\bar{a}}.\]

\item If $F_a(\tau\mid\bfx)=i/m_a$ for some $i\in [m_a-1]$, then $t_{ia}\leq \tau<t_{(i+1)a}$. For any such $\tau$, we must have $\left|F_{a}(\tau\mid\bfx)-F_{\bar{a}}(\tau\mid\bfx)\right|\leq\hat{\delta} $. That is, we must have
\[\left(\frac{i}{m_a}-\hat{\delta}\right)_+\leq F_{\bar{a}}(\tau\mid\bfx) \leq \max\left(\frac{i}{m_a}+\hat{\delta},1\right)\]
or equivalently $t_{(\lfloor m_{\bar{a}}(i/m_a-\hat{\delta})_+\rfloor)\bar{a}}\leq \tau<t_{(\lfloor m_{\bar{a}}\max(i/m_a+\hat{\delta},1)\rfloor+1)\bar{a}}$. Therefore, the following inequalities must hold
\[t_{(\lfloor m_{\bar{a}}(i/m_a-\hat{\delta})_+\rfloor)\bar{a}}\leq t_{ia}\leq t_{(i+1)a} \leq t_{(\lfloor m_{\bar{a}}\max(i/m_a+\hat{\delta},1)\rfloor+1)\bar{a}}.\]

\item If $F_a(\tau\mid\bfx)=1$, then $t_{(m_a)a}\leq \tau<+\infty$. For any such $\tau$, we must have $\left|F_{a}(\tau\mid\bfx)-F_{\bar{a}}(\tau\mid\bfx)\right|\leq\hat{\delta} $. That is, we must have
\[\left(1-\hat{\delta}\right)_+\leq F_{\bar{a}}(\tau\mid\bfx) \leq 1,\]
or equivalently $t_{(\lfloor m_{\bar{a}}(1-\hat{\delta})_+\rfloor)\bar{a}}\leq \tau<+\infty$. Therefore, the following inequalities must hold
\[t_{(\lfloor m_{\bar{a}}(i/m_a-\hat{\delta})_+\rfloor)\bar{a}}\leq t_{(m_a)a}.\]
\end{enumerate}

If $F(\tau\mid\bfx)=i/m$, then we must have $i/m-\hat{\delta}\leq F_{a}(t_i\mid\bfx)\leq i/m+\hat{\delta}$ for all $a\in A$ to solve \eqref{model_KSD_eq1}. 

Finally, we observe that since all $\{\Pr_a\}_{a\in A}$ are equiprobable  discrete distributions, we must have $\hat{\delta}\in \{|i/m_a-j/m_{\bar a}|\}_{i\in [0,m_a], j\in [0,m_{\bar{a}}],a< 
\bar{a}\in A}$.
\QEDA

We remark that to solve \eqref{eq_sp_fair_KSD} to optimality, we can run the binary search to find the optimal $\hat{\delta}\in \{|i/m_a-j/m_{\bar a}|\}_{i\in [0,m_a], j\in [0,m_{\bar{a}}],a< 
\bar{a}\in A}$. That is, given a current $\hat{\delta}$ value, we optimize the total cost $\E[Q(\bfx,\tilde{\bm\xi})]$ subject to the set $\KSD_{\hat{\delta}}(\bm{x})$.  
Next, we check whether the optimal value is no larger than $V^*+\epsilon|V^*|$ or not. If yes, we decrease $\hat{\delta}$; otherwise, we increase it.

Alternatively, we can perform difference-of-convex (DC) method 
at each binary search step, where we solve a continuous relaxation by relaxing the binary variables to be continuous and rewrite the set $\Omega_a(k)$ to 
\begin{equation*}\label{set_Qa_v2}
\Omega_a(k)=\left\{(\bfx,t_{ka})\in \mathcal X\times \Re: \begin{aligned}
& \pi_{ika}\in [0,1], z_{ika}\in [0,1],\pi_{ika}\leq z_{ika},\forall i\in C_a,\\
& \sum_{i\in C_a}z_{ika}=k, \sum_{i\in C_a}\pi_{ika}=1,
t_{ka}=\sum_{i\in C_a}\hat{t}_{ika},\\
&(\bfx,\bar{w}_{i})\in X_i,
z_{ika}(t_{ka}-\bar{w}_{i})\geq 0,
\hat{t}_{ika}\leq \pi_{ika}\bar{w}_{i},\forall i\in C_a
\end{aligned}\right\}.
\end{equation*}
Here, we can rewrite each bilinear term as a difference between two convex functions.

\newpage

\section{Proofs}\label{proofs}

\subsection{Proof of \Cref{prop_comono}}\label{proof_prop_comono}

\propcomono*
\noindent\textit{Proof. } Given a standard uniform distribution $U$, let us define a joint distribution $\Qe_{a\bar a}$ such that $(f(\bfx,\tilde{\bm\xi}_a),f(\bfx,\tilde{\bm\xi}_{\bar a}))\stackrel{\text{d}}{=} 
(F_{a}^{-1}(U\mid\bfx),F_{\bar a}^{-1}(U\mid\bfx))$ 
for any pair $a<\bar a\in A$, and a fixed decision $\bfx \in\mathcal X$. Then, we have
\[\sqrt[^q]{\int_{\Xi\times\Xi}\left\|\bm{\zeta}_{1}-\bm{\zeta}_{2}\right\|^{q} \Qe_{a\bar a}(d\bm{\zeta}_{1},d\bm{\zeta}_{2})}=
\sqrt[^q]{\int_{0}^1\left|F_{a}^{-1}(u\mid\bfx)-F_{\bar a}^{-1}(u\mid\bfx)\right|^{q} du}\geq W_q\left(\Pr_{f(\bfx,\tilde{\bm\xi}_a)},\Pr_{f(\bfx,\tilde{\bm\xi}_{\bar a})}\right)\]
where the inequality is because $\Qe$ is an admissible joint distribution. 
According to \Cref{lem_icf_WD}, $\Qe$ is the ideal joint distribution when computing the Wasserstein fairness measure in \ref{eq_sp_fair}. 
This completes the proof.
\QEDA

\subsection{Proof of \Cref{prop_bernoulli}}\label{proof_prop_bernoulli}
\propbernoulli*
\noindent\textit{Proof. }
Since $f(\bfx,\tilde{\bm\xi})\in\{0,1\}$, we observe that 
$$F_a^{-1}(y\mid\bfx)=\begin{cases}
0& y\in(0, F_a(0\mid\bfx)], \\
1& y\in(F_a(0\mid\bfx), 1],
\end{cases}$$
for each $a\in A$. According to \Cref{lem_icf_WD}, the Wasserstein fairness measure $\WD_q(\bfx)$ can be simplified as
\begin{align*}
\WD_q(\bfx) =&\max_{a<\bar a\in A} \sqrt[^q]{\int_{0}^{1}\left|F^{-1}_{a}(y\mid\bfx)-F_{\bar a}^{-1}(y\mid\bfx)\right|^{q}dy}
=\max_{a<\bar a\in A} \sqrt[^q]{\left|F_{a}(0\mid\bfx)-F_{\bar a}(0\mid\bfx)\right|^{q}},\\
=& \max_{a<\bar a\in A} \left|F_{a}(0\mid\bfx)-F_{\bar a}(0\mid\bfx)\right|
= \max_{a<\bar a\in A} \left|\Pr_a\{f(\bfx,\tilde{\bm\xi}_{a})=0\}-\Pr_{\bar a}\{f(\bfx,\tilde{\bm\xi}_{\bar a})=0\}\right|=\DP(\bfx),
\end{align*}
where the second and third equalities are due to $\Pr(f(\bfx,\tilde{\bm\xi})\in\{0,1\})=1$ and the observation above, while the fourth one follows from the definition of cumulative distribution functions.
\QEDA

\subsection{Proof of \Cref{prop_bound_ksd_WF}}\label{proof_prop_bound_ksd_WF}
\propboundksdWF*
\noindent\textit{Proof. } We split the proof into two steps.

\noindent\textbf{Step 1.} We first derive the relationship between $\WD_1(\bfx)$ and $\KSD(\bfx)$.

We let $t_{1a\bar a}(\bfx)=\min\{\min_{t}\{t: F_{a}(t\mid\bfx)>0\}, \min_{t} \{t: F_{\bar a}(t\mid\bfx)>0\}\}$ and $t_{2a\bar a}(\bfx)=\max\{\max_{t} \{t: F_{a}(t\mid\bfx)<1\}, \max_{t} \{t: F_{\bar a}(t\mid\bfx)<1\}\}$. Then, according to \Cref{lem_icf_WD}, we have 
\begin{align*}
\WD_1(\bfx) =&\max_{a<\bar a\in A} \int_{t} \left|F_{a}(t\mid\bfx)-F_{\bar a}(t\mid\bfx)\right|dt,\\
=&\max_{a<\bar a\in A} \int_{t_{1a\bar a}(\bfx)}^{t_{2a\bar a}(\bfx)} \left|F_{a}(t\mid\bfx)-F_{\bar a}(t\mid\bfx)\right|dt,\\
\leq&\max_{a<\bar a\in A}(t_{2a\bar a}(\bfx)-t_{1a\bar a}(\bfx))\max_{a<\bar a\in A} \sup_{t} \left|F_{a}(t\mid\bfx)-F_{\bar a}(t\mid\bfx)\right|=\max_{a<\bar a\in A}(t_{2a\bar a}(\bfx)-t_{1a\bar a}(\bfx))\KSD(\bfx),
\end{align*}
which yields the lower bound on $\KSD(\bfx)$.

To establish the upper bound on $\KSD(\bfx)$, we let 
$\Delta_{a\bar a}(\bfx)=\{\bar{t}: |F_{a}(\bar{t}\mid\bfx)-F_{\bar a}(\bar{t}\mid\bfx)|=\sup_{t} |F_{a}(t\mid\bfx)-F_{\bar a}(t\mid\bfx)|\}$ and use $\mu(\Delta_{a\bar a}(\bfx))$ 
to denote the Lebesgue measure of the set $\Delta_{a\bar a}(\bfx)$, which is positive since all the groups are finitely distributed. Then, we have 
\begin{align*}
\WD_1(\bfx)=&\max_{a<\bar a\in A} \int_{t} \left|F_{a}(t\mid\bfx)-F_{\bar a}(t\mid\bfx)\right|dt,\\
\geq&\max_{a<\bar a\in A} \int_{\Delta_{a\bar a}(\bfx)} \left|F_{a}(t\mid\bfx)-F_{\bar a}(t\mid\bfx)\right|dt,\\
\geq&\min_{a<\bar a\in A}\mu(\Delta_{a\bar a}(\bfx))\max_{a<\bar a\in A} \sup_{t} \left|F_{a}(t\mid\bfx)-F_{\bar a}(t\mid\bfx)\right|=\min_{a<\bar a\in A}\mu(\Delta_{a\bar a}(\bfx))\KSD(\bfx).
\end{align*}
Combining both lower and upper bounds, we obtain the desired inequalities
\[\frac{1}{\max_{a<\bar a\in A}(t_{2a\bar a}(\bfx)-t_{1a\bar a}(\bfx))}\WD_1(\bfx) \leq \KSD(\bfx) \leq \frac{1}{\min_{a<\bar a\in A}\mu(\Delta_{a\bar a}(\bfx))}\WD_1(\bfx).\]

\noindent\textbf{Step 2.} Next, we derive the bounds between $\WD_1(\bfx)$ and $\WD_q(\bfx)$ for any $q\in[1,\infty]$. 
According to \Cref{lem_icf_WD}, we know that
\begin{align*}
\WD_q(\bfx)=&\max_{a<\bar a\in A} \sqrt[^q]{\int_{0}^{1}\left|F_{a}^{-1}(y\mid\bfx)-F_{\bar a}^{-1}(y\mid\bfx)\right|^{q}dy},\\
=&\max_{a<\bar a\in A} \sqrt[^q]{\sum_{j\in J_{a\bar a}\setminus \{1\}} w_{ja\bar a}(\bfx) \left|F_{a}^{-1}(b_{ja\bar a}(\bfx)\mid\bfx)-F_{\bar a}^{-1}(b_{ja\bar a}(\bfx)\mid\bfx)\right|^{q}},\\
\leq&\max_{a< \bar a\in A} \sqrt[^q]{\sum_{j\in J_{a\bar a}\setminus \{1\}} \frac{w_{ja\bar a}(\bfx)^q}{\eta(\bfx)^{q-1}} \left|F_{a}^{-1}(b_{ja\bar a}(\bfx)\mid\bfx)-F_{\bar a}^{-1}(b_{ja\bar a}(\bfx)\mid\bfx)\right|^{q}},\\
\leq&\eta(\bfx)^{\frac{1-q}{q}}\max_{a< \bar a\in A}\sum_{j\in J_{a\bar a}\setminus \{1\}} w_{ja\bar a} (\bfx)\left|F_{a}^{-1}(b_{ja\bar a}(\bfx)\mid\bfx)-F_{\bar a}^{-1}(b_{ja\bar a}(\bfx)\mid\bfx)\right|=\eta(\bfx)^{\frac{1-q}{q}}\WD_1(\bfx),
\end{align*}
where 
the second equality is due to \Cref{def_break_point}, the first inequality is due to $\eta(\bfx)=\max_{a<\bar a\in A, j\in J_a\setminus \{1\}}w_{ja\bar a}(\bfx)$, and the second one is because $\|\cdot\|_{q} \leq \|\cdot\|_{1}$. 

Meanwhile, by using H\"older's inequality,  we have
\[\sqrt[^1]{\int_{0}^{1}\left|F_{a}^{-1}(y\mid\bfx)-F^{-1}(y\mid\bfx)\right|^{1}dy} \leq \sqrt[^q]{\int_{0}^{1}\left|F_{a}^{-1}(y\mid\bfx)-F^{-1}(y\mid\bfx)\right|^{q}dy} \sqrt[^p]{\int_{0}^{1}dy}\]
for any $p, q\in [1,\infty)$ with $1/p+1/q=1$,
Thus, the following inequality holds
\[\int_{0}^{1}\left|F_{a}^{-1}(y\mid\bfx)-F_{\bar a}^{-1}(y\mid\bfx)\right|^{1}dy \leq \sqrt[^q]{\int_{0}^{1}\left|F_{a}^{-1}(y\mid\bfx)-F_{\bar a}^{-1}(y\mid\bfx)\right|^{q}dy},\]
i.e., $\WD_1(\bfx)\leq \WD_q(\bfx)$. This concludes the proof.
\QEDA

\subsection{Proof of \Cref{thm_np_hard_dfso}}\label{proof_thm_np_hard_dfso}
\thmnpharddfso*
\noindent\textit{Proof. }
We derive a reduction from the chance constrained optimization problem, which is strongly NP-hard \citep{ahmed2018relaxations}. 
Let us consider the following feasibility problem of the generic chance-constrained stochastic program:
\begin{quote}
Does there exist a feasible solution to the following chance-constrained set 
\[\left\{(\bfx,\bfz)\in \mathcal X:\bfz\in \{0,1\}^{m'}, \sum_{i\in [m']}z_i=m'-k, \right\}\]
where $\mathcal X=\{(\bfx,\bfz)\in\Re^{n}\times[0,1]^{m'}: \hat{\bm A}_i \bfx \geq \bm b_i -\bm M(1-z_i), \forall i\in [m']\}$ with large coefficients $\bm M$ for each $i\in [m']$?
\end{quote}
Let us consider a special case of \ref{eq_sp_fair} with $|A|=2$, $m_a=m_{\bar{a}}=m'$, $C_1=[m']$, $C_2=[m'+1,2m']$, and ${\bm\xi}\in \Re^{n+m'+1}$ with
\begin{align*}
f((\bfx,\bfz),{\bm\xi})={\bm\xi}_{1:n}^\top \bfx+ {\bm\xi}_{[n+1:n+m']}^\top \bfz+\xi_{n+m'+1}.
\end{align*}
Specifically, for each individual $i\in [2m']$, we design their corresponding scenario $\bm{\xi}_i$ such that
\begin{align*}
& f((\bfx,\bfz),{\bm\xi}_i)=z_i, \forall i\in [m'],\\
& f((\bfx,\bfz),{\bm\xi}_j)=0, \forall j\in [m'+1:m'+k],\\
& f((\bfx,\bfz),{\bm\xi}_j)=1, \forall j\in [m'+k+1:2m'].
\end{align*}
Assuming that $\epsilon=\infty$, this particular \ref{eq_sp_fair} becomes
\begin{align*}
v^*(q)=\min_{(\bfx, \bfz)\in \mathcal X} \quad & \WD_q^q((\bfx, \bfz)).
\end{align*}
Hence, we see that $v^*(q)=0$ if and only if there exists a point $(\bfx,\bfz)\in\mathcal X$ such that $\bfz\in \{0,1\}^{m'}, \sum_{i\in [m']}z_i=m'-k$. In other words, $v^*(q)=0$ if and only if $(\bfx,\bfz)$ is feasible to the chance-constrained stochastic program. This completes the proof.
\QEDA

\subsection{Proof of \Cref{prop_quant_set}}\label{proof_prop_quant_set}
\propquantset*
\noindent\textit{Proof. } By definition of the inverse distribution function, we have
\begin{align*}
\Omega_a(k)&=\left\{(\bfx,t_{ka})\in \mathcal X\times \Re: F^{-1}_{a}(k/m_a\mid\bfx)=t_{ka}\right\}\\
&=\left\{(\bfx,t_{ka})\in \mathcal X\times \Re: \begin{aligned}
& \pi_{ika}\in\{0,1\}, z_{ika}\in\{0,1\},\pi_{ika}\leq z_{ika},(\bfx,\bar{w}_{i})\in X_i,\forall i\in C_a,\\
& \sum_{i\in C_a}z_{ika}=k, \sum_{i\in C_a}\pi_{ika}=1,t_{ka}=\sum_{i\in C_a}\bar{w}_{i}\pi_{ika},\\
& t_{ka}\geq \bar{w}_{i}-(M_i+M_{(k)})(1-z_{ika}), 
t_{ka}\leq \bar{w}_{i}+(M_i+M_{(k)})z_{ika},\\ 
\end{aligned}\right\},
\end{align*}
for each $k\in [m_a]$ and $a\in A$. Next, we arrive at the conclusion by linearizing the bilinear terms with the McCormick representation.
\QEDA

\subsection{Proof of \Cref{thm_M3}}\label{proof_thm_M3}

\thmMthird*

\noindent\textit{Proof. }
We split the proof into two steps.

\noindent\textbf{Step 1.} First, we will reformulate $\WD_q^q(\bfx)$.
According to \eqref{eq_lem_inv_cdf_WD1}, we can represent $\WD_q^q(\bfx)\leq \nu$ as
\begin{equation}\label{eq_lem_inv_cdf_WD1_epsilon}
\sum_{i\in [\hat{m}_{a\bar a}-1]}(\hat{b}_{(i+1)a\bar a}-\hat{b}_{ia\bar a})\left|F^{-1}_{a}(\hat{b}_{(i+1)a\bar a}\mid\bfx)-F^{-1}_{\bar a}(\hat{b}_{(i+1)a\bar a}\mid\bfx)\right|^{q}\leq\nu, \forall a<\bar a\in A.
\end{equation}
Suppose $F^{-1}_{a}(y\mid\bfx)=t_{ja}$ for $j\in [m_a], a\in A$. Let $\delta_{ija\bar a 1}=\I( (\hat{b}_{ia\bar a},\hat{b}_{(i+1)a\bar a}]\subseteq ((j-1)/m_a, j/m_a])$ for each $i\in [\hat{m}_{a\bar a}-1], j\in [m_a], a<\bar a\in A$ and $\delta_{ija\bar a 2}=\I( (\hat{b}_{ia\bar a},\hat{b}_{(i+1)a\bar a}]\subseteq ((j-1)/m_{\bar a}, j/m_{\bar a}])$ for each $i\in [\hat{m}_{a\bar a}-1], j\in [m_{\bar a}], a<\bar a\in A$.
Then, the constraints \eqref{eq_lem_inv_cdf_WD1_epsilon} are equivalent to
\begin{equation}
\sum_{i\in [\hat{m}_{a\bar a}-1]}(\hat{b}_{(i+1)a\bar a}-\hat{b}_{ia\bar a}) \eta_{ia\bar a}^{q}\leq\nu,\forall a<\bar a\in A,
\end{equation}
where
\begin{equation}\label{model3_WD_t}
\left|\sum_{j\in[m_a] }\delta_{ija\bar a 1}t_{ja}-\sum_{j\in [m_{\bar a}]}\delta_{ija\bar a 2}t_{j\bar a}\right|\leq \eta_{ia\bar a}, \forall i\in [\hat{m}_{a\bar a}-1], a<\bar a\in A.
\end{equation}

\noindent\textbf{Step 2.} 
By choosing $(\bfx, t_{ja})\in \Omega_a(j)$ for all $j\in [m_a], a\in A$, we have the formulation \eqref{model3_WD} for the set $\F_q$.

\QEDA

\subsection{Proof of \Cref{prop_M1_M2_conti}}\label{proof_prop_M1_M2_conti}
\propMoneMtwoconti*

\noindent\textit{Proof. }
\begin{enumerate}[(i)]
\item Since the optimal value of the continuous relaxation of the Discretized Formulation is at least zero, it suffices to show that there exists a feasible solution for the continuous relaxation of the Discretized Formulation such that its objective value is zero.

In the Discretized Formulation, we choose any $\bm{x}\in\X$ and $\bar{w}_i=f(\bm{x},\bm{\xi}_i)$ for any $i\in [m]$. We also let  $\nu=0,\hat{w}_{ijka\bar{a}}=\bar{z}_{ijka\bar a1}=\bar{z}_{ijka\bar a2}=0$ and 
$z_{ijka\bar a}=2^{1-k}/\bar{\Omega}_{a\bar a},$
for all $i\in C_a, j\in C_{\bar a},$ and $k\in[\bar{\Omega}_{a\bar a}]$. Then we have
\[\sum_{i\in C_a}\sum_{k\in[\bar{\Omega}_{a\bar a}]} 2^{k-1} z_{ijka\bar a} = m_a, \forall j\in C_{\bar a}, \sum_{j\in C_{\bar a}}\sum_{k\in[\bar{\Omega}_{a\bar a}]} 2^{k-1} z_{ijka\bar a} = m_{\bar a}, \forall i\in C_a,\]
for all $a<\bar a\in A$. 
It remains to show that
\[(\bar{z}_{ijka\bar a1}, z_{ijka\bar a}, \bar{w}_{i})\in \MC(0, 1, -M_i, M_{i}), (\bar{z}_{ijka\bar a2}, z_{ijka\bar a}, \bar{w}_{j})\in \MC(0, 1, -M_j, M_{j}),\]
for all $i\in C_a, j\in C_{\bar a}, k\in[\bar{\Omega}_{a\bar a}], a<\bar a\in A$. This is true by choosing 
\[M_{i}\geq \max_{a<\bar a\in A}\bar{\Omega}_{a\bar a}/(\bar{\Omega}_{a\bar a}-1)\max_{\bm{x}\in \mathcal X,\eqref{eq_tol}} |f(\bfx,{\bm\xi}_i)|\]
for each $i\in[m]$. Hence, we see that $(\bm{x},\bar{\bm w},\nu,\bm{z},\bar{\bm{z}},\hat{\bm w})$ satisfies the constraints in \eqref{WD_q_trans_micpr}, which yields an objective value zero. 

\item Similarly, in the Complementary Formulation, we choose any $\bm{x}\in\X$ and $\bar{w}_i=f(\bm{x},\bm{\xi}_i)$ for any $i\in [m]$. We let $z_{ija\bar{a}}=\min\{m_a,m_{\bar a}\}/(m_a m_{\bar{a}}), w_{ij}= \hat{w}_{ij}^q, \hat{w}_{ij}= |\bar{w}_{i}-\bar{w}_{j}|,$ and $\pi_{ija\bar{a}}=1/(m_a m_{\bar{a}})$ for all $i\in C_a, j\in C_{\bar a}, a<\bar a\in A$. We also let $\mu_{ia\bar a}=0$ and $\lambda_{ja\bar a}=0$ for all $i\in C_a, j\in C_{\bar a}, a<\bar a\in A$.  We note that
\[\sum_{i\in C_a} \pi_{ija\bar a}=\frac{1}{m_{\bar a}}, 
 \sum_{j\in C_{\bar a}} \pi_{ija\bar a}=\frac{1}{m_a}, \pi_{ija\bar a}\leq \min\{m_a^{-1},m_{\bar a}^{-1}\}{z}_{ija\bar a}, \pi_{ija\bar a}\geq 0,\]
for all $i\in C_a, j\in C_{\bar a}, a<\bar a\in A$.

It remains to show that
\[w_{ij}\leq\hat{M}_{a\bar a}\left(1-{z}_{ija\bar a}\right),\]
for all $i\in C_a, j\in C_{\bar a}, a<\bar a\in A$. This is true by choosing $\hat{M}_{a\bar a}\geq (m_a m_{\bar{a}})/(m_a m_{\bar{a}}-\min\{m_a,m_{\bar a}\})\sum_{(i,j)\in C_a\times C_{\bar a}}(M_i+M_j)^q$ for each $a<\bar a\in A$.
Hence, we see that $(\bm{x},\bar{\bm w},\hat{\bm w},\bm{w},\bm{\mu},\bm{\lambda},\nu,\bm{z},\bar{\bm{z}},\bm{\pi})$ satisfies the constraints in \eqref{model2_WD1}, which yields an objective value zero.

\item We observe that for any feasible solution of the Quantile Formulation, we must have $z_{i(m_a)a}=1$ for each $i\in C_a$ and $a\in A$. Therefore, $t_{(m_a)a}=F_{a}^{-1}(1\mid\bfx),t_{(m_{\bar a})\bar a}=F_{\bar a}^{-1}(1\mid\bfx)$ and $\eta_{\hat{m}_{a\bar a}-1}\geq |t_{(m_a)a}-t_{(m_{\bar a})\bar a}|$. That is,
\begin{align*}
\nu\geq \max_{a<\bar{a}\in A}\left(\hat{b}_{(\hat{m}_{a\bar a})a\bar a}-\hat{b}_{(\hat{m}_{a\bar a}-1)a\bar a}\right)|F_{a}^{-1}(1\mid\bfx)-F_{\bar a}^{-1}(1\mid\bfx)|^q.
\end{align*}

We show that this bound is tight by constructing a solution such that $t_{ka}=0$ for any $k\in [m_a-1]$ and $a\in A$. In fact, for any $\bm x\in\chi$, let $\bar{w}_i=f(\bm{x},\bm{\xi}_i)$ for each $i \in [m]$. Then for any $i\in C_a, k\in [m_a-1],$ and $a\in A$, we let $z_{ika}=k/m_a$, $\pi_{ika}=1/m_a$, $\hat{t}_{ika}=0$. It remains to show that for any $i\in C_a, k\in [m_a-1],$ and $a\in A$,
\[(\hat{t}_{ika},\pi_{ika},\bar{w}_{i})\in\MC(0,1,-M_i,M_i)\]
which must hold when $M_{i}\geq \max_{a\in A}m_a\max_{\bm{x}\in \mathcal X,\eqref{eq_tol}} |f(\bfx,{\bm\xi}_i)|$.

\item We observe that for any feasible solution of the Aggregate Quantile Formulation, we must have $z_{i(m_a)a}=1$ for each $i\in C_a$ and $a\in A$. Therefore, $\bar t_{(m_a)a}=\sum_{i\in C_{a}}F_{a}^{-1}(i/m_a\mid\bfx),\bar t_{(m_{\bar a})\bar a}=\sum_{i\in C_{\bar a}}F_{\bar a}^{-1}(i/m_{\bar a}\mid\bfx)$. According to the first two constraints in \eqref{model4_WD}, we have
\begin{align*}
\nu &\geq \max_{a<\bar{a}\in A}\sum_{i\in [\hat{m}_{a\bar a}-1]}\left(\hat{b}_{(i+1)a\bar a}-\hat{b}_{ia\bar a}\right) \left|\sum_{j\in[m_a] }\delta_{ija\bar a 1}t_{ja}-\sum_{j\in [m_{\bar a}]}\delta_{ija\bar a 2}t_{j\bar a}\right|^{q}\\
&\geq \max_{a<\bar{a}\in A}\left|\sum_{i\in [\hat{m}_{a\bar a}-1]}\left(\hat{b}_{(i+1)a\bar a}-\hat{b}_{ia\bar a}\right) \left[\sum_{j\in[m_a] }\delta_{ija\bar a 1}t_{ja}-\sum_{j\in [m_{\bar a}]}\delta_{ija\bar a 2}t_{j\bar a}\right]\right|^{q}\\
& = \max_{a<\bar{a}\in A}\left|m_a^{-1}\sum_{i\in C_{a}}F_{a}^{-1}(i/m_a\mid\bfx)-m_{\bar a}^{-1}\sum_{i\in C_{\bar a}}F_{\bar a}^{-1}(i/m_{\bar a}\mid\bfx)\right|^q
\end{align*}
where the second inequality is due to Jensen's inequality. 
\QEDA
\end{enumerate}

\subsection{Proof of \Cref{thm_gel_NP}}\label{proof_thm_gel_NP}

\thmgelNP*

\noindent\textit{Proof. }
We show a reduction from the integer programming feasibility problem, which is known to be strongly NP-complete. Consider the following feasibility problem of a binary integer program:
\begin{quote}
Does there exist a feasible solution to the binary program $X=\{\bfx\in \{-1,1\}^{n-1}: \bm{A} \bfx\geq \bm{b} \}$?
\end{quote}

We consider a special case of the Gelbrich bound \eqref{gelbrich_Cholesky} by letting $\epsilon=\infty$, $|A|=2$, $\bm{r}(\bfx,y)=(\bfx^\top,y)^\top$, $\bm{\mu}_a = \bm{\mu}_{\bar{a}} = \bm 0$, $\bm L_{a} = \begin{bmatrix} I_{m_{\bar{a}}} & \bm{0} \\ \bm{0} & 0 \end{bmatrix}$, $\bm L_{\bar{a}} = \begin{bmatrix} \bm{0} & \bm{0} \\ \bm{0} & 1 \end{bmatrix}$, and defining
\[\mathcal X:=\left\{(\bfx,y)\in [-1,1]^{n-1}\times\{\sqrt{n-1}\}: \bm{A} \bfx \geq \bm{b} \right\}.\]
Under this setting, the Gelbrich bound \eqref{gelbrich_Cholesky} simplifies to
\begin{subequations}\label{gelbrich_np_hard_ex}
\begin{align}
v_G=\min_{\bfx} \quad & \left(\sqrt{\bfx^{\top}\bfx} - \sqrt{n-1}\right)^2,\\
\text{s.t.} \quad
& \bfx\in [-1,1]^{n-1}, \bm{A} \bfx\geq \bm{b}.
\end{align}
\end{subequations}
We see that $v_G=0$ in the formulation \eqref{gelbrich_np_hard_ex}  if and only if there exists a binary feasible solution to the set $X$. Thus, the claim follows.
\QEDA

\subsection{Proof of \Cref{thm_gel_SDP}}\label{proof_thm_gel_SDP}

\thmgelSDP*

 \noindent\textit{Proof. }
We split the proof into two steps.\par
\noindent\textbf{Step 1.} We know that $\bm{Z}_{a\bar{a}}$ is positive semidefinite and $\left(Z_{a\bar{a}11} - 2Z_{a\bar{a}1(n+2)} +  Z_{a\bar{a}(n+2) (n+2)}\right)\geq0$. Therefore, by removing the second term in the left-hand side of \eqref{gelbrich_SDP_rank1_eq1}, we have
\begin{subequations}\label{gelbrich_SDP_relax_beta}
\begin{align}
v_{\underline G} \stackrel{}{\geq } \min_{\bfx\in\mathcal X,\bm s,\bm Z,\nu} \quad &\nu,\\
\text{s.t.} \quad
& \left|\bm\mu_{a}^{\top}\bm{r}(\bfx) -\bm\mu_{\bar{a}}^{\top}\bm{r}(\bfx) \right|^2 \leq \nu, \forall a<\bar{a}\in A,\\
& \eqref{eq_tol},\eqref{gelbrich_SDP_rank1_eq2}-\eqref{gelbrich_SDP_rank1_eq6}, \eqref{gelbrich_SDP_rank1_schur}
\nonumber
\end{align}
\end{subequations}
where the optimal value of the minimization problem is equal to $v_J(2)$. Thus, we have $v_{\underline G}\geq v_J(2)$.

\noindent\textbf{Step 2.} Let $(\bfx^*,\nu^*)$ denote the optimal solution of \eqref{lower_bound_Wq} for $q=2$. We compute $\bfz_{a}^* = \bm  L_{a}^{\top}\bm{r}(\bfx^*)$, $\sigma_{a}^* = \|\bfz_{a}^*\|_2$, $\bm s^*_{a\bar{a}} = \begin{bmatrix} \sigma^*_{a} & \bfz^*_{a} & \sigma^*_{\bar{a}}  & \bfz^*_{\bar{a}} \end{bmatrix}^{\top}$, and $\bm{Z}_{a\bar{a}}^* = \bm s_{a\bar{a}}^* \cdot \bm s_{a\bar{a}}^{*\top}$ for all $a<\bar{a}\in A$.
Suppose $\sigma_{a}^* \geq \sigma_{\bar{a}}^*$ and $\gamma = \sigma_{a}^*/\sigma_{\bar{a}}^*\geq1$. Let $\hat{\bm s}_{a\bar{a}} = \begin{bmatrix} \hat{\sigma}_{a} & \hat{\bfz}_{a} & \hat{\sigma}_{\bar{a}} & \hat{\bfz}_{\bar{a}} \end{bmatrix}^{\top} = \begin{bmatrix} \sigma_{a}^* & \bfz_{a}^* & \gamma\sigma_{\bar{a}}^*  & \gamma\bfz_{\bar{a}}^* \end{bmatrix}^{\top}$. For any given $\alpha\geq \bar{\alpha}$ (we will specify $\bar{\alpha}$ later), we construct $\hat{\bm{Z}}_{a\bar{a}}$ as
\[\hat{\bm{Z}}_{a\bar{a}} = \begin{bmatrix} \hat{\sigma}_{a}^{2}  & \bm{0} & \hat{\sigma}_{a}\hat{\sigma}_{\bar{a}} & \bm{0}\\ 
\bm{0} & \hat{\bfz}_{a}^{2} & \bm{0} & \bm{0} \\  
\hat{\sigma}_{a}\hat{\sigma}_{\bar{a}}  & \bm{0} &\hat{\sigma}_{\bar{a}}^{2} & \bm{0} \\ 
\bm{0}  & \bm{0} & \bm{0} & \hat{\bfz}_{\bar{a}}^{2} \end{bmatrix}.\]
It is evident that $(\bfx^*,\bm s^*,\alpha\hat{\bm Z},\nu^*)$ satisfies constraints \eqref{eq_tol}, \eqref{gelbrich_SDP_rank1_eq1}-\eqref{gelbrich_SDP_rank1_eq2}, \eqref{gelbrich_SDP_rank1_eq6}.

We next compute 
$\alpha \hat{\bm{Z}}_{a\bar{a}} - \bm s_{a\bar{a}}^*\cdot {\bm s_{a\bar{a}}^{*}}^{\top}$. Without loss of generality, we can assume that $z_{ia}^*\neq 0$ for all $i\in [n]$ and $a\in A$. We first observe that by dropping the first row and column, the resulting principal submatrix has rank $2n+1$. Thus, the rank of  $\hat{\bm{Z}}_{a\bar{a}}$ is at least $2n+1$. On the other hand, let $\bm{u}=[1,\bm{0},-1,\bm{0}]^\top$. Then we have $\bm{u}^\top \hat{\bm{Z}}_{a\bar{a}}\bm{u}=0$. Thus, the rank of  $\hat{\bm{Z}}_{a\bar{a}}$ is $2n+1$. Let $\{\bm{v}_i\}_{i\in [2n+1]}$ denote the nonzero eigenvectors of $\hat{\bm{Z}}_{a\bar{a}}$ with the corresponding positive eigenvalues $\{\lambda_i\}_{i\in [2n+1]}\subseteq \Re_{++}$. Then we can rewrite the matrix as
\[\hat{\bm{Z}}_{a\bar{a}}=\bm{V}\bm{\Lambda}\bm{V}^\top,\]
where $\bm{V}=[\bm{v}_1,\ldots, \bm{v}_{2n+1}]$ and $\bm{\Lambda}=\textrm{diag}(\bm\lambda)$. On the other hand, since $\bm{s}^*\in \textrm{span}(\{\bm{v}_i\}_{i\in [2n+1]}\cup\{\bm u\})$, we have
\[\bm{s}^*= t_0 \bm u+\sum_{i\in [2n+1]}q_i\bm{v}_i. \]
Thus, we have
\begin{align*}
\alpha \hat{\bm{Z}}_{a\bar{a}} - \bm s_{a\bar{a}}^*\cdot {\bm s_{a\bar{a}}^{*}}^{\top}
=\bm{V}(\alpha \bm{\Lambda}-\bm q\bm q^\top)\bm{V}^\top-t_0^2 \bm u \bm u^\top,
\end{align*}
where $t_0=(\bm s_{a\bar{a}}^*)^\top \bm u/2$. 
Since $\bm{\Lambda}\succ \bm 0$, we can choose $\bar \alpha=\max\{\bm q^\top \bm q/\min_{i\in [2n+1]}\lambda_i, 2t_0^2/(\beta\min_{i\in [2n+1]}\lambda_i) \}$ such that for any $\alpha\geq \bar \alpha$, we have
\begin{align*}
\alpha \hat{\bm{Z}}_{a\bar{a}} - \bm s_{a\bar{a}}^*\cdot {\bm s_{a\bar{a}}^{*}}^{\top}
=\bm{V}(\alpha \bm{\Lambda}-\bm q\bm q^\top)\bm{V}^\top-t_0^2 \bm u \bm u^\top\succeq -t_0^2 \bm u \bm u^\top \succeq -\beta\lambda_{\min}^+(\alpha \hat{\bm{Z}}_{a\bar{a}} ) \bm{I}_{2n+2}.
\end{align*}
This completes the proof.\QEDA

\subsection{Proof of \Cref{thm_exact_gelbrich}}\label{proof_thm_exact_gelbrich}

\thmexactgelbrich*
  \noindent\textit{Proof. }
Note that for any $a<\bar{a}\in A$, and an optimal comonotonic joint distribution $\Qe_{a,\bar{a}}$ of $f(\bfx,\tilde{\bm\xi}_{a})$ and $f(\bfx,\tilde{\bm\xi}_{\bar{a}})$ with marginals $\Pr_a,\Pr_{\bar a}$ such that $(f(\bfx,\tilde {\bm\xi}_{a})-\bm\mu_{a}^{\top}\bm{r}(\bfx)-s(\bfx), f(\bfx,\tilde{\bm\xi}_{\bar a})-\bm\mu_{\bar a}^{\top}\bm{r}(\bfx)-s(\bfx))\xrightarrow{m_{a}\rightarrow \infty,m_{\bar a}\rightarrow \infty}(\sqrt{\bm{r}(\bfx)^{\top}\bm\Sigma_{a}\bm{r}(\bfx)}\tilde u, \sqrt{\bm{r}(\bfx)^{\top}\bm\Sigma_{\bar a}\bm{r}(\bfx)}\tilde u)$, we have
\begin{align*}
&\E_{\Qe_{a,\bar{a}}}[|f(\bfx,\tilde{\bm\xi}_{a})-f(\bfx,\tilde{\bm\xi}_{\bar{a}})|^2]
 = \E_{\Pr_{a}}\left[\left(f(\bfx,\tilde{\bm\xi}_{a})-\bm\mu_{a}^{\top}\bm{r}(\bfx)-s(\bfx)\right)^2\right]+\E_{\Pr_{\bar a}}\left[\left(f(\bfx,\tilde{\bm\xi}_{\bar{a}})-\bm\mu_{\bar a}^{\top}\bm{r}(\bfx)-s(\bfx)\right)^2\right] 
 \\
 &-2\E_{\Pr_{a}}\left[\left(f(\bfx,\tilde{\bm\xi}_{a})-\bm\mu_{a}^{\top}\bm{r}(\bfx)-s(\bfx)\right)\left(\bm\mu_{a}^{\top}\bm{r}(\bfx) -\bm\mu_{\bar{a}}^{\top}\bm{r}(\bfx)\right)\right]\\
   &-2\E_{\Pr_{\bar a}}\left[\left(f(\bfx,\tilde{\bm\xi}_{\bar a})-\bm\mu_{\bar a}^{\top}\bm{r}(\bfx)-s(\bfx)\right)\left(\bm\mu_{a}^{\top}\bm{r}(\bfx) -\bm\mu_{\bar{a}}^{\top}\bm{r}(\bfx)\right)\right]\\
  &+ \left(\bm\mu_{a}^{\top}\bm{r}(\bfx) -\bm\mu_{\bar{a}}^{\top}\bm{r}(\bfx) \right)^2-2\E_{\Qe_{a,\bar{a}}}\left[\left(f(\bfx,\tilde{\bm\xi}_{a})-\bm\mu_{a}^{\top}\bm{r}(\bfx)-s(\bfx)\right)\left(f(\bfx,\tilde{\bm\xi}_{\bar{a}})-\bm\mu_{\bar a}^{\top}\bm{r}(\bfx)-s(\bfx)\right)\right]\\
  &\xrightarrow{m_{a}\rightarrow \infty,m_{\bar a}\rightarrow \infty}
  \left(\bm\mu_{a}^{\top}\bm{r}(\bfx) -\bm\mu_{\bar{a}}^{\top}\bm{r}(\bfx) \right)^2 + \left(\sqrt{\bm{r}(\bfx)^{\top}\bm\Sigma_{a}\bm{r}(\bfx)} - \sqrt{\bm{r}(\bfx)^{\top}\bm\Sigma_{\bar a}\bm{r}(\bfx)}\right)^2.
\end{align*}
Thus, the claim follows. 
\QEDA

\subsection{Proof of \Cref{thm_exact_gelbrich_gap}}\label{proof_thm_exact_gelbrich_gap}

\thmexactgelbrichgap*

 \noindent\textit{Proof. }
Note that the inequality $v_G\leq v^*$ is due to the derivation of the Gelbrich bound. It remains to establish the other direction. 
For notational convenience, we define $\hat{\Pr}_a$ as the true distribution of random variable $s(\bfx)+\bm\mu_{a}^{\top}\bm{r}(\bfx)+\sqrt{\bm{r}(\bfx)^{\top}\bm\Sigma_{a}\bm{r}(\bfx)}\tilde{u}_a$ for each $a\in A$, and let $\bar{u}_a$ be the discrete random variable with a uniform distribution on the samples $\{u_i\}_{i\in C_a}$, where  $\bar{\mu}_a$ and ${\rm var}(\bar{u}_a)$ are the corresponding sample mean and variance, respectively.   We also let $\Pr_a'$ be the empirical distribution of $s(\bfx)+\bm\mu_{a}^{\top}\bm{r}(\bfx)+\sqrt{\bm{r}(\bfx)^{\top}\bm\Sigma_{a}\bm{r}(\bfx)}\bar{u}_a$ for each $a\in A$, and $\bar{\Pr}$ be the joint distribution of all the random variables $\{\tilde{u}_a\}_{a\in A}$.
In this case, for any $a<\bar{a}\in A$, the Gelbrich bound reduces to
\begin{align*}
&\left(\bm\mu_{a}^{\top}\bm{r}(\bfx) -\bm\mu_{\bar{a}}^{\top}\bm{r}(\bfx) +
\sqrt{\bm{r}(\bfx)^{\top}\bm\Sigma_{a}\bm{r}}\bar{\mu}_a- \sqrt{\bm{r}(\bfx)^{\top}\bm\Sigma_{\bar a}\bm{r}(\bfx) }\bar{\mu}_{\bar a}\right)^2 + \\
&\left(\sqrt{\bm{r}(\bfx)^{\top}\bm\Sigma_{a}\bm{r}(\bfx){\rm var}(\bar{u}_a) }- \sqrt{\bm{r}(\bfx)^{\top}\bm\Sigma_{\bar a}\bm{r}(\bfx){\rm var}(\bar{u}_{\bar a}) }\right)^2.
\end{align*}

We split the proof into five steps. 

\noindent\textbf{Step 1.} 
Following the proof in \Cref{thm_exact_gelbrich}, for any $a<\bar{a}\in A$, we have
\begin{align*}
W_{2}^2\left(\hat{\Pr}_{a},\hat{\Pr}_{\bar{a}}\right)=\left(\bm\mu_{a}^{\top}\bm{r}(\bfx) -\bm\mu_{\bar{a}}^{\top}\bm{r}(\bfx) \right)^2+\left(\sqrt{\bm{r}(\bfx)^{\top}\bm\Sigma_{a}\bm{r}(\bfx)} - \sqrt{\bm{r}(\bfx)^{\top}\bm\Sigma_{\bar a}\bm{r}(\bfx)}\right)^2.
\end{align*}
According to the triangle inequality, for any $a<\bar{a}\in A$, we have
\begin{align*}
W_{2}\left(\Pr_{a}',\Pr_{\bar{a}}'\right)
\leq W_{2}\left(\Pr_{a}',\hat{\Pr}_{a}\right)+
W_{2}\left(\hat{\Pr}_{a},\hat{\Pr}_{\bar{a}}\right)
+W_{2}\left(\hat{\Pr}_{\bar{a}},\Pr_{\bar{a}}'\right).
\end{align*}

\noindent\textbf{Step 2.} In view of Theorem 2 in \cite{fournier2015rate}, for each $a\in A$, there exist two constants $c_{1a}>0,c_{2a}>0$ such that for any $\hat{\eta}_a>0$, we have
\begin{align*}
\bar{\Pr}\left\{W_{2}^2\left(\Pr_{a}',\hat{\Pr}_{a}\right)>\hat{\eta}_a\right\}
\leq c_{1a} \exp(-c_{2a}m_{a}\hat{\eta}_a).
\end{align*}
By letting the right-hand side probability be no larger than $\hat{\delta}\in (0,0.1)$, we obtain $\hat{\eta}_a:=-\log(c_{1a}/\hat{\delta})/(c_{2a}m_{a})$.

\noindent\textbf{Step 3.} We have 
\begin{align*}
\left(\bm\mu_{a}^{\top}\bm{r}(\bfx) -\bm\mu_{\bar{a}}^{\top}\bm{r}(\bfx) \right)^2\leq& 
\left(\bm\mu_{a}^{\top}\bm{r}(\bfx) -\bm\mu_{\bar{a}}^{\top}\bm{r}(\bfx) +
\sqrt{\bm{r}(\bfx)^{\top}\bm\Sigma_{a}\bm{r}}\bar{\mu}_a- \sqrt{\bm{r}(\bfx)^{\top}\bm\Sigma_{\bar a}\bm{r}(\bfx) }\bar{\mu}_{\bar a}\right)^2\\
&+\left|\bm\mu_{a}^{\top}\bm{r}(\bfx) -\bm\mu_{\bar{a}}^{\top}\bm{r}(\bfx) \right| \left|
\sqrt{\bm{r}(\bfx)^{\top}\bm\Sigma_{a}\bm{r}}\bar{\mu}_a- \sqrt{\bm{r}(\bfx)^{\top}\bm\Sigma_{\bar a}\bm{r}(\bfx) }\bar{\mu}_{\bar a}\right|\\
\leq&\left(\bm\mu_{a}^{\top}\bm{r}(\bfx) -\bm\mu_{\bar{a}}^{\top}\bm{r}(\bfx) +
\sqrt{\bm{r}(\bfx)^{\top}\bm\Sigma_{a}\bm{r}}\bar{\mu}_a- \sqrt{\bm{r}(\bfx)^{\top}\bm\Sigma_{\bar a}\bm{r}(\bfx) }\bar{\mu}_{\bar a}\right)^2+4M^2\left|\bar{\mu}_a- \bar{\mu}_{\bar a}\right|,
\end{align*}
where $M:=\max_{a\in A}\max_{\bfx\in\mathcal X}\left\{\max\{\sqrt{\bm{r}(\bfx)^{\top}\bm\Sigma_{a}\bm{r}(\bfx)},|\bm\mu_{a}^{\top}\bm{r}(\bfx)|\}:\eqref{eq_tol}\right\}$.
According to
the Chebyshev inequality, for each $a\in A$, the following probabilistic bound holds: 
\begin{align*}
\bar{\Pr}\left\{|\bar{\mu}_a|>\frac{1}{\sqrt{\hat{\eta}m_{a}}}\right\}
\leq \hat{\eta}.
\end{align*}
Thus, using the union bound,  we obtain 
\begin{align*}
\left(\bm\mu_{a}^{\top}\bm{r}(\bfx) -\bm\mu_{\bar{a}}^{\top}\bm{r}(\bfx) \right)^2\leq &
\left(\bm\mu_{a}^{\top}\bm{r}(\bfx) -\bm\mu_{\bar{a}}^{\top}\bm{r}(\bfx) +
\sqrt{\bm{r}(\bfx)^{\top}\bm\Sigma_{a}\bm{r}}\bar{\mu}_a- \sqrt{\bm{r}(\bfx)^{\top}\bm\Sigma_{\bar a}\bm{r}(\bfx) }\bar{\mu}_{\bar a}\right)^2\\
&+4M^2\left(\sqrt{\frac{1}{\sqrt{\hat{\eta}m_{a}}}}+\sqrt{\frac{1}{\sqrt{\hat{\eta}m_{\bar{a}}}}}\right)^2,
\end{align*}
that is,
\begin{align*}
&\left|\bm\mu_{a}^{\top}\bm{r}(\bfx) -\bm\mu_{\bar{a}}^{\top}\bm{r}(\bfx) \right|\leq 
\left|\bm\mu_{a}^{\top}\bm{r}(\bfx) -\bm\mu_{\bar{a}}^{\top}\bm{r}(\bfx) +
\sqrt{\bm{r}(\bfx)^{\top}\bm\Sigma_{a}\bm{r}}\bar{\mu}_a- \sqrt{\bm{r}(\bfx)^{\top}\bm\Sigma_{\bar a}\bm{r}(\bfx) }\bar{\mu}_{\bar a}\right|\\
&+2M\left(\sqrt{\frac{1}{\sqrt{\hat{\eta}m_{a}}}}+\sqrt{\frac{1}{\sqrt{\hat{\eta}m_{\bar{a}}}}}\right),
\end{align*}
with probability at least $1-2\hat{\eta}$.

\noindent\textbf{Step 4.} In view of
Theorem 4.7.1 in \cite{vershynin2018high} and the Markov inequality, for each $a\in A$, there exists a constant $c_{3a}>0$ such that
\begin{align*}
\bar{\Pr}\left\{|{\rm var}(\bar{u}_a)-1|>\frac{2c_{3a}}{\hat{\eta}\sqrt{m_{a}}}\right\}
\leq \hat{\eta}.
\end{align*}

Thus, according to the union bound, with probability at least $1-2\hat{\eta}$, we have
\begin{align*}
&\left(\sqrt{\bm{r}(\bfx)^{\top}\bm\Sigma_{a}\bm{r}(\bfx)} - \sqrt{\bm{r}(\bfx)^{\top}\bm\Sigma_{\bar a}\bm{r}(\bfx)}\right)^2\leq 
\left(\sqrt{\bm{r}(\bfx)^{\top}\bm\Sigma_{a}\bm{r}(\bfx){\rm var}(\bar{u}_a)} - \sqrt{\bm{r}(\bfx)^{\top}\bm\Sigma_{\bar a}\bm{r}(\bfx){\rm var}(\bar{u}_{\bar a})}\right)^2\\
&+4M^2\left(\sqrt{\frac{2c_{3a}}{\hat{\eta}\sqrt{m_{a}}}}+\sqrt{\frac{2c_{3\bar a}}{\hat{\eta}\sqrt{m_{\bar{a}}}}}\right)^2.
\end{align*}
That is,
\begin{align*}
&\left|\sqrt{\bm{r}(\bfx)^{\top}\bm\Sigma_{a}\bm{r}(\bfx)} - \sqrt{\bm{r}(\bfx)^{\top}\bm\Sigma_{\bar a}\bm{r}(\bfx)}\right|\leq 
\left|\sqrt{\bm{r}(\bfx)^{\top}\bm\Sigma_{a}\bm{r}(\bfx){\rm var}(\bar{u}_a)} - \sqrt{\bm{r}(\bfx)^{\top}\bm\Sigma_{\bar a}\bm{r}(\bfx){\rm var}(\bar{u}_{\bar a})}\right|\\
&+2M\left(\sqrt{\frac{2c_{3a}}{\hat{\eta}\sqrt{m_{a}}}}+\sqrt{\frac{2c_{3\bar a}}{\hat{\eta}\sqrt{m_{\bar{a}}}}}\right).
\end{align*}

\noindent\textbf{Step 5.} Combining all the steps together and using the union bound again, we have with probability at most $1-6\hat{\eta}$, 
\begin{align*}
&W_{2}^2\left(\Pr_{a},\Pr_{\bar{a}}\right)\leq 
\left(\bm\mu_{a}^{\top}\bm{r}(\bfx) -\bm\mu_{\bar{a}}^{\top}\bm{r}(\bfx) +
\sqrt{\bm{r}(\bfx)^{\top}\bm\Sigma_{a}\bm{r}}\bar{\mu}_a- \sqrt{\bm{r}(\bfx)^{\top}\bm\Sigma_{\bar a}\bm{r}(\bfx) }\bar{\mu}_{\bar a}\right)^2 + \\ 
&+ \left(\sqrt{\bm{r}(\bfx)^{\top}\bm\Sigma_{a}\bm{r}(\bfx){\rm var}(\bar{u}_a) }- \sqrt{\bm{r}(\bfx)^{\top}\bm\Sigma_{\bar a}\bm{r}(\bfx){\rm var}(\bar{u}_{\bar a}) }\right)^2+\bar{C}_1(\hat{\eta}\min_{a\in A}\sqrt{m_{a}})^{-1},
\end{align*}
where $\bar{C}_1$ is a positive constant depending on $\{c_{1a},c_{2a},c_{3a}\}_{a\in A}$ and $M$. Letting $\bm{x}$ be an optimal solution of \ref{eq_sp_fair} and redefining $\hat{\eta}:=6\hat{\eta}$, the conclusion follows.
\QEDA

\subsection{Proof of \Cref{thm_exact_gelbrich_v2}}\label{proof_thm_exact_gelbrich_v2}

\thmexactgelbrichvtwo*

\noindent\textit{Proof. }
Note that for any $a<\bar{a}\in A$, and an optimal comonotonic joint distribution $\Qe_{a,\bar{a}}$ of $f(\bfx,\tilde{\bm\xi}_{a})$ and $f(\bfx,\tilde{\bm\xi}_{\bar{a}})$ with marginals $\Pr_a,\Pr_{\bar a}$ such that $(f(\bfx,\tilde{\bm\xi}_{a})-\bm\mu_{a}^{\top}\bm{r}(\bfx)-s(\bfx), f(\bfx,\tilde{\bm\xi}_{\bar a})-\bm\mu_{\bar a}^{\top}\bm{r}(\bfx)-s(\bfx))\xrightarrow{m_{a}\rightarrow \infty,m_{\bar a}\rightarrow \infty}(\hat{\bm\xi}_{ a}^{\top}\bm{r}(\bfx), \hat{\bm\xi}_{\bar a}^{\top}\bm{r}(\bfx))$ and the random vectors $\hat{c}_a^{-1}\hat{\bm\xi}_{ a},\hat{c}_{\bar{a}}^{-1}\hat{\bm\xi}_{\bar a}$ obey the same distribution with zero mean and covariance matrix $\bm{\Sigma}_{a\bar a}$, we have
\begin{align*}
&\E_{\Qe_{a,\bar{a}}}[|f(\bfx,\tilde{\bm\xi}_{a})-f(\bfx,\tilde{\bm\xi}_{\bar{a}})|^2]
 = \E_{\Pr_{a}}\left[\left(f(\bfx,\tilde{\bm\xi}_{a})-\bm\mu_{a}^{\top}\bm{r}(\bfx)-s(\bfx)\right)^2\right]+\E_{\Pr_{\bar a}}\left[\left(f(\bfx,\tilde{\bm\xi}_{\bar{a}})-\bm\mu_{\bar a}^{\top}\bm{r}(\bfx)-s(\bfx)\right)^2\right] \\
 &-2\E_{\Pr_{a}}\left[\left(f(\bfx,\tilde{\bm\xi}_{a})-\bm\mu_{a}^{\top}\bm{r}(\bfx)-s(\bfx)\right)\left(\bm\mu_{a}^{\top}\bm{r}(\bfx) -\bm\mu_{\bar{a}}^{\top}\bm{r}(\bfx)\right)\right]\\
   &-2\E_{\Pr_{\bar a}}\left[\left(f(\bfx,\tilde{\bm\xi}_{\bar a})-\bm\mu_{\bar a}^{\top}\bm{r}(\bfx)-s(\bfx)\right)\left(\bm\mu_{a}^{\top}\bm{r}(\bfx) -\bm\mu_{\bar{a}}^{\top}\bm{r}(\bfx)\right)\right]\\
  &+ \left(\bm\mu_{a}^{\top}\bm{r}(\bfx) -\bm\mu_{\bar{a}}^{\top}\bm{r}(\bfx) \right)^2-2\E_{\Qe_{a,\bar{a}}}\left[\left(f(\bfx,\tilde{\bm\xi}_{a})-\bm\mu_{a}^{\top}\bm{r}(\bfx)-s(\bfx)\right)\left(f(\bfx,\tilde{\bm\xi}_{\bar{a}})-\bm\mu_{\bar a}^{\top}\bm{r}(\bfx)-s(\bfx)\right)\right]\\
  &\xrightarrow{m_{a}\rightarrow \infty,m_{\bar a}\rightarrow \infty}
\left(\bm\mu_{a}^{\top}\bm{r}(\bfx) -\bm\mu_{\bar{a}}^{\top}\bm{r}(\bfx) \right)^2 + \left(\hat{c}_{a}\sqrt{\bm{r}(\bfx)^{\top}\bm{\Sigma}_{a\bar a}\bm{r}(\bfx)} - \hat{c}_{\bar a}\sqrt{\bm{r}(\bfx)^{\top}\bm{\Sigma}_{a\bar a}\bm{r}(\bfx)}\right)^2.
\end{align*}
Thus, the claim follows. 
\QEDA

\newpage

\section{Not MICP-R Functions and Their Piecewise Linear Approximations}\label{sec_not_micpr}
When the utility functions are exponential or logarithmic, their corresponding $\F_q$ sets are typically not MICP-R. Hence we propose to approximate them using piecewise linear functions.
\begin{restatable}{proposition}{propnotmicprexp}\label{prop_not_micpr_exp}
If $f(\bfx, {\bm\xi})=\exp\left(\bm\xi^\top \bm{r}(\bfx)+ s(\bfx)\right)$, where $\bm{r}(\bfx)$ and $s(\bfx)$ are linear functions, then set $\F_q$, in general, is not MICP-R even when $m=2$ and $|A|=2$.
\end{restatable}
  \noindent\textit{Proof. }Let us consider a special case of \ref{eq_sp_fair} with $n=2$, $\mathcal X=[0,1]^2$, $m=2$, $|A|=2$, $m_a=m_{\bar{a}}=1$, $\epsilon =0.1$ and $f({\bfx,\bm\xi}_1)= \exp(x_1), f({\bfx,\bm\xi}_2)= \exp(x_2)$. 
Under this setting, we have
\[
\F_q=\left\{(\bfx,\nu)\in [0,1]^2\times \Re_+: |\exp(x_1)-\exp(x_2)|^q \leq \nu\right\}.
\]
It suffices to show that the sublevel set of the function $|\exp(x_1)-\exp(x_2)|^q$ is not MICP-R. Specifically, letting $\nu=0.1$ in $\F_q$, we consider the set
\[
\hat{\F}_q=\left\{(\bfx,\nu)\in [0,1]^2\times \Re_+: |\exp(x_1)-\exp(x_2)|^q \leq 0.1\right\}.
\]
Then for any two distinct points $\bm{x}_1,\bm{x}_2$ satisfying $x_{11}=\log(\exp(x_{21})+\sqrt[q]{0.1}),x_{12}=\log(\exp(x_{22})+\sqrt[q]{0.1})$, and $x_{21}, x_{22}\in (0,\log(e-\sqrt[q]{0.1}))$, one can show that their midpoint $(\bm{x}_1+\bm{x}_2)/2\notin \hat{\F}_q$. Indeed,
we have 
\begin{align*}
&\exp\left(\frac{1}{2}\left(\log(\exp(x_{21})+\sqrt[q]{0.1})+\log(\exp(x_{22})+\sqrt[q]{0.1})\right)\right)-\exp\left(\frac{1}{2}\left(x_{21}+x_{22}\right)\right)-\sqrt[q]{0.1}>0\\
(\Leftrightarrow) &\sqrt{\exp(x_{21})+\sqrt[q]{0.1}}\sqrt{\exp(x_{22})+\sqrt[q]{0.1}}-\sqrt{\exp(x_{21})}\sqrt{\exp(x_{22})}-\sqrt[q]{0.1}>0\\
(\Leftrightarrow) &\left(\exp(x_{21})+\sqrt[q]{0.1}\right)\left(\exp(x_{22})+\sqrt[q]{0.1}\right)-\left(\sqrt{\exp(x_{21})}\sqrt{\exp(x_{22})}+\sqrt[q]{0.1}\right)^2>0\\
(\Leftrightarrow) &\sqrt[q]{0.1}\exp(x_{21})+\sqrt[q]{0.1}\exp(x_{22})-2\sqrt[q]{0.1}\sqrt{\exp(x_{21})\exp(x_{22})}>0\\
(\Leftrightarrow) &\left(\sqrt{\exp(x_{21})}-\sqrt{\exp(x_{22})}\right)^2>0
\end{align*}
where the last inequality holds due to $x_{21}\neq x_{22}$.

Since there is an infinite number of these points, according to \Cref{lem_not-micpr}, the set $\hat{\F}_q$ is not MICP-R. Hence, the set $\F_q$ is not MICP-R. 
\QEDA

Therefore, when $f(\bfx, {\bm\xi})=\exp\left(\bm\xi^\top \bm{r}(\bfx)+ s(\bfx)\right)$, we propose to use an iterative discretization method to approximate it. 
Suppose that there are $T$ candidate points $\{\hat{\bm{x}}_\tau\}_{\tau\in [T]}$ and their corresponding $g_\tau= \bm\xi^\top \bm{r}(\hat{\bm{x}}_\tau)+ s(\hat{\bm{x}}_\tau)$, then we approximate $f(\bfx, {\bm\xi})\approx\max_{\tau\in [T]}\exp(g_\tau)+\exp(g_\tau) (\bm\xi^\top \bm{r}(\bfx)+ s(\bfx)-g_\tau) $. 
Thus, we obtain the following approximate MICP set of $X_i$ for each $i\in [m]$ as
\begin{equation}\nonumber
\begin{aligned}
X_i\approx \hat{X}_{i}&=\left\{(\bfx, \bar{w}_{i})\in \mathcal X\times\Re: \begin{aligned}
&\bar{w}_{i}\geq \exp(g_{\tau})+\exp(g_{\tau}) (\bm\xi^\top \bm{r}(\bfx)+ s(\bfx)-g_{\tau}),\forall \tau \in [T], \\
&\bar{w}_{i}\leq \exp(g_{\tau})+\exp(g_{\tau}) (\bm\xi^\top \bm{r}(\bfx)+ s(\bfx)-g_{\tau i})+M_{i\tau}(1-z_{i\tau}),\forall \tau \in [T],\\
&\sum_{\tau\in [T]} z_{i\tau} = 1, z_{i\tau}\in\{0,1\}, \forall \tau\in [T]
\end{aligned}\right\},
\end{aligned}
\end{equation}
where $M_{i\tau}=\max_{\bfx\in\mathcal X}\{\exp\left(\bm\xi_i^\top \bm{r}(\bfx)+ s(\bfx)\right)-\exp(g_{\tau})- \exp(g_{\tau}) (\bm\xi^\top \bm{r}(\bfx)+ s(\bfx)-g_{\tau i})\}$ for each $\tau\in [T]$.

A similar negative result holds when $f(\bfx, {\bm\xi})=\log\left(\bm\xi^\top \bm{r}(\bfx)+ s(\bfx)\right)$.
\begin{restatable}{proposition}{propnotmicprlog}\label{prop_not_micpr_log}
If $f(\bfx, {\bm\xi})=\log\left(\bm\xi^\top \bm{r}(\bfx)+ s(\bfx)\right)$, where $\bm{r}(\bfx)$ and $s(\bfx)$ are linear functions, and $\bm\xi^\top \bm{r}(\bfx)+ s(\bfx)>0$ for all $\bm{x}\in \mathcal X$, then set $\F_q$ is not MICP-R even when $m=2$ and $|A|=2$.
\end{restatable}
 \noindent\textit{Proof. }Let us consider a special case of \ref{eq_sp_fair} with $n=2$, $\mathcal X=[1,+\infty)^2$, $m=2$, $A=|2|$, $m_a=m_{\bar{a}}=1$ and $f({\bfx,\bm\xi}_1)= \log(x_1)$, $f({\bfx,\bm\xi}_2)= \log(x_2)$.  Under this setting, we have
\[
\F_q=\left\{(\bfx,\nu)\in [1,+\infty)^2\times\Re_+: |\log(x_1)-\log(x_2)|^q \leq \nu\right\}.
\]
Consider any two distinct points $(\bm{x}_1,\nu_1),(\bm{x}_2,\nu_2)\in [1,+\infty)^2\times \Re_+$ such that $x_{11}=x_{12}=1$,  $x_{21}=\exp(\sqrt[q]{\nu_1})>\exp(q-1)$, $x_{22}=\exp(\sqrt[q]{\nu_2})>\exp(q-1)$. Note that when $\nu\geq (q-1)^q$, the function $\exp(\sqrt[q]{\nu})$ is convex in $\nu$.
It remains to show that their midpoint $((\bm{x}_1,\nu_1)+(\bm{x}_2,\nu_2))/2\notin \F_q$. Indeed, we have
\begin{align*}
&\frac{1}{2}\left(\exp(\sqrt[q]{\nu_1})+\exp(\sqrt[q]{\nu_2})\right)-\exp\left(\sqrt[q]{\left(\frac{1}{2}\left(\nu_{1}+\nu_{2}\right)\right)}\right)>0
\end{align*}
due to the convexity of the function $\exp(\sqrt[q]{\nu})$.

Since there are infinitely many of these points, according to \Cref{lem_not-micpr}, the set $\F_q$ is not MICP-R.
\QEDA

Therefore, when $f(\bfx, {\bm\xi})=\log\left(\bm\xi^\top \bm{r}(\bfx)+ s(\bfx)\right)$, we also propose to use an iterative discretization method to approximate it. 
Suppose that there are $T$ candidate points $\{\hat{\bm{x}}_\tau\}_{\tau\in [T]}$ and their corresponding $g_\tau= \hat{\bm{x}}_\tau^\top \bm{r}(\hat{\bm{x}}_\tau)+ s(\hat{\bm{x}}_\tau)$, then we approximate $f(\bfx, {\bm\xi})\approx\min_{\tau\in [T]}\log(g_\tau)+g_\tau^{-1} (\bm\xi^\top \bm{r}(\bfx)+ s(\bfx)-g_\tau) $. Thus, we obtain the following approximate MICP set of $X_i$ for each $i\in [m]$ as
\begin{equation}\nonumber
\begin{aligned}
X_i\approx \hat{X}_{i}&=\left\{(\bfx, \bar{w}_{i})\in \mathcal X\times\Re: \begin{aligned}
&\bar{w}_{i}\leq \log(g_{\tau})+g_{\tau}^{-1} (\bm\xi^\top \bm{r}(\bfx)+ s(\bfx)-g_{\tau}),\forall \tau \in [T], \\
&\bar{w}_{i}\geq \log(g_{\tau})+g_{\tau}^{-1}  (\bm\xi^\top \bm{r}(\bfx)+ s(\bfx)-g_{\tau i})-M_{i\tau}(1-z_{i\tau}),\forall \tau \in [T],\\
&\sum_{\tau\in [T]} z_{i\tau} = 1, z_{i\tau}\in\{0,1\}, \forall \tau\in [T]
\end{aligned}\right\},
\end{aligned}
\end{equation}
where $M_{i\tau}=\max_{\bfx\in\mathcal X}\{-\log\left(\bm\xi_i^\top \bm{r}(\bfx)+ s(\bfx)\right)+\log(g_{\tau})+g_{\tau}^{-1} (\bm\xi^\top \bm{r}(\bfx)+ s(\bfx)-g_{\tau i})\}$ for each $\tau\in [T]$.

\newpage

\section{Additional Numerical Results}
\subsection{Results of Exact MICP Formulations With Jensen Inequality}
\Cref{result_m0_m1_m2_jensen} and \Cref{result_m3_m4_jensen} display the numerical results for the \hyperref[set_F_model1]{Vanilla Formulation}, \hyperref[WD_q_trans_micpr]{Discretized Formulation}, \hyperref[model2_WD1]{Complementary Formulation}, \hyperref[model3_WD]{Quantile Formulation}, \hyperref[model4_WD]{Aggregate Quantile Formulation} by adding Jensen bound as a valid inequality. We see that adding the Jensen bound improves the lower bounds for most instances. However, it often yields worse upper bounds and does not help solve the hard instances to optimality. Therefore, the main paper only reports the exact formulation comparison results without having Jensen bound. 
\begin{table}[htbp]
\centering
\caption{Results of Exact MICP Formulations (With Jensen Inequality)}\label{result_m0_m1_m2_jensen}
\tiny
\begin{tabular}{c|rrrr|rrrr|rrrr}
\hline
\multirow{2}{*}{m} & \multicolumn{4}{c|}{\hyperref[set_F_model1]{Vanilla Formulation}} & \multicolumn{4}{c|}{\hyperref[WD_q_trans_micpr]{Discretized Formulation}} & \multicolumn{4}{c}{\hyperref[model2_WD1]{Complementary Formulation}} \\ \cline{2-13}
& \multicolumn{1}{c}{Obj.Val} & \multicolumn{1}{c}{LB} & \multicolumn{1}{c}{Gap (\%)} & \multicolumn{1}{c|}{Time} & \multicolumn{1}{c}{Obj.Val} & \multicolumn{1}{c}{LB} & \multicolumn{1}{c}{Gap (\%)} & \multicolumn{1}{c|}{Time} & \multicolumn{1}{c}{Obj.Val} & \multicolumn{1}{c}{LB} & \multicolumn{1}{c}{Gap (\%)} & \multicolumn{1}{c}{Time} \\ \hline
15  & 342.44 & 207.12 & 39.51 & 3600.00 & 342.43 & 342.43 & 0.00  & 57.36   & 345.38  & 207.12 & 40.03 & 3600.00 \\
20  & 232.67 & 89.13  & 61.69 & 3600.00 & 236.27 & 192.01 & 18.73 & 3600.00 & 245.28  & 89.13  & 63.66 & 3600.00 \\
25  & 139.71 & 46.03  & 67.05 & 3600.00 & 158.40 & 75.62  & 52.26 & 3600.00 & 248.53  & 46.03  & 81.48 & 3600.00 \\
30  & 177.52 & 109.56 & 38.28 & 3600.00 & 218.84 & 109.56 & 49.94 & 3600.00 & 217.74  & 109.56 & 49.68 & 3600.00 \\
35  & 140.54 & 92.95  & 33.86 & 3600.00 & 326.31 & 92.95  & 71.51 & 3600.00 & 963.45  & 92.95  & 90.35 & 3600.00 \\
40  & 256.82 & 187.66 & 26.93 & 3600.00 & 583.48 & 187.66 & 67.84 & 3600.00 & 294.49  & 187.66 & 36.27 & 3600.00 \\
45  & 223.49 & 138.17 & 38.18 & 3600.00 & 480.48 & 138.17 & 71.24 & 3600.00 & 693.37  & 138.17 & 80.07 & 3600.00 \\
50  & 170.04 & 116.18 & 31.67 & 3600.00 & 933.76 & 116.18 & 87.56 & 3600.00 & 915.49  & 116.18 & 87.31 & 3600.00 \\
55  & 209.61 & 141.53 & 32.48 & 3600.00 & 576.83 & 141.53 & 75.46 & 3600.00 & 636.08  & 141.53 & 77.75 & 3600.00 \\
60  & 131.32 & 69.87  & 46.80 & 3600.00 & 546.72 & 69.87  & 87.22 & 3600.00 & 857.62  & 69.87  & 91.85 & 3600.00 \\
65  & 153.22 & 83.21  & 45.70 & 3600.00 & 548.73 & 83.21  & 84.84 & 3600.00 & 1125.53 & 83.21  & 92.61 & 3600.00 \\
70  & 136.70 & 83.84  & 38.67 & 3600.00 & 754.04 & 83.84  & 88.88 & 3600.00 & 962.53  & 83.84  & 91.29 & 3600.00 \\
75  & 173.96 & 85.22  & 51.01 & 3600.00 & 703.20 & 85.22  & 87.88 & 3600.00 & 881.76  & 85.22  & 90.34 & 3600.00 \\
80  & 159.68 & 112.03 & 29.84 & 3600.00 & 567.92 & 112.03 & 80.27 & 3600.00 & 746.74  & 112.03 & 85.00 & 3600.00 \\
85  & 194.68 & 104.78 & 46.18 & 3600.00 & 599.51 & 104.78 & 82.52 & 3600.00 & 563.55  & 104.78 & 81.41 & 3600.00 \\
90  & 176.64 & 126.73 & 28.26 & 3600.00 & 807.08 & 126.73 & 84.30 & 3600.00 & 1713.39 & 126.73 & 92.60 & 3600.00 \\
95  & 231.74 & 132.55 & 42.80 & 3600.00 & 749.70 & 132.55 & 82.32 & 3600.00 & 1144.84 & 132.55 & 88.42 & 3600.00 \\
100 & 147.47 & 91.76  & 37.78 & 3600.00 & 683.09 & 91.76  & 86.57 & 3600.00 & 1427.90 & 91.76  & 93.57 & 3600.00 \\ \hline
\end{tabular}
\end{table}

\begin{table}[htbp]
\centering
\caption{Results of Exact MICP Formulations (With Jensen Inequality)}\label{result_m3_m4_jensen}
\tiny
\begin{tabular}{c|rrrr|rrrr}
\hline
\multirow{2}{*}{m} & \multicolumn{4}{c|}{\hyperref[model3_WD]{Quantile Formulation}} & \multicolumn{4}{c}{\hyperref[model4_WD]{Aggregate Quantile Formulation}} \\ \cline{2-9}
& \multicolumn{1}{c}{Obj.Val} & \multicolumn{1}{c}{LB} & \multicolumn{1}{c}{Gap (\%)} & \multicolumn{1}{c|}{Time} & \multicolumn{1}{c}{Obj.Val} & \multicolumn{1}{c}{LB} & \multicolumn{1}{c}{Gap (\%)} & \multicolumn{1}{c}{Time} \\ \hline
15  & 342.43 & 342.43 & 0.00  & 0.56    & 342.43 & 342.40 & 0.01 & 0.54    \\
20  & 230.62 & 230.62 & 0.00  & 5.26    & 230.62 & 230.61 & 0.01 & 0.49    \\
25  & 135.03 & 135.03 & 0.00  & 18.58   & 135.03 & 135.03 & 0.00 & 1.65    \\
30  & 172.21 & 172.21 & 0.00  & 53.27   & 172.21 & 172.20 & 0.00 & 5.92    \\
35  & 133.54 & 133.54 & 0.00  & 911.07  & 133.54 & 133.54 & 0.00 & 8.81    \\
40  & 252.42 & 252.42 & 0.00  & 3068.64 & 252.42 & 252.42 & 0.00 & 15.39   \\
45  & 219.63 & 192.58 & 12.32 & 3600.00 & 219.17 & 219.17 & 0.00 & 21.90   \\
50  & 170.01 & 120.79 & 28.95 & 3600.00 & 169.99 & 169.99 & 0.00 & 41.36   \\
55  & ---    & 143.47 & 29.93 & 3600.00 & 204.77 & 204.76 & 0.00 & 75.79   \\
60  & ---    & 79.65  & 39.12 & 3600.00 & 130.84 & 130.84 & 0.00 & 199.56  \\
65  & ---    & 83.25  & 41.59 & 3600.00 & 142.55 & 142.54 & 0.00 & 327.37  \\
70  & ---    & 83.94  & 38.24 & 3600.00 & 136.33 & 134.50 & 1.34 & 3600.00 \\
75  & ---    & 85.22  & 38.21 & 3600.00 & ---    & 135.69 & 1.61 & 3600.00 \\
80  & ---    & 112.03 & 29.35 & 3600.00 & ---    & 155.35 & 2.03 & 3600.00 \\
85  & ---    & 104.78 & 28.12 & 3600.00 & ---    & 144.59 & 0.81 & 3600.00 \\
90  & ---    & 126.73 & 26.15 & 3600.00 & ---    & 169.96 & 0.95 & 3600.00 \\
95  & ---    & 132.55 & 23.74 & 3600.00 & ---    & 171.93 & 1.09 & 3600.00 \\
100 & ---    & 91.76  & 35.32 & 3600.00 & ---    & 139.46 & 1.68 & 3600.00 \\ \hline
\end{tabular}
\end{table}

\subsection{Fair Knapsack}\label{sec_fair_knapsack}
This subsection extends \ref{eq_sp_fair} to the classic knapsack problem. Given weights $\bm{w}\in\Ze_+^{m}$ and values $\bm{\xi}\in\Ze_+^{m}$ of a set of $m$ items, the objective of the knapsack problem is to select a subset of items to maximize the total value given that the total weight does not exceed the capacity $C$. That is, the knapsack problem can be formulated as
\begin{align}\label{model_knapsack}
V^*=\max_{\bfx} \left\{ \sum_{i\in[m]} \xi_{i} x_i: \sum_{i\in[m]} w_{i} x_i \leq C,x_i\in\{0,1\}, \forall i\in[m]\right\},
\end{align}
where the binary variable $x_i$ indicates whether item $i\in [m]$ is selected or not. 
We define the distributional fairness for the knapsack problem, where we choose the utility function of the fair knapsack problem to be the value of item $f(\bfx,\xi_i) = \xi_i x_i$ if being selected for each $i\in [m]$. 

In this experiment, we use $A=\{a,\bar{a}\}$ and generate the hypothetical data in the following manner. The weights $w_{i}$ are drawn from $\text{Unif}\{1, 100\}$. The first $\lceil m/2\rceil$ data points are assigned with the sensitive attribute $a$, where their values $\{\xi_{i}\}_{i\in [\lceil m/2\rceil]}$ are drawn independently from $\{\text{Unif}\{w_{i}+10, w_{i}+30\}\}_{i\in [\lceil m/2\rceil]}$. The remaining data points are assigned with $\bar{a}$, where their values $\{\xi_{i}\}_{i\in [\lceil m/2\rceil+1,m]}$ are drawn independently from $\{\text{Unif}\{w_{i}+20, w_{i}+60\}\}_{_{i\in [\lceil m/2\rceil+1,m]}}$. We generate a dataset of $m=1{,}000$ items and choose the capacity $C = 0.5\sum_{i\in[m]}w_{i}$. We set the inefficiency level parameter to $\epsilon\in\{0.01,0.05,0.1\}$. We choose type $q=2$ Wasserstein fairness and solve \ref{eq_sp_fair} using its AM algorithm in Section \ref{sec_AM}. For ease of illustration, we only display the histograms of each group's utility for selected items $(x=1)$. The probabilities of selection $\Pr_{a}(x=1)$ and $\Pr_{\bar{a}}(x=1)$ are reported.

\begin{figure}[htbp]
\begin{subfigure}[Vanilla]{{
\centering\includegraphics[width=0.23\textwidth]{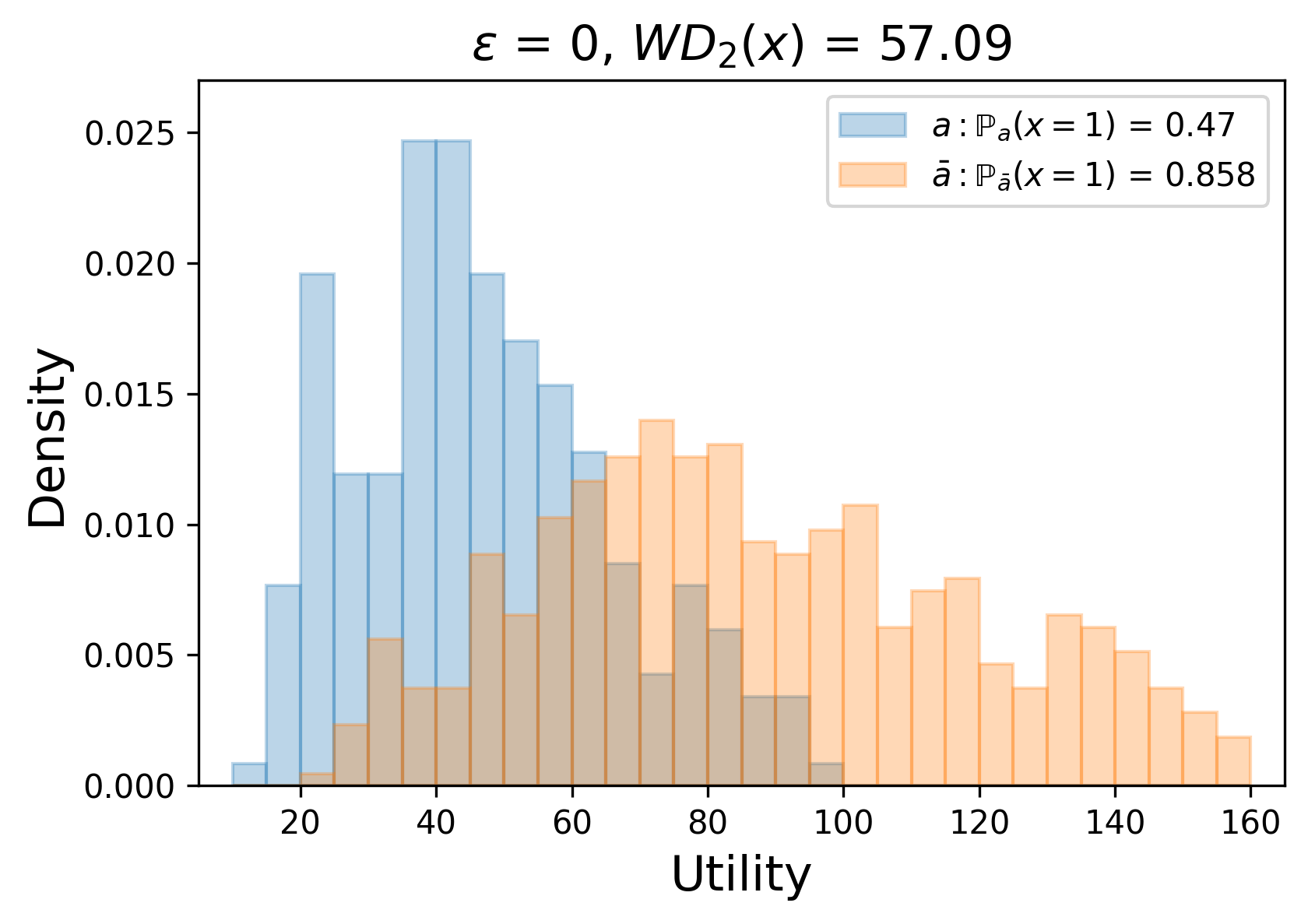}}\label{fig_knapsack_eps0}}
\end{subfigure}
\hfill
\begin{subfigure}[DFSO with $\epsilon=0.01$]{{
\centering\includegraphics[width=0.23\textwidth]{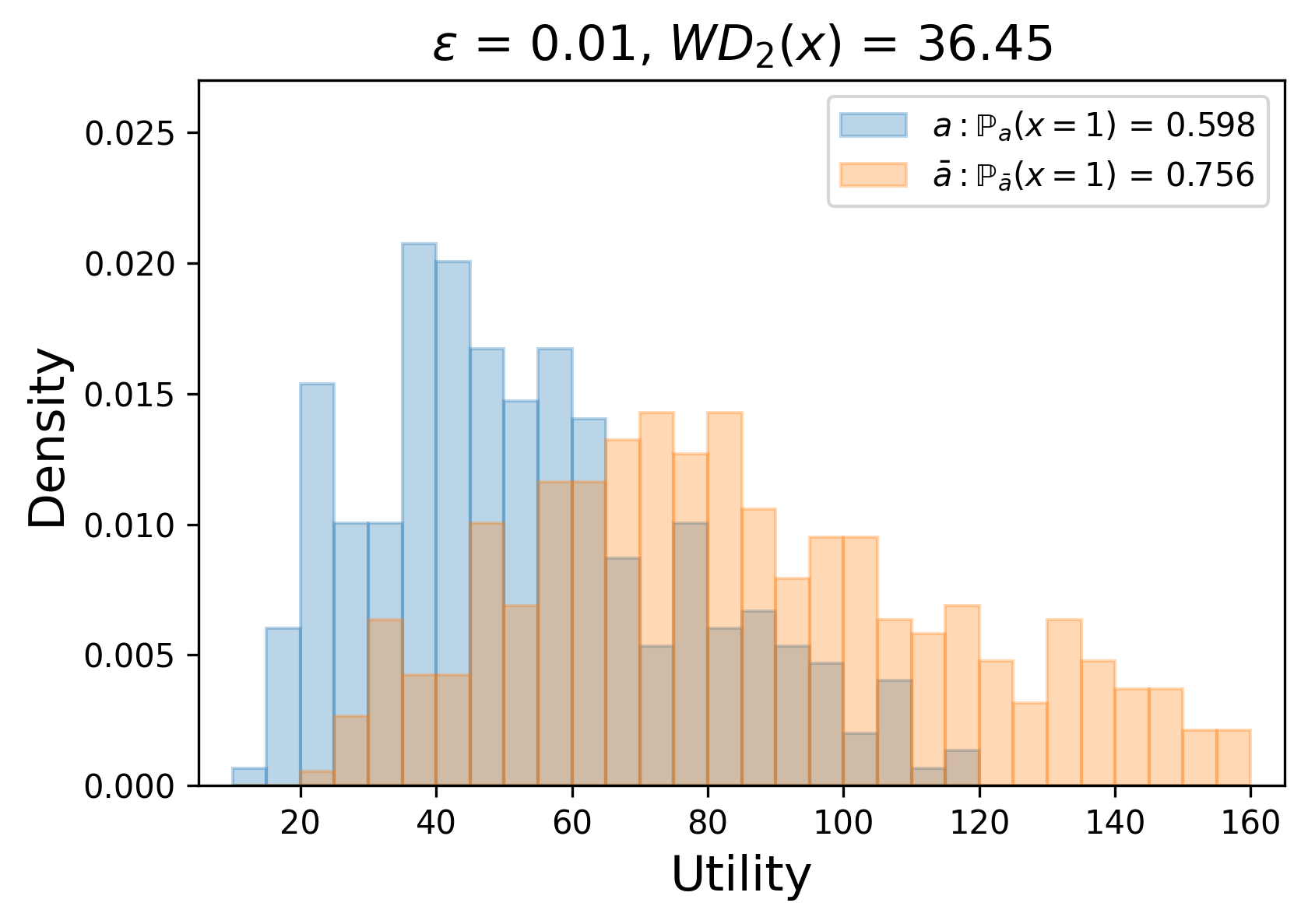}}\label{fig_knapsack_eps3}}
\end{subfigure}
\hfill
\begin{subfigure}[DFSO with $\epsilon=0.05$]{{
\centering\includegraphics[width=0.23\textwidth]{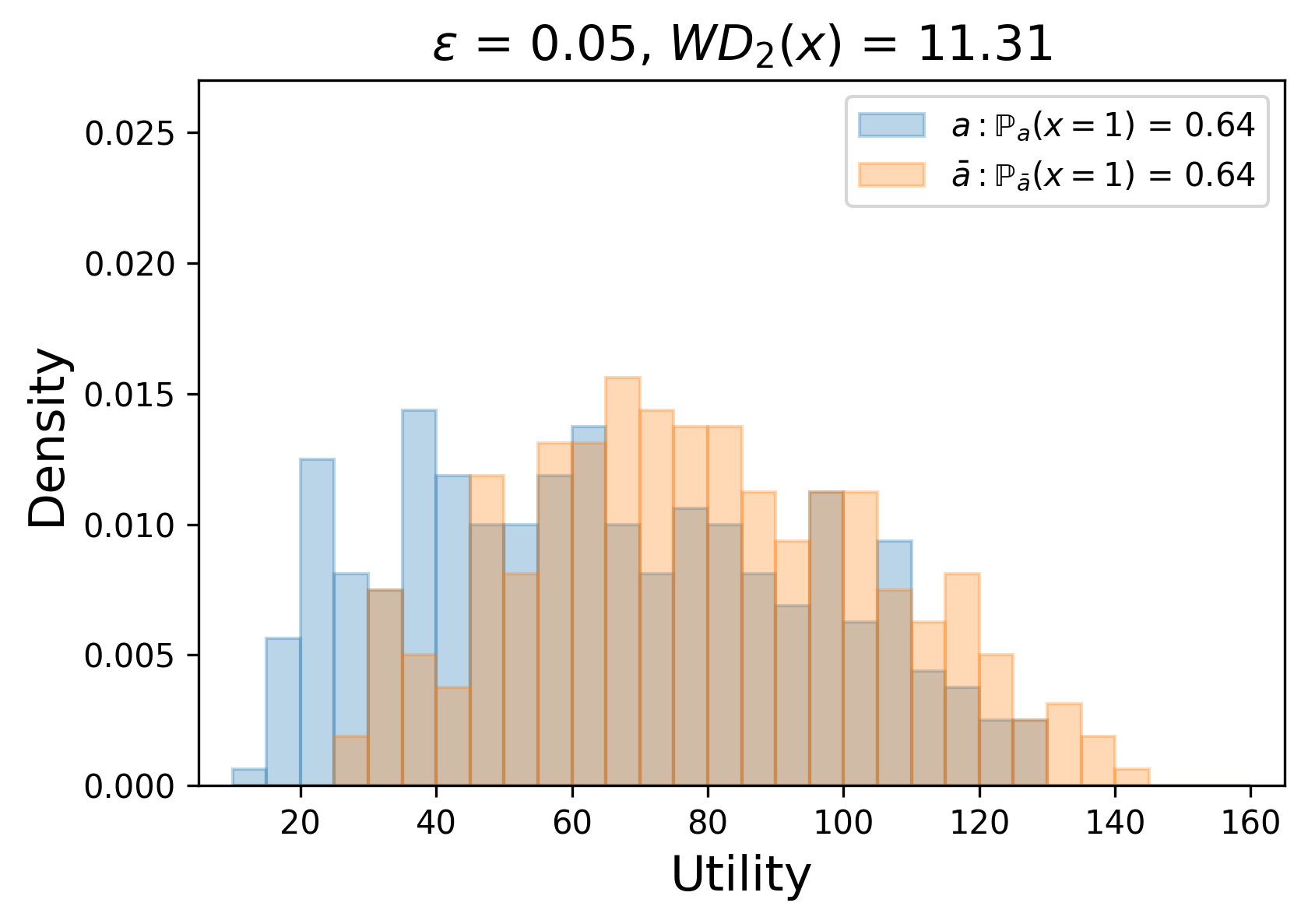}}\label{fig_knapsack_eps5}}
\end{subfigure}
\hfill
\begin{subfigure}[DFSO with $\epsilon=0.1$]{{
\centering\includegraphics[width=0.24\textwidth]{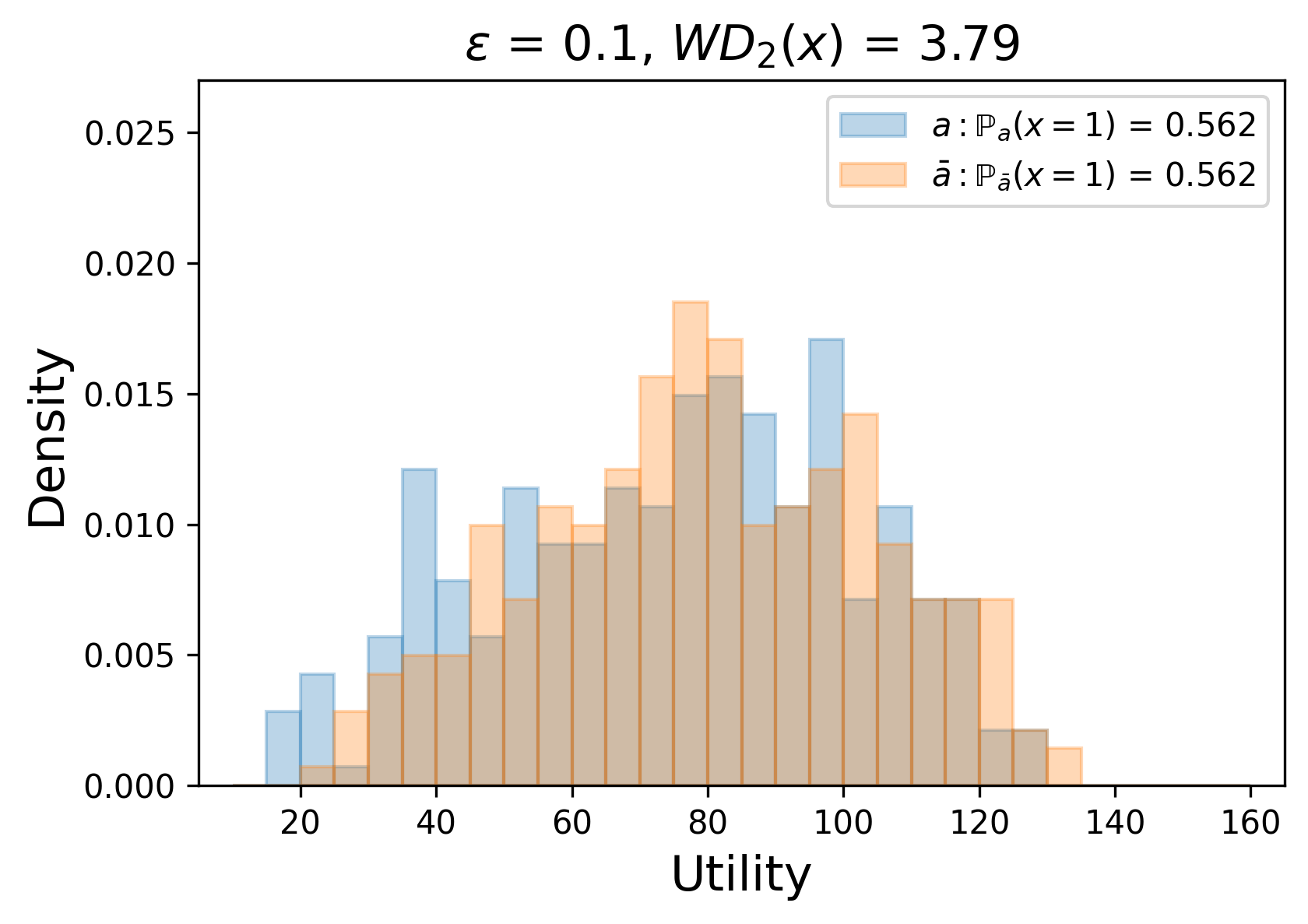}}\label{fig_knapsack_eps10}}
\end{subfigure}
\caption{Histograms of Utility for Fair Knapsack}\label{fig_knapsack_dfso}
\end{figure}

Figure \ref{fig_knapsack_dfso} presents the histograms of utility for fair knapsack. The vanilla knapsack problem \eqref{model_knapsack} has a Wasserstein fairness score of $57.09$, where it selects $47\%$ and $85.8\%$ of items from groups $a$ and $\bar{a}$, respectively. In Figures \ref{fig_knapsack_eps3}-\ref{fig_knapsack_eps10}, \ref{eq_sp_fair} improves the Wasserstein fairness score significantly between groups $a$ and $\bar{a}$ by slightly reducing the efficiency. The two groups have the same probability of selected items when $\epsilon\in\{0.05,0.1\}$. 
This shows that the proposed approach can achieve distributional fairness for the knapsack problem while maintaining high efficiency.

\color{black}

\end{appendices}

\end{document}